\documentclass[10pt,francais]{format}
\usepackage[francais,english]{babel}
\usepackage{amssymb,rotating}
\usepackage{pstricks,pst-plot}
\usepackage{xypic,amsfonts,amsmath}

\title[Indice et d{\'e}composition de Cartan]{Indice et d{\'e}composition de Cartan d'une alg{\`e}bre de Lie
  semi-simple r{\'e}elle}
\alttitle{The index and the Cartan decomposition of a real semisimple
  Lie algebra}
\author[Anne Moreau]{Anne Moreau}

\address{Universit{\'e} Paris 7\\
Institut de Math{\'e}matiques de Jussieu \\
Th{\'e}orie des groupes \\
Case 7012 \\ 2 Place Jussieu \\
75251 Paris Cedex 05, France }
\email{moreaua@math.jussieu.fr}



\altkeywords{index, Iwasawa decomposition,  Lie algebra,
  quasi-reductive, Cartan involution, stable form, Cayley transform}

\def\g{\mathfrak{g}}
\def\h{\mathfrak{h}}
\def\r{\mathfrak{r}}
\def\z{\mathfrak{z}}
\def\p{\mathfrak{p}}
\def\q{\mathfrak{q}}
\def\m{\mathfrak{m}}
\def\u{\mathfrak{u}}
\def\t{\mathfrak{t}}
\def\n{\mathfrak{n}}
\def\a{\mathfrak{a}}
\def\b{\mathfrak{b}}

\def\qed{\hfill \rule{.2cm}{.2cm} \\}

\newtheorem{theoreme}{Th{\'e}or{\`e}me}[section]

\newtheorem{lemme}[theoreme]{Lemme}
\newtheorem{proposition}[theoreme]{Proposition}

\newtheorem{definition}[theoreme]{D{\'e}finition\rm}
\newtheorem{remarque}{\it Remarque}

\renewcommand{\theequation}{\arabic{equation}}
\begin{document}

\newtheorem{theo_alpha}{Th{\'e}or{\`e}me}
\def\theo_alpha{\Alph{theo_alpha}}

\def\dv#1#2{\langle {#1},{#2}\rangle}
\def\an#1#2{\def\deux{#2} \ifx\deux\empty {\mathcal O}_{#1} 
\else {\mathcal O}_{#1,#2} \fi }
\def\tk#1#2{{#2}\otimes _{#1}}
\def\tens{\raisebox{.3mm}{\scriptsize$\otimes$}}

\def\N{\mathbb{N}}
\def\Z{\mathbb{Z}}
\def\R{\mathbb{R}}
\def\C{\mathbb{C}}
\def\Q{\mathbb{Q}}
\def\A{\mathbb{A}}

\def\g{\mathfrak{g}}
\def\h{\mathfrak{h}}
\def\r{\mathfrak{r}}
\def\z{\mathfrak{z}}
\def\p{\mathfrak{p}}
\def\q{\mathfrak{q}}
\def\m{\mathfrak{m}}
\def\s{\mathfrak{sl}}
\def\u{\mathfrak{u}}
\def\t{\mathfrak{t}}
\def\n{\mathfrak{n}}
\def\a{\mathfrak{a}}
\def\b{\mathfrak{b}}

\begin{abstract}
La d{\'e}composition d'Iwasawa $\mathfrak{g}_0=\mathfrak{k}_0 \oplus \widehat{\mathfrak{a}}_0 \oplus
\mathfrak{n}_0$ issue de la d{\'e}composition de Cartan $\mathfrak{g}_0=\mathfrak{k}_0 \oplus
\mathfrak{p}_0$ d'une alg{\`e}bre de Lie semisimple r{\'e}elle permet d'{\'e}crire $\mathfrak{g}_0$
sous la forme $\mathfrak{g}_0=\mathfrak{k}_0 \oplus \mathfrak{b}_0$, avec $\mathfrak{b}_0=\widehat{\mathfrak{a}}_0 \oplus
\mathfrak{n}_0$. La question de savoir si l'indice est additif dans la
d{\'e}composition $\mathfrak{g}_0=\mathfrak{k}_0 \oplus \mathfrak{b}_0$ a
{\'e}t{\'e} soulev{\'e}e par M. Ra{\"\i}s dans \cite{Rais}. Dans \cite{Moreau3}, il est
{\'e}crit que l'indice est toujours additif pour cette
d{\'e}composition. Pr{\'e}cis{\'e}ment, j'y affirme que
l'indice de $\mathfrak{b}$ est donn{\'e} par la relation\,:
${\rm ind} \; \mathfrak{b} = {\rm rg \ } \mathfrak{g} - {\rm rg \ } \mathfrak{k}$,  o{\`u} $\mathfrak{b}$ est
le complexifi{\'e} de $\mathfrak{b}_0$. Ce r{\'e}sultat n'est en fait pas vrai
en g{\'e}n{\'e}ral. On dispose en fait de l'in{\'e}galit{\'e}\,: ${\rm ind} \;
\mathfrak{b} \geq {\rm rg \ } \mathfrak{g} - {\rm rg \ }
\mathfrak{k}$, et l'{\'e}galit{\'e} a lieu si, et seulement si, une certaine
condition est satisfaite. Cet article a pour but de corriger cette
erreur. On reprend la d{\'e}marche de \cite{Moreau3} pour obtenir
cette fois l'in{\'e}galit{\'e} pr{\'e}c{\'e}dente. 
On donne alors une carat{\'e}risation des
alg{\`e}bres de Lie simples r{\'e}elles $\mathfrak{g}_0$ pour lesquelles
l'{\'e}galit{\'e} a lieu. On {\'e}tudie ou outre dans cet article le caract{\`e}re quasi-r{\'e}ductif de certaines sous-alg{\`e}bres de
$\mathfrak{g}$. Cette partie est nouvelle par rapport {\`a}  \cite{Moreau3}.
\end{abstract}
\begin{altabstract}
The Iwasawa decomposition $\mathfrak{g}=\mathfrak{k}_0 \oplus \widehat{\mathfrak{a}}_0 \oplus \mathfrak{n}_0$ of the
real semisimple Lie algebra $\mathfrak{g}_0$ comes from its Cartan
decomposition $\mathfrak{g}_0=\mathfrak{k}_0 \oplus \mathfrak{p}_0$. Then we get
$\mathfrak{g}_0=\mathfrak{k}_0 \oplus \mathfrak{b}_0$ where
$\mathfrak{b}_0=\widehat{\mathfrak{a}}_0 \oplus \mathfrak{n}_0$. The
question of knowing if the index were additive in the decomposition
$\mathfrak{g}_0=\mathfrak{k}_0 \oplus \mathfrak{b}_0$ goes back
M. Ra{\"\i}s \cite{Rais}. In \cite{Moreau3}, I wrote that the index always
is additive for this decomposition. Precisly, I claim that the index
of $\mathfrak{b}$ is given by the following formula\,:
${\rm ind} \; \mathfrak{b} = {\rm rk \ } \mathfrak{g} - {\rm rk \ } \mathfrak{k}$, where $\mathfrak{b}$  is the
complexification of $\mathfrak{b}_0$. This result is false in
general. We actually have an inequality\,: ${\rm ind} \;
\mathfrak{b} \geq {\rm rg \ } \mathfrak{g} - {\rm rg \ }
\mathfrak{k}$. The goal of this paper is  to correct this
mistake. We resume the approach of \cite{Moreau3} to obtain this time the
previous inequality. Then we give in more a characterization of the semisimple real
Lie algebra $\mathfrak{g}_0$ for which the index is additive in  the decomposition
$\mathfrak{g}_0=\mathfrak{k}_0 \oplus \mathfrak{b}_0$. Moreover, we
study in this paper the quasi-reductive character of  some subalgebras of
$\mathfrak{g}$. This is a new part in comparison with  \cite{Moreau3}.
\end{altabstract} 

\maketitle

\section*{Introduction}

Soit $\mathfrak{g}_0$ une alg{\`e}bre de Lie semi-simple r{\'e}elle et soit
$\mathfrak{g}_0=\mathfrak{k}_0 \oplus \mathfrak{p}_0$ une d{\'e}composition de Cartan de
$\mathfrak{g}_0$. On note $\theta$ l'involution de Cartan correspondante. Soit
$\widehat{\mathfrak{a}}_0$ un sous-espace ab{\'e}lien maximal de
$\mathfrak{p}_0$. Le sous-espace $\widehat{\mathfrak{a}}_0$ est form{\'e}
d'{\'e}l{\'e}ments semi-simples dans $\mathfrak{g}_0$ et, pour
$\lambda$ dans $\widehat{\mathfrak{a}}_0^{*}$, on pose
$$\mathfrak{g}_0^{\lambda}=\{X \in \mathfrak{g}_0 \ | \
\lbrack H,X]=\lambda(H)X, \ \forall H \in
\widehat{\mathfrak{a}}_0\} .$$ 
L'ensemble $\Sigma$ constitu{\'e} des formes lineaires non nulles
$\lambda$ sur $\widehat{\mathfrak{a}}_0$ pour lesquelles le sous-espace $\mathfrak{g}_0^{\lambda}$ est non
nul est un syst{\`e}me de racines dans $\widehat{\mathfrak{a}}_0^{ *}$. Soit $\Sigma_{+}$ un syst{\`e}me de racines
positives de $\Sigma$. On pose  
$$\mathfrak{n}_0=\bigoplus\limits_{\lambda \in \Sigma^+}
\mathfrak{g}_0^{\lambda},$$ 
de sorte qu'on obtient la d{\'e}composition d'Iwasawa de $\mathfrak{g}_0$ suivante\,:
$$\mathfrak{g}_0=\mathfrak{k}_0 \oplus \widehat{\mathfrak{a}}_0 \oplus \mathfrak{n}_0 .$$ 
En posant $\mathfrak{b}_0= \widehat{\mathfrak{a}}_0\ \oplus \mathfrak{n}_0$, on obtient la d{\'e}composition\,: 
$$\mathfrak{g}_0=\mathfrak{k}_0 \oplus\mathfrak{b}_0 .$$

Dans tout ce qui suit, on note sans indice $0$ les complexifi{\'e}s des
alg{\`e}bres de Lie r{\'e}elles not{\'e}es, elles, avec un indice $0$. Ainsi,
$\mathfrak{g}=(\mathfrak{g}_0)^{\mathbb{C}}$, $\mathfrak{k}=(\mathfrak{k}_0)^{\mathbb{C}}$,
$\widehat{\mathfrak{a}}=(\widehat{\mathfrak{a}}_0)^{\mathbb{C}}$,
$\mathfrak{n}=(\mathfrak{n}_0)^{\mathbb{C}}$,
$\mathfrak{b}=(\mathfrak{b}_0)^{\mathbb{C}}$, etc. L'indice d'une alg{\`e}bre de Lie $\mathfrak{q}$, not{\'e} ${\rm ind} \; \mathfrak{q}$, est la dimension
minimale des stabilisateurs pour l'action coadjointe des {\'e}l{\'e}ments de $\mathfrak{q}^*$. Si $\mathfrak{q}$ est le complexifi{\'e} d'une
alg{\`e}bre de Lie r{\'e}elle $\mathfrak{q}_0$, alors on a\,: 
${\rm ind} \; \mathfrak{q}={\rm ind} \; \mathfrak{q}_0$.\\

La question de savoir si l'indice est additif dans la
d{\'e}composition $\mathfrak{g}_0=\mathfrak{k}_0 \oplus \mathfrak{b}_0$ a
{\'e}t{\'e} soulev{\'e}e par M. Ra{\"\i}s dans \cite{Rais}. Dans \cite{Moreau3}, il est
{\'e}crit que l'indice est toujours additif pour cette d{\'e}composition. Pr{\'e}cis{\'e}ment, j'affirme que
l'indice de $\mathfrak{b}$ est donn{\'e} par la relation\,:
${\rm ind} \; \mathfrak{b} = {\rm rg \ } \mathfrak{g} - {\rm rg \ } \mathfrak{k}$,  o{\`u} $\mathfrak{b}$ est
le complexifi{\'e} de $\mathfrak{b}_0$. Ce r{\'e}sultat n'est  pas vrai
en g{\'e}n{\'e}ral, comme le prouve le cas de l'alg{\`e}bre de Lie simple r{\'e}elle $\mathfrak{g}=\mathfrak{sl}(2p,1)$, pour $p
\geq 1$. En effet, si l'indice {\'e}tait additif pour cette alg{\`e}bre, la formule pr{\'e}c{\'e}dente donnerait  ${\rm ind} \;
\mathfrak{b}=0$, car  ici ${\rm rg \ } \mathfrak{g} = {\rm rg \ }
\mathfrak{k}$. Or, d'apr{\`e}s \cite{Carmona}, l'alg{\`e}bre $\mathfrak{b}$ ne
poss{\`e}de pas, dans ce cas, d'orbite ouverte dans son
dual, ce qui signifie\,: ${\rm ind} \;
\mathfrak{b} > 0$. \\

Le but de cet article
est de corriger l'erreur commise dans \cite{Moreau3}. On reprend exactement la d{\'e}marche de 
\cite{Moreau3} pour construire des formes lin{\'e}aires r{\'e}guli{\`e}res sur
$\mathfrak{b}$. Mais on obtient cette fois l'in{\'e}galit{\'e}\,:
\begin{equation}\label{relation2} 
{\rm ind} \; \mathfrak{g}_0 \leq   {\rm ind} \; \mathfrak{k}_0 + {\rm ind} \; \mathfrak{b}_0. 
\end{equation}
On donne alors ici une caract{\'e}risation des alg{\`e}bres de Lie simples
r{\'e}elles pour lesquelles il y a {\'e}galit{\'e} dans la relation pr{\'e}c{\'e}dente.\\

La premi{\`e}re partie concerne la structure
de $\mathfrak{g}_0$. On rappelle dans la deuxi{\`e}me partie la construction
\guillemotleft en cascade\guillemotright \ de Kostant. Cette
construction servira {\`a}  d{\'e}finir dans les parties suivantes 
une forme lin{\'e}aire r{\'e}guli{\`e}re sur $\mathfrak{b}$ ainsi qu'une forme stable sur une certaine sous-alg{\`e}bre de
Borel de $\mathfrak{g}$ contenant $\mathfrak{b}$. On utilise dans
la partie
\ref{part_cayley} les transformations de Cayley pour obtenir des
relations utiles pour la suite (Proposition \ref{calcul_k}). On
prouve dans la partie $4$ la relation  (\ref{relation2}) et on obtient
une condition n{\'e}cessaire et suffisante pour qu'il y ait {\'e}galit{\'e} dans
la relation (\ref{relation2}) (Th{\'e}or{\`e}me \ref{formule_indice}). La partie
\ref{part_quasi} concerne les alg{\`e}bres de Lie quasi-r{\'e}ductives. Parmi les alg{\`e}bres
de Lie r{\'e}elles semi-simples $\mathfrak{g}_0$, on caract{\'e}rise  celles pour lesquelles la sous-alg{\`e}bre
$\mathfrak{b}$ est quasi-r{\'e}ductive. On donne en outre une
description pr{\'e}cise de certaines sous-alg{\`e}bres paraboliques
quasi-r{\'e}ductives de $\mathfrak{g}$. La derni{\`e}re partie repose sur la classification des alg{\`e}bres de Lie
simples r{\'e}elles. On compare certaines propri{\'e}t{\'e}s portant sur
$\mathfrak{b}$\,: \guillemotleft l'indice est additif dans la
d{\'e}composition $\mathfrak{g}=\mathfrak{k} \oplus \mathfrak{b}$\guillemotright, \guillemotleft $\mathfrak{b}$ est
quasi-r{\'e}ductive\guillemotright, \guillemotleft $\mathfrak{b}$ poss{\`e}de une
orbite ouverte dans son dual\guillemotright, etc. Pour chaque type d'alg{\`e}bre de Lie simple r{\'e}elle, on pr{\'e}cise
quelles sont les propri{\'e}tes satisfaites par $\mathfrak{b}$.\\

Notons que la partie \ref{part_quasi} est nouvelle par rapport {\`a}
\cite{Moreau3}. Les r{\'e}sultats principaux de cette partie (Th{\'e}or{\`e}me
\ref{car_sea} et Table \ref{colored_root}) sont, pour une large part, 
ind{\'e}pendants du reste. La partie \ref{part_calcul} est, quant {\`a} elle,
largement modifi{\'e}e par rapport {\`a} la partie correspondante de
\cite{Moreau3} compte tenu de l'erreur annonc{\'e}e.
   
\section{Quelques pr{\'e}cisions sur la structure de $\mathfrak{g}_0$
}\label{part_structure}

On regroupe dans cette partie quelques r{\'e}sultats concernant
la structure de $\mathfrak{g}_0$. On trouve les preuves de ces r{\'e}sultats dans
\cite{Dixmier} et \cite{Knapp}. Soit $\mathfrak{h}_0$ une
sous-alg{\`e}bre de Cartan de $\mathfrak{g}_0$, stable par $\theta$. Puisque $\mathfrak{h}_0$ est stable, elle s'{\'e}crit sous la
forme
$$ \mathfrak{h}_0=\mathfrak{a}_0 \oplus \mathfrak{t}_0,$$
avec $\mathfrak{a}_0$ dans $\mathfrak{p}_0$ et $\mathfrak{t}_0$ dans
$\mathfrak{k}_0$. La sous-alg{\`e}bre $\mathfrak{h}$ est une sous-alg{\`e}bre de
Cartan de $\mathfrak{g}$ et on note $\Delta$ le syst{\`e}me de racines
associ{\'e} au couple $(\mathfrak{g},\mathfrak{h})$. Les {\'e}l{\'e}ments de $\Delta$ sont {\`a}
valeurs r{\'e}elles sur $\mathfrak{a}_{0} \oplus i \mathfrak{t}_{0}$. On choisit un syst{\`e}me de racines positives $\Delta_{+}$ dans  $\Delta$
 en prenant $\mathfrak{a}_{0}$ avant $i \mathfrak{t}_{0}$ pour former l'ordre
lexicographique sur $(\mathfrak{a}_{0} \oplus i \mathfrak{t}_{0})^*$. Ainsi,
pour $\alpha$ une racine de $\Delta$ non nulle sur $\mathfrak{a}_{0}$, la positivit{\'e} de
$\alpha$ ne d{\'e}pend que de sa restriction {\`a} $\mathfrak{a}_{0}$. On note $\Pi$ la
base de $\Delta_+$. \\

On note encore $\theta$ l'extension $\mathbb{C}$-lin{\'e}aire de $\theta$ {\`a}
$\mathfrak{g}$. La transpos{\'e}e de $\theta$ est {\'e}galement not{\'e}e $\theta$. On note
$\Delta^{'}$  (respectivement $\Delta^{''}$) l'ensemble des racines de $\Delta$ 
qui s'annulent sur $\mathfrak{a}$ (respectivement qui ne s'annulent pas sur
$\mathfrak{a}$). On a $\theta(\Delta)=\Delta$ et $\Delta^{'}$ est l'ensemble des {\'e}l{\'e}ments de $\Delta$ invariants par
$\theta$. On pose $\Delta_{+}^{''}=\Delta^{''} \cap \Delta_{+}$,
$\Delta_{-}^{''}=\Delta^{ ''} \cap (-\Delta_{+})$. On
a $\theta(\Delta_{+}^{''})=\Delta_{-}^{''}$ et
$(\Delta_{+}^{''}+\Delta_{+}^{''}) \cap \Delta \subset
\Delta_{+}^{''}$. 

Pour chaque {\'e}l{\'e}ment $\alpha$ de $\Delta$, on fixe $X_{\alpha}$ un {\'e}l{\'e}ment non nul
de $\mathfrak{g}^{\alpha}$ et on note $H_{\alpha}$ l'unique
{\'e}l{\'e}ment de $[\mathfrak{g}^{\alpha},\mathfrak{g}^{-\alpha}]$ tel que
$\alpha(H_{\alpha})=2$. Une racine est dite
{\it r{\'e}elle} si elle prend des valeurs r{\'e}elles sur $\mathfrak{h}_0$
(i.e. si elle s'annule sur $\mathfrak{t}_0$), {\it imaginaire} si elle prend des valeurs imaginaires sur $\mathfrak{h}_0$
(i.e. si elle s'annule sur $\mathfrak{a}_0$), et {\it complexe}
sinon. Le lemme suivant est connu et ne pr{\'e}sente pas de difficult{\'e}\,:\\

\begin{lemme}\label{struct1}
\begin{description}
\item[{\rm (i)}] Pour $\alpha$ dans $\Delta$, on a\,: $\theta X_{\alpha}
  \in \mathfrak{g}^{\theta \alpha}$,

\item[{\rm (ii)}] Pour $\alpha$ dans $\Delta$, on a\,: $\theta H_{\alpha}=H_{\theta \alpha}$,

\item[{\rm (iii)}] Soit $\alpha$ une racine de $\Delta$. On a les
  {\'e}quivalences suivantes\,:
\begin{eqnarray*}
& & \alpha \textrm{ est r{\'e}elle } \iff \theta \alpha=-\alpha \iff H_{\alpha}
  \in \mathfrak{a}, \\
& & \alpha \textrm{ est imaginaire } \iff \theta \alpha=\alpha \iff H_{\alpha}
  \in \mathfrak{t} .
\end{eqnarray*}

\end{description}
\end{lemme}
~\\

La {\it dimension compacte} est par d{\'e}finition la dimension, $\dim
\mathfrak{t}_0$, de l'intersection de $\mathfrak{h}_0$ avec $\mathfrak{k}_0$ et
la {\it dimension non-compacte} est par d{\'e}finition la dimension, $\dim
\mathfrak{a}_0$, de l'intersection de $\mathfrak{h}_0$ avec $\mathfrak{p}_0$. On dit que
$\mathfrak{h}_0$ est {\it maximalement compacte} si la dimension
compacte est la plus grande possible et on dit que
$\mathfrak{h}_0$ est {\it maximalement non-compacte} si la dimension
non-compacte est la plus grande possible.\\

\begin{remarque} Il se peut que $\mathfrak{h}_0$ soit {\`a} la fois maximalement
  compacte et maximalement non-compacte. C'est le cas si $\mathfrak{g}_0$ est
  l'alg{\`e}bre de Lie r{\'e}elle sous-jacente {\`a} une alg{\`e}bre de Lie simple
  complexe, ou si $\mathfrak{g}_0$ est isomorphe {\`a} $\mathfrak{sl}(n,\mathbb{H})$
  ou {\`a} l'alg{\`e}bre de Lie simple exceptionelle $E IV$, comme on peut le
  voir {\`a} l'aide de la table \ref{tableau} pr{\'e}sent{\'e}e {\`a} la fin de ce chapitre.\\  
\end{remarque}

Si $\alpha$ est une racine imaginaire de $\Delta$, alors $\theta
\alpha=\alpha$ donc $\mathfrak{g}^{\alpha}$ est stable par $\theta$, et on a
$\mathfrak{g}^{\alpha}=(\mathfrak{g}^{\alpha} \cap \mathfrak{k}) \oplus (\mathfrak{g}^{\alpha} \cap
\mathfrak{p})$. Puisque $\mathfrak{g}^{\alpha}$ est de dimension $1$, on a
$\mathfrak{g}^{\alpha} \subseteq \mathfrak{k}$ ou $\mathfrak{g}^{\alpha} \subseteq
\mathfrak{p}$. On dit que la racine imaginaire $\alpha$ est {\it
compacte} si $\mathfrak{g}^{\alpha} \subseteq \mathfrak{k}$ et {\it
non-compacte} si $\mathfrak{g}^{\alpha} \subseteq \mathfrak{p}$. Le r{\'e}sultat
suivant est d{\'e}montr{\'e} en
\cite{Knapp}, proposition 6.70.\\

\begin{lemme}\label{max_compacte} Soit $\mathfrak{h}_0$ une
  sous-alg{\`e}bre de Cartan de $\mathfrak{g}_0$ stable
  par $\theta$. Alors, il n'existe pas de racine imaginaire non-compacte si,
  et seulement si, 
$\mathfrak{h}_0$ est
  maximalement non-compact et, il n'existe pas de racine r{\'e}elle si,
  et seulement si, 
$\mathfrak{h}_0$ est
  maximalement compact.\\

\end{lemme}

Soit $\widehat{\mathfrak{t}}_{0}$ un sous-espace ab{\'e}lien maximal du centralisateur
$\mathfrak{m}_0=\mathfrak{z}_{\mathfrak{k}_0}(\widehat{\mathfrak{a}}_{0})$ de $\widehat{\mathfrak{a}}_{0}$ dans
$\mathfrak{k}_0$. Le sous-espace $\widehat{\mathfrak{t}}_{0}$ est
form{\'e} d'{\'e}l{\'e}ments semi-simples dans $\mathfrak{g}_0$ et $\widehat{\mathfrak{h}}_{0}=\widehat{\mathfrak{a}}_{0} \oplus \widehat{\mathfrak{t}}_{0}$ est une
sous-alg{\`e}bre de Cartan de $\mathfrak{g}$ stable par $\theta$ qui est
maximalement non-compacte. On surmonte d'un chapeau les
ensembles d{\'e}finis pr{\'e}c{\'e}demment relatifs {\`a}
$\widehat{\mathfrak{h}}$. On a\,:
$$(\mathfrak{g}_{0})^{\lambda}=\mathfrak{g}_0 \cap (\bigoplus\limits_{\alpha \in
  \widehat{\Delta}_{+}^{ ''} \atop \alpha_{|\widehat{\mathfrak{a}}}=\lambda} \mathfrak{g}^{\alpha})$$
et
$$\mathfrak{n}=(\mathfrak{n}_{0})^{\mathbb{C}}=\bigoplus\limits_{\alpha \in
\widehat{\Delta}_{+}^{ ''}} \mathfrak{g}^{\alpha} .$$
D'o{\`u}
$$\mathfrak{b}=\widehat{\mathfrak{a}} \oplus (\bigoplus\limits_{\alpha \in
   \widehat{\Delta}_{+}^{ ''}} \mathfrak{g}^{\alpha}) .$$

On termine par un lemme\,:\\

\begin{lemme}\label{struct2}
L'ensemble
  $\widehat{\Delta}^{'}_{|\widehat{\mathfrak{t}}}$ est le syst{\`e}me de racines
  associ{\'e} au couple $(\mathfrak{m},\widehat{\mathfrak{t}})$. Si $\alpha$ appartient {\`a} $\widehat{\Delta}^{ '}$, on
  a 
$\mathfrak{m}^{\alpha_{|\widehat{\mathfrak{t}}}}=\mathfrak{g}^{\alpha}$. Enfin,
  on a\,:
$$\mathfrak{m}=\widehat{\mathfrak{t}} \oplus (\bigoplus\limits_{\alpha \in
    {\widehat{\Delta}}^{'} } \mathfrak{g}^{\alpha}) .$$
~\\

\end{lemme}

\section{Construction \guillemotleft en
  cascade\guillemotright \ de Kostant}\label{part_kostant} 

\subsection{} Soit $\mathfrak{h}_0$ une sous-alg{\`e}bre de Cartan de $\mathfrak{g}_0$
stable par $\theta$. On reprend les notations de la partie
pr{\'e}c{\'e}dente. On utilise en outre les notations introduites dans \cite{TauvelYu} et \cite{TauvelYu2}. 
 Si $\lambda$ appartient {\`a} $\mathfrak{h}^{*}$, on {\'e}crit $\langle
\lambda, \alpha^{\vee} \rangle$ pour $\lambda(H_{\alpha})$. Pour
toute partie $S$ de $\Pi$, on note $\Delta^{S}$ le syst{\`e}me de
racines engendr{\'e} par $S$, et $\Delta_{+}^{S}$ le syst{\`e}me de racines
positives correspondant. Si $S$ est une partie connexe de $\Pi$, le
syst{\`e}me de racines $\Delta^{S}$ est irreductible et on note
$\varepsilon_{S}$ la plus grande racine de $\Delta_{+}^{S}$.\\

Supposons que $S$ est une partie connexe de $\Pi$. Les r{\'e}sultats qui
suivent vont {\^e}tre utilis{\'e}s {\`a} plusieurs reprises dans la suite; ils sont
d{\'e}montr{\'e}s dans \cite{Bourbaki} et rappel{\'e}s dans \cite{Tauvel}, \cite{TauvelYu} et \cite{TauvelYu2}. Pour toute racine
$\alpha$ de $\Delta_{+}^{S} \setminus\{\varepsilon_S\}$, on a\,: $\langle
\alpha, {\varepsilon_S}^{\vee} \rangle \in \{0,1\}$. Si $T$ est l'ensemble des
racines $\alpha$ de $\Delta^{S}$ qui v{\'e}rifient $\langle
\alpha, {\varepsilon_S}^{\vee} \rangle =0$, alors $T$ est un syst{\`e}me de
racines dans le sous-espace de $\mathfrak{h}^{*}$ qu'il engendre et l'ensemble
$\{ \alpha \in S \ | \ \langle
\alpha, {\varepsilon_S}^{\vee} \rangle =0 \}$ forme une base de $T$. De
plus, si $\alpha$ appartient {\`a} $T \cap \Delta_{+}^{S}$, alors on a\,:
$\alpha \pm \varepsilon_S \not\in \Delta$. Ainsi, pour $\alpha$ dans $T$, les racines $\alpha$ et $\varepsilon_S$ sont fortement orthogonales.\\  

On rappelle la construction et quelques propri{\'e}t{\'e}s d'un
ensemble de racines deux {\`a} deux fortement orthogonales dans
$\Delta$. Par r{\'e}currence sur le cardinal de $S$, on d{\'e}finit un
sous-ensemble $\mathcal{K}(S)$ de l'ensemble des parties de $\Pi$ de la
mani{\`e}re suivante\,: \begin{description}
\item[a)] $\mathcal{K}(\emptyset)=\emptyset$,
\item[b)] Si $S_1$,\ldots,$S_r$ sont les composantes connexes de $S$,
  on a\,:
$$\mathcal{K}(S)=\mathcal{K}(S_1) \cup \cdots \cup \mathcal{K}(S_r),$$
\item[c)] Si $S$ est connexe, alors\,:
$$\mathcal{K}(S)=\{S\} \cup \mathcal{K}(\{ \alpha \in S \ 
| \ \langle \alpha, \varepsilon_{S}^{\vee} \rangle =0 \}) .$$
\end{description}
~\\

Les deux lemmes suivants regroupent des propri{\'e}t{\'e}s de cette
construction utiles pour la suite. Ils sont {\'e}nonc{\'e}s dans
\cite{Tauvel}, \cite{TauvelYu} ou \cite{TauvelYu2}.\\

\begin{lemme}\label{kostant1} \begin{description}
\item[{\rm (i)}] Tout {\'e}l{\'e}ment $K$ de $\mathcal{K}(S)$ est une partie connexe de
  $\Pi$.
\item[{\rm (ii)}] Si $K,K'$ appartiennent {\`a} $\mathcal{K}(S)$, alors ou
  bien $K \subset K'$, ou bien $K' \subset K$, ou bien $K$ et $K'$ sont des parties disjointes de $S$
  telles que $\alpha + \beta$ n'appartient pas {\`a} $\Delta$, pour
  $\alpha$ dans $\Delta^{K}$ et $\beta$ dans $\Delta^{K'}$.
\item[{\rm (iii)}] Si $K$ et $K'$ sont des {\'e}l{\'e}ments distincts de
  $\mathcal{K}(S)$, alors $\varepsilon_{K}$ et $\varepsilon_{K'}$ sont
  fortement orthogonales.\\

\end{description}
\end{lemme}
 
Si $K \in \mathcal{K}(\Pi)$, on pose\,:

$$\Gamma^{K}=\{ \alpha \in \Delta^{K} \ | \ 
 \langle \alpha, \varepsilon_{K}^{\vee} \rangle > 0 \}, \hspace{0.5cm}
\Gamma_{0}^{K}=\Gamma^{K} \setminus\{\varepsilon_{K}\}, \hspace{0.5cm}
\mathcal{H}_{K}=\bigoplus\limits_{\alpha \in \Gamma^{K}} \mathfrak{g}^{\alpha} .$$
~\\

\begin{lemme}\label{kostant2} Soit $K,K'$ dans $\mathcal{K}(\Pi)$, $\alpha,\beta$ dans $\Gamma^{K}$ et
 $\gamma$ dans $\Gamma^{K'}$.
\begin{description}
\item[{\rm (i)}] On a $\Gamma^{K}=\Delta_{+}^{K} \setminus \{ \delta
  \in \Delta_{+}^{K} \ | \ \langle \delta, \varepsilon_{S}^{\vee} \rangle =0\}$.
\item[{\rm (ii)}] L'ensemble $\Delta_+$ est la r{\'e}union disjointe des
  $\Gamma^{K''}$ pour $K''$ dans $\mathcal{K}(\Pi)$, et $\mathcal{H}_{K}$
  est une alg{\`e}bre de Heisenberg de centre $\mathfrak{g}^{\varepsilon_{K}}$.
\item[{\rm (iii)}] Si $\alpha + \beta$ appartient {\`a} $\Delta$, alors $\alpha
  +\beta=\varepsilon_{K}$.
\item[{\rm (iv)}]  Si $\alpha + \gamma$ appartient {\`a} $\Delta$, alors
  ou bien $K
  \subset K'$ et $\alpha + \gamma$ appartient {\`a} $\Gamma^{K'}$, ou bien $K'
  \subset K$ et $\alpha + \gamma$ appartient {\`a} $\Gamma^{K}$.\\

\end{description}

\end{lemme}

\begin{remarque}\label{heisenberg} Notons que si $K$ est un {\'e}l{\'e}ment de $\mathcal{K}(\Pi)$, alors pour
toute racine $\alpha$ de $\Gamma_{0}^{K}$, il existe une unique racine
$\beta$ de $\Gamma_{0}^{K}$ telle que $\alpha +\beta=\varepsilon_K$ et on
a\,:
$$\langle \alpha , {\varepsilon_K}^{\vee} \rangle=
\langle \beta , {\varepsilon_K}^{\vee} \rangle =1 .$$ 
Cela r{\'e}sulte du point (ii) du lemme \ref{kostant2} et des r{\'e}sultats de
\cite{Bourbaki} rappel{\'e}s pr{\'e}c{\'e}demment.\\

\end{remarque}

Le cardinal de $\mathcal{K}(\Pi)$ d{\'e}pend de $\mathfrak{g}$ mais pas de
$\mathfrak{h}$ ou de $\Pi$. On note $k_{\mathfrak{g}}$ cet entier. La table \ref{kg} donne la valeur de $k_{\mathfrak{g}}$
pour les diff{\'e}rents types d'alg{\`e}bres de Lie simples.\\
{\footnotesize
\begin{table}[h]
\begin{center}
\begin{tabular}{|c|c|c|c|c|c|c|c|c|c|}
\hline &
\textbf{ $A_{l}, l \geq 1$ } &
\textbf{ $B_{l}, l \geq 2$ } &  
\textbf{ $C_{l}, l \geq 3$ } & 
\textbf{ $D_{l}, l \geq 4$ } &
\textbf{ $E_6$ } &
\textbf{ $E_7$ } & 
\textbf{ $E_8$ } &
\textbf{ $F_4$ } &
\textbf{ $G_2$ } \\
\hline
 & & & & & & & & & \\
$k_{\mathfrak{g}}$ & $\left[ \displaystyle{\frac{l+1}{2}} \right]$ & $l$ & $l$ &
$2 \left[ \displaystyle{\frac{l}{2}} \right]$ &
$4$ & $7$ & $8$ & $4$ & $2$ \\
 & & & & & & & & & \\
\hline
\end{tabular}
~\\
\caption{\label{kg} $k_{\mathfrak{g}}$ pour les alg{\`e}bres de Lie simples.}

\end{center}
\end{table}}

On descrit dans les tables \ref{classic} et \ref{except} l'ensemble $ \{\varepsilon_K
\  , \ K \in \mathcal{K}(\Pi)\}$ pour les diff{\'e}rents types d'alg{\`e}bres
de Lie simples complexes.\\

\hspace{-4cm}
{\footnotesize
\begin{sidewaystable}[ht]
\begin{center}
\begin{tabular}{|c|c|c|}
\hline
$A_l$, $l \geq 1$ &  
\begin{pspicture}(-2.5,-0)(3.5,1)
\pscircle(-2,0){1mm}
\pscircle(-1,0){1mm}
\pscircle(0,0){1mm}
\pscircle(1,0){1mm}
\pscircle(2,0){1mm}
\pscircle(3,0){1mm}
\psline(-1.9,0)(-1.1,0)
\psline(-0.1,0)(-0.9,0)
\psline[linestyle=dotted](0.1,0)(0.9,0)
\psline(1.1,0)(1.9,0)
\psline(2.1,0)(2.9,0)

\rput[b](-2,0.2){$\beta_1$}
\rput[b](-1,0.2){$\beta_2$}
\rput[b](2,0.2){$\beta_{l-1}$}
\rput[b](3,0.2){$\beta_l$}
\end{pspicture}
& $\left\lbrace \beta_i + \cdots + \beta_{i+(l-2i +1)} \ , \ 1 \leq i \leq
    \left[ \displaystyle{\frac{l+1}{2}}\right] \right\rbrace.$\\
&&\\
\hline
$B_l$,  &
\begin{pspicture}(-2.5,-0)(3.5,1)
\pscircle(-2,0){1mm}
\pscircle(-1,0){1mm}
\pscircle(0,0){1mm}
\pscircle(1,0){1mm}
\pscircle(2,0){1mm}
\pscircle(3,0){1mm}
\psline(-1.9,0)(-1.1,0)
\psline(-0.1,0)(-0.9,0)
\psline[linestyle=dotted](0.1,0)(0.9,0)
\psline(1.1,0)(1.9,0)
\psline(2.1,0.05)(2.9,0.05)
\psline(2.1,-0.05)(2.9,-0.05)
\rput[b](-2,0.2){$\beta_1$}
\rput[b](-1,0.2){$\beta_2$}
\rput[b](2,0.2){$\beta_{l-1}$}
\rput[b](3,0.2){$\beta_l$}
\rput(2.5,0){$>$}
\end{pspicture}
&$\left\lbrace \beta_i + 2 \beta_{i+1} + \cdots + 2 \beta_{l} \ , \ 1
  \leq i \leq l-1, \ i \textrm{ impair } \right\rbrace$ \\
$l \geq 2$ &&$\cup  \left\lbrace
  \beta_i \  , \ 1
  \leq i \leq l, \ i \textrm{ impair } \right\rbrace.$\\
&&\\
\hline
$C_l$, $l \geq 3$ &
\begin{pspicture}(-2.5,-0)(3.5,1)
\pscircle(-2,0){1mm}
\pscircle(-1,0){1mm}
\pscircle(0,0){1mm}
\pscircle(1,0){1mm}
\pscircle(2,0){1mm}
\pscircle(3,0){1mm}
\psline(-1.9,0)(-1.1,0)
\psline(-0.1,0)(-0.9,0)
\psline[linestyle=dotted](0.1,0)(0.9,0)
\psline(1.1,0)(1.9,0)
\psline(2.1,0.05)(2.9,0.05)
\psline(2.1,-0.05)(2.9,-0.05)
\rput[b](-2,0.2){$\beta_1$}
\rput[b](-1,0.2){$\beta_2$}
\rput[b](2,0.2){$\beta_{l-1}$}
\rput[b](3,0.2){$\beta_l$}
\rput(2.5,0){$<$}
\end{pspicture}
& $\left\lbrace 2 \beta_i + \cdots +2 \beta_{l-1} + \beta_{l} \ , \ 1
  \leq i \leq l-1 \right\rbrace$\\
$l \geq 3$ & & $\cup  \ \left\lbrace
  \beta_l \right\rbrace.$ \\
&&\\
\hline 
$D_l$, $l$ pair &
\begin{pspicture}(-2.5,-0.2)(3.6,1.1)
\pscircle(-2,0){1mm}
\pscircle(-1,0){1mm}
\pscircle(0,0){1mm}
\pscircle(1,0){1mm}
\pscircle(2,0){1mm}
\pscircle(3,0.7){1mm}
\pscircle(3,-0.7){1mm}
\psline(-1.9,0)(-1.1,0)
\psline(-0.1,-0)(-0.9,0)
\psline[linestyle=dotted](0.1,0)(0.9,0)
\psline(1.1,0)(1.9,0)
\psline(2.07,0.05)(2.92,0.68)
\psline(2.07,-0.05)(2.92,-0.68)
\rput[b](-2,0.2){$\beta_1$}
\rput[b](-1,0.2){$\beta_2$}
\rput[b](2,0.2){$\beta_{l-2}$}
\rput[l](3.1,0.8){$\beta_l$}
\rput[l](3.1,-0.8){$\beta_{l-1}$}
\end{pspicture}
& $\left\lbrace \beta_i + 2 \beta_{i+1} + \cdots +2 \beta_{l-2} +
  \beta_{l-1} +\beta_{l} , 1
  \leq i \leq l-3, \ i \textrm{ impair } \right\rbrace $\\
$l \geq 3$ & & $\cup  \left\lbrace
  \beta_i \  , \ 1
  \leq i \leq l, \ i \textrm{ impair } \right\rbrace  \ \cup  \ \left\lbrace
  \beta_l \right\rbrace.$ \\
&&\\
\hline
$D_l$, $l$ impair &
\begin{pspicture}(-2.5,-0.2)(3.6,1.1)
\pscircle(-2,0){1mm}
\pscircle(-1,0){1mm}
\pscircle(0,0){1mm}
\pscircle(1,0){1mm}
\pscircle(2,0){1mm}
\pscircle(3,0.7){1mm}
\pscircle(3,-0.7){1mm}
\psline(-1.9,0)(-1.1,0)
\psline(-0.1,-0)(-0.9,0)
\psline[linestyle=dotted](0.1,0)(0.9,0)
\psline(1.1,0)(1.9,0)
\psline(2.07,0.05)(2.92,0.68)
\psline(2.07,-0.05)(2.92,-0.68)
\rput[b](-2,0.2){$\beta_1$}
\rput[b](-1,0.2){$\beta_2$}
\rput[b](2,0.2){$\beta_{l-2}$}
\rput[l](3.1,0.8){$\beta_l$}
\rput[l](3.1,-0.8){$\beta_{l-1}$}
\end{pspicture}
 &
$\left\lbrace \beta_i + 2 \beta_{i+1} + \cdots +2 \beta_{l-2} +
  \beta_{l-1} +\beta{l} \ , \ 1
  \leq i \leq l-3, \ i \textrm{ impair } \right\rbrace$ \\
$l \geq 5$ & &$\cup  \left\lbrace
  \beta_i \  , \ 1
  \leq i \leq l-2, \ i \textrm{ impair } \right\rbrace  \ \cup  \ \left\lbrace
  \beta_{l-2} + \beta_{l-1} +\beta_{l} \right\rbrace.$ \\
&&\\
\hline
\end{tabular}
~\\
\caption{\label{classic} $\{\varepsilon_K
\  , \ K \in \mathcal{K}(\Pi)\}$ pour les alg{\`e}bres de Lie simples classiques.}
\end{center}
\end{sidewaystable}

}

\hspace{-4cm}
{\footnotesize
\begin{sidewaystable}[ht]

\begin{center}
\hspace{-4cm}
\begin{tabular}{|c|c|c|}
\hline
$E_6$ &
\begin{pspicture}(-2.5,0)(2.5,1)

\pscircle(-2,0){1mm}
\pscircle(-1,0){1mm}
\pscircle(0,0){1mm}
\pscircle(1,0){1mm}
\pscircle(2,0){1mm}
\pscircle(0,-1){1mm}

\psline(-1.9,0)(-1.1,0)
\psline(-0.1,0)(-0.9,0)
\psline(0.1,0)(0.9,0)
\psline(1.1,0)(1.9,0)
\psline(0,-0.1)(0,-0.9)

\rput[b](-2,0.2){$\beta_1$}
\rput[b](-1,0.2){$\beta_3$}
\rput[b](0,0.2){$\beta_4$}
\rput[b](1,0.2){$\beta_5$}
\rput[b](2,0.2){$\beta_6$}
\rput[r](-0.2,-1){$\beta_2$}
\end{pspicture}
& $\left(\begin{array}{ccccc}
1&2&3&2&1\\
&&2&&
\end{array} \right) \ , \ \left(\begin{array}{ccccc}
1&1&1&1&1\\
&&0&&
\end{array}\right)  \ , \
\left(\begin{array}{ccccc}
0&1&1&1&0\\
&&0&&
\end{array}\right) \ , \  \left(\begin{array}{ccccc}
0&0&1&0&0\\
&&0&&
\end{array}\right).$\\
&&\\
&&\\
&&\\
\hline
$E_7$ &
\begin{pspicture}(-2.5,0)(3.5,1)

\pscircle(-2,0){1mm}
\pscircle(-1,0){1mm}
\pscircle(0,0){1mm}
\pscircle(1,0){1mm}
\pscircle(2,0){1mm}
\pscircle(3,0){1mm}

\pscircle(0,-1){1mm}

\psline(-1.9,0)(-1.1,0)
\psline(-0.1,0)(-0.9,0)
\psline(0.1,0)(0.9,0)
\psline(1.1,0)(1.9,0)
\psline(2.1,0)(2.9,0)
\psline(0,-0.1)(0,-0.9)

\rput[b](-2,0.2){$\beta_1$}
\rput[b](-1,0.2){$\beta_3$}
\rput[b](0,0.2){$\beta_4$}
\rput[b](1,0.2){$\beta_5$}
\rput[b](2,0.2){$\beta_6$}
\rput[b](3,0.2){$\beta_7$}
\rput[r](-0.2,-1){$\beta_2$}
\end{pspicture}
& $\left(\begin{array}{cccccc}
2&3&4&3&2&1\\
&&2&&&
\end{array}\right) \ , \ \left(\begin{array}{cccccc}
0&1&2&2&2&1\\
&&1&&&
\end{array}\right)$\\
&&\\
&&  $\left(\begin{array}{cccccc}
0&1&2&1&0&0\\
&&1&&&
\end{array}\right) \ , \ \left(\begin{array}{cccccc}
0&0&0&0&0&1\\
&&0&&&
\end{array}\right)$\\
&&\\
&&$\left(\begin{array}{cccccc}
0&0&0&0&0&0\\
&&1&&&
\end{array}\right)\ , \  \left(\begin{array}{cccccc}
0&1&0&0&0&0\\
&&0&&&
\end{array}\right) \ , \ \left(\begin{array}{cccccc}
0&0&0&1&0&0\\
&&0&&&
\end{array}\right) $\\
&&\\
\hline
$E_8$ & 
\begin{pspicture}(-2.5,0)(4.5,1)

\pscircle(-2,0){1mm}
\pscircle(-1,0){1mm}
\pscircle(0,0){1mm}
\pscircle(1,0){1mm}
\pscircle(2,0){1mm}
\pscircle(3,0){1mm}
\pscircle(4,0){1mm}

\pscircle(0,-1){1mm}

\psline(-1.9,0)(-1.1,0)
\psline(-0.1,0)(-0.9,0)
\psline(0.1,0)(0.9,0)
\psline(1.1,0)(1.9,0)
\psline(2.1,0)(2.9,0)
\psline(3.1,0)(3.9,0)
\psline(0,-0.1)(0,-0.9)

\rput[b](-2,0.2){$\beta_1$}
\rput[b](-1,0.2){$\beta_3$}
\rput[b](0,0.2){$\beta_4$}
\rput[b](1,0.2){$\beta_5$}
\rput[b](2,0.2){$\beta_6$}
\rput[b](3,0.2){$\beta_7$}
\rput[b](4,0.2){$\beta_8$}
\rput[r](-0.2,-1){$\beta_2$}
\end{pspicture}
& $ \left(\begin{array}{ccccccc}
2&4&6&5&4&3&2\\
&&3&&&&
\end{array}\right)\ , \ \left(\begin{array}{ccccccc}
2&3&4&3&2&1&0\\
&&2&&&&
\end{array}\right)\ , \ \left(\begin{array}{ccccccc}
0&1&2&2&2&1&0\\
&&1&&&&
\end{array}\right)$\\
&&\\
&&$\left(\begin{array}{ccccccc}
0&1&2&1&0&0&0\\
&&1&&&&
\end{array}\right) \ , \  \left(\begin{array}{ccccccc}
0&0&0&0&0&1&0\\
&&0&&&&
\end{array}\right)\ , \ \left(\begin{array}{ccccccc}
0&0&0&0&0&0&0\\
&&1&&&&
\end{array}\right)$\\
&&\\
&& $\left(\begin{array}{ccccccc}
0&1&0&0&0&0&0\\
&&0&&&&
\end{array}\right) \ , \ \left(\begin{array}{ccccccc}
0&0&0&1&0&0&0\\
&&0&&&&
\end{array}\right)$\\
&&\\
\hline
$F_4$&
\begin{pspicture}(-2.5,-0)(1.5,1)

\pscircle(-2,0){1mm}
\pscircle(-1,0){1mm}
\pscircle(0,0){1mm}
\pscircle(1,0){1mm}

\psline(-1.9,0)(-1.1,0)
\psline(-0.09,0.05)(-0.91,0.05)
\psline(-0.09,-0.05)(-0.91,-0.05)
\psline(0.1,0)(0.9,0)

\rput[b](-2,0.2){$\beta_{1}$}
\rput[b](-1,0.2){$\beta_{2}$}
\rput[b](0,0.2){$\beta_{3}$}
\rput[b](1,0.2){$\beta_{4}$}

\rput(-0.5,0){$>$}
\end{pspicture}
& 
$\left(\begin{array}{cccc}
2&3&4&2\\
\end{array}\right) \ , \  \left(\begin{array}{cccc}
0&1&2&2
\end{array}\right)\ , \ \left(\begin{array}{cccc}
0&1&2&0
\end{array}\right)\ , \ 
\left(\begin{array}{cccc}
0&1&0&0
\end{array}\right) 
.$\\
&&\\
\hline
$G_2$ &
\begin{pspicture}(-1.5,-0)(1.5,1)

\pscircle(-1,0){1mm}
\pscircle(0,0){1mm}

\psline(-0.09,0.06)(-0.91,0.06)
\psline(-0.09,0)(-0.91,0)
\psline(-0.09,-0.06)(-0.91,-0.06)

\rput[b](-1,0.2){$\beta_{1}$}
\rput[b](0,0.2){$\beta_{2}$}

\rput(-0.5,0){$>$}
\end{pspicture}
&$\left(\begin{array}{cccc}
2&3\end{array}\right)\ , \ \left(\begin{array}{cc}
1&2\\
\end{array}\right) .$\\
&&\\
\hline
\end{tabular} 
~\\
\caption{\label{except} $\{\varepsilon_K
\  , \ K \in \mathcal{K}(\Pi)\}$ pour les alg{\`e}bres de Lie simples exceptionnelles.}

\end{center}
\end{sidewaystable}
}

\subsection{} On {\'e}tudie dans ce paragraphe la fa{\c c}on dont
l'involution $\theta$ agit sur l'ensemble $\mathcal{K}(\Pi)$. Posons
\begin{eqnarray*}
\mathcal{K}^{''}(\Pi) & = & \{K \in \mathcal{K}(\Pi) \ | \
{\varepsilon_{K}}_{|_{\mathfrak{a}}} \not= 0\}\\
& = & \{K \in \mathcal{K}(\Pi) \ | \
\varepsilon_{K} \in \Delta_{+}^{''}\},\\  
\end{eqnarray*}
et 
\begin{eqnarray*}
\mathcal{K}^{'}(\Pi) & = & \{K \in \mathcal{K}(\Pi) \ | \
{\varepsilon_{K}}_{|_{\mathfrak{a}}} = 0\}\\
& = & \{K \in \mathcal{K}(\Pi) \ | \
\varepsilon_{K} \in \Delta_{+}^{'}\}  .
\end{eqnarray*}
~\\

Si $K$ appartient {\`a} $\mathcal{K}^{'}(\Pi)$, alors $\theta
\varepsilon_K=\varepsilon_K$. On cherche maintenant {\`a} {\'e}tudier la fa{\c c}on dont
l'involution $\theta$ agit sur l'ensemble $\mathcal{K}^{''}(\Pi)$. On
introduit pour cela la d{\'e}finition suivante\,:

\begin{definition}\label{prop_q} On dira que le couple
  $(\mathfrak{h}_0,\Pi)$ {\rm a la propri{\'e}t{\'e} $(P)$} si, pour tout $K$
  dans $\mathcal{K}^{''}(\Pi)$, il existe un unique {\'e}l{\'e}ment $L$ dans
  $\mathcal{K}^{''}(\Pi)$ tel que, $-\theta
\varepsilon_K=\varepsilon_L$.
\end{definition}

On introduit les deux sous-ensembles suivants de $\mathcal{K}^{''}(\Pi)$\,:

\begin{eqnarray*}
\mathcal{K}_{{\rm r\acute{e}el}}(\Pi) & = & \{K \in \mathcal{K}^{''}(\Pi) \ | \
-\theta \varepsilon_K=\varepsilon_K\}\\
& = & \{K \in \mathcal{K}(\Pi) \ | \
\varepsilon_{K} \textrm{ est r{\'e}elle } \}
\end{eqnarray*}
et 
\begin{eqnarray*}
\mathcal{K}_{{\rm comp}}(\Pi) & = & \{K \in \mathcal{K}^{''}(\Pi) \ | \
-\theta \varepsilon_K \not = \varepsilon_K\} .
\end{eqnarray*}
~\\

Si $(\mathfrak{h}_0,\Pi)$ a la propri{\'e}t{\'e} $(P)$, alors $\theta$ induit
  une involution dans $\mathcal{K}^{''}(\Pi)$, que l'on note encore
  $\theta$. L'ensemble $\mathcal{K}_{{\rm r\acute{e}el}}(\Pi)$ est
  alors l'ensemble des points fixes pour cette involution et, puisque $\theta$ est une involution, le cardinal
de $\mathcal{K}_{{\rm comp}}(\Pi)$ est pair et l'on peut choisir un
sous-ensemble $\mathcal{K}_{{\rm comp}}^{+}(\Pi)$ de
$\mathcal{K}_{{\rm comp}}(\Pi)$ de sorte que, 
$$\mathcal{K}_{{\rm comp}}(\Pi)=\mathcal{K}_{{\rm comp}}^{+}(\Pi)
\cup \theta \mathcal{K}_{{\rm comp}}^{+}(\Pi) .$$
~\\

\begin{proposition}\label{Kostant_Cartan1} Soit $\mathfrak{h}_0$ une sous-alg{\`e}bre de Cartan de $\mathfrak{g}_0$
stable par $\theta$. Si $\mathfrak{h}_0$ est maximalement
non-compacte, alors  $(\mathfrak{h}_0,\Pi)$ a la propri{\'e}t{\'e} $(P)$.
\end{proposition}

\begin{proof} Il s'agit de prouver que, pour tout $K$
  dans $\mathcal{K}^{''}(\Pi)$, il existe un unique {\'e}l{\'e}ment $L$ dans
  $\mathcal{K}^{''}(\Pi)$ tel que, $-\theta
\varepsilon_K=\varepsilon_L$. On d{\'e}montre cette assertion par
  induction sur l'inclusion dans $\mathcal{K}^{''}(\Pi)$.\\
\\
1) Soit $K$ dans $\mathcal{K}^{''}(\Pi)$ qui est maximal pour
  l'inclusion et supposons par l'absurde l'assertion fausse. Alors il
  existe un {\'e}l{\'e}ment $L$ dans $\mathcal{K}^{''}(\Pi)$ tel que $- \theta
\varepsilon_K$ appartienne {\`a} $\Gamma_{0}^{L}$. Soit $\alpha$ dans
  $\Gamma_{0}^{L}$ tel que,
\begin{eqnarray}\label{p1} \alpha + (- \theta \varepsilon_K)= \varepsilon_L.
\end{eqnarray}
On a alors,
\begin{eqnarray}\label{p2} (\theta \alpha) + (- \theta \varepsilon_L)= \varepsilon_K.
\end{eqnarray}

Si $\theta \alpha$ est une racine positive, alors $\alpha$ appartient
{\`a} $\Delta_{+}^{'}$ et $\theta \alpha=\alpha$. La relation (\ref{p2})
entraine alors que $\alpha$ appartient {\`a} $\Gamma_{0}^{K}$, d'apr{\`e}s le
lemme \ref{kostant2}, (iii) et (iv). On en
d{\'e}duit que $K=L$, car l'intersection $\Gamma^{L} \cap \Gamma^{K}$ est
non vide, puis que $\varepsilon_K+\theta \varepsilon_K$ est une
racine. C'est une racine imaginaire non-compacte. En effet, c'est
clairement une racine imaginaire et, d'apr{\`e}s le lemme \ref{struct1} (i), le crochet
$[X_{\varepsilon_K},\theta X_{\varepsilon_K}]$ est un {\'e}l{\'e}ment non nul de
$\mathfrak{g}^{\varepsilon_K+\theta \varepsilon_K}$, qui est contenu dans
$\mathfrak{p}$, donc la racine imaginaire $\varepsilon_K+\theta
\varepsilon_K$ est non-compacte. Ceci est alors en contradiction avec
le lemme \ref{max_compacte}, car $\mathfrak{h}_0$ est maximalement
non-compacte.

Par suite, $-\theta \alpha$ est une racine positive et la
relation (\ref{p2}) donne\,:
$$\varepsilon_K + (-\theta \alpha) = - \theta \varepsilon_L.$$
Puisque $K$ est maximal, il r{\'e}sulte du lemme \ref{kostant2} (iv), que
les racines $-\theta \alpha$ et $- \theta \varepsilon_L$ appartiennent
{\`a} $\Delta_{+}^{K}$. La relation pr{\'e}c{\'e}dente contredit alors que
$\varepsilon_K$ est la plus grande racine de $\Delta_{+}^{K}$.\\
\\
2) Soit $K$ dans $\mathcal{K}^{''}(\Pi)$. Supposons
l'assertion d{\'e}montr{\'e}e pour tout {\'e}l{\'e}ment dans $\mathcal{K}^{''}(\Pi)$
contenant strictement $K$ et supposons par l'absurde l'assertion
fausse pour $K$. Alors il
  existe un {\'e}l{\'e}ment $L$ dans $\mathcal{K}^{''}(\Pi)$ tel que $- \theta
\varepsilon_K$ appartienne {\`a} $\Gamma_{0}^{L}$. Soit $\alpha$ comme
pr{\'e}c{\'e}demment dans
  $\Gamma_{0}^{L}$ tel que,
\begin{eqnarray}\label{p3} \alpha + (- \theta \varepsilon_K)= \varepsilon_L.
\end{eqnarray}
En raisonnant comme dans le cas (1), on obtient que $-\theta \alpha$
est une racine positive et on a,
\begin{eqnarray}\label{p4} \varepsilon_K + (-\theta \alpha) = - \theta \varepsilon_L.
\end{eqnarray}
Soit $M$ dans $\mathcal{K}^{''}(\Pi)$ tel que $-\theta \alpha$
appartienne {\`a} $\Gamma_{0}^{M}$. D'apr{\`e}s le lemme \ref{kostant2} (iv),
il y a deux cas\,: si $M \subseteq K$, on obtient une contradition en
raisonnant comme dans le cas (1). Sinon, alors $- \theta
\varepsilon_L$ appartient {\`a}  $\Gamma^{M}$ et d'apr{\`e}s l'hypoth{\`e}se de
r{\'e}currence, on en d{\'e}duit que $ -\theta \varepsilon_L=
\varepsilon_M$. La relation (\ref{p4}) contredit alors que les racines $\varepsilon_K$ et $\varepsilon_M$ sont fortement
orthogonales.\\

Par induction sur l'inclusion, la proposition est d{\'e}montr{\'e}e.       
\end{proof}

En particulier, puisque $\widehat{\mathfrak{h}}_0$ est maximalement non-compacte, la
proposition pr{\'e}c{\'e}dente assure que le couple
$(\widehat{\mathfrak{h}}_0,\widehat{\Pi})$ a la propri{\'e}t{\'e} $(P)$. 

\subsection{} On s'int{\'e}resse dans ce paragraphe {\`a} la condition
suivante\,:\\

$({\tt *}) \hspace{4cm} \mathcal{K}(\Pi^{'})=\mathcal{K}^{'}(\Pi),$\\
\\
qui {\'e}quivaut a l'inclusion\,:
$\mathcal{K}(\Pi^{'}) \subset \mathcal{K}(\Pi)$. Notons que cette
condition est en particulier remplie si $\Delta^{'}$ est vide, ce qui
est le cas d{\`e}s que $\mathfrak{m}$ est ab{\'e}lienne. 

\begin{lemme}\label{Kostant_Cartan2} Soit $\mathfrak{h}_0$ une sous-alg{\`e}bre de Cartan de $\mathfrak{g}_0$
stable par $\theta$. On suppose que $(\mathfrak{h}_0,\Pi)$ a la
propri{\'e}t{\'e} $(P)$ et satisfait {\`a} la condition $({\tt *})$. Alors\,:
\begin{itemize} \item[{\rm (i)}] $\forall \alpha \in \Delta_{+}^{'}$, $\forall K \in \mathcal{K}^{''}(\Pi)$,
  on a\,: $\alpha \pm \varepsilon_K \not\in \Delta$,
\item[{\rm (ii)}] $\forall K \in \mathcal{K}^{''}(\Pi)$, $\Gamma^{K}
   \subset \Delta_{+}^{''}$,
\item[{\rm (ii)}] Si $\mathfrak{h}_0$ n'est pas maximalement compacte,
  alors $\mathcal{K}_{{\rm
      r\acute{e}el}}(\Pi)$ n'est pas vide.
\end{itemize}
\end{lemme}

\begin{proof} (i) Supposons par
  l'absurde qu'il existe $\alpha$ dans $\Delta_{+}^{'}$ et $K$ dans $\mathcal{K}^{''}(\Pi)$ tel que
  $\alpha + \varepsilon_K \in \Delta$. Puisque $\alpha$ appartient {\`a}
  $\Delta_{+}^{'}$, il existe $L$ dans
  $\mathcal{K}(\Pi^{'})=\mathcal{K}^{'}(\Pi)$ tel que $\alpha$
  appartient {\`a} $\Gamma_{0}^{L}$. La relation $\alpha + \varepsilon_K
  \in \Delta$ entraine alors l'inclusion $K \subset L$, d'apr{\`e}s le
  lemme \ref{kostant2} (iv). Ceci est absurde car les racines de
  $\Delta_{+}^{''}$ sont plus grandes que celles de $\Delta_{+}^{'}$
  par construction de $\Delta_+$. 

Supposons maintenant par
  l'absurde qu'il existe $\alpha$ dans $\Delta_{+}^{'}$ et $K$ dans $\mathcal{K}^{''}(\Pi)$ tel que
  $\alpha - \varepsilon_K \in \Delta$. Si $\varepsilon_K-\alpha
  =\beta$ appartient {\`a} $\Delta_{+}$, alors il r{\'e}sulte du lemme
  \ref{kostant2} (iv) que $\alpha$ et $\beta$ appartiennent {\`a} $\Gamma_{0}^{K}$. Comme
  par ailleurs, $\alpha$ appartient {\`a} $\Delta_{+}^{'}$, il existe un
  unique $L$ dans
  $\mathcal{K}(\Pi^{'})=\mathcal{K}^{'}(\Pi)$ tel que $\alpha$
  appartient {\`a} $\Gamma_{0}^{L}$, d'o{\`u} la contradiction. On en d{\'e}duit
  que $\alpha - \varepsilon_K=\beta$ est une racine positive. Mais
  alors $\alpha$ est plus grande que $\varepsilon_K$, ce qui est
  absurde, pusique $\alpha$ est une racine imaginaire.\\
\\
(ii) Soit $K  \in \mathcal{K}^{''}(\Pi)$ et $\alpha \in
\Gamma_{0}^{K}$. Il existe $\beta \in \Gamma_{0}^{K}$ tel que $\alpha
+ \beta=\varepsilon_K$.  Supposons par l'absurde que $\alpha$ appartient {\`a}
$\Delta_{+}^{'}$. Puisque $(\mathfrak{h}_0,\Pi)$ a la
propri{\'e}t{\'e} $(P)$, on a alors\,:
$$(-\theta \beta)=\varepsilon_{\theta K} +\alpha,$$
ce qui contredit (i).\\
\\
(iii) Puisque $\mathfrak{h}_0$ n'est pas maximalement
compacte, il r{\'e}sulte du lemme \ref{max_compacte} qu'il existe un
{\'e}l{\'e}ment $K$ dans $\mathcal{K}^{''}(\Pi)$ tel que $\Gamma^{K}$
poss{\`e}de une racine r{\'e}elle $\alpha$. Il suffit de montrer que $\varepsilon_K$ est
une racine r{\'e}elle. On peut supposer que $\alpha$ appartient  {\`a}
$\Gamma_{0}^{K}$. Il existe $\beta$ dans $\Gamma_{0}^{K}$ tel
que,
\begin{eqnarray}\label{q1} \alpha + \beta =\varepsilon_K.
\end{eqnarray}
Puisque $(\mathfrak{h}_0,\Pi)$ a la
propri{\'e}t{\'e} $(P)$, on en d{\'e}duit la relation suivante\,:
\begin{eqnarray}\label{q2} \alpha + (- \theta \beta)
  =\varepsilon_{\theta K}.
\end{eqnarray}
D'apr{\`e}s (ii), $- \theta \beta$ est une
racine positive et la relation (\ref{q2}) implique que $\alpha$ et  $-
\theta \beta$ appartiennent {\`a} $\Gamma_{0}^{\theta K}$. L'intersection
$\Gamma_{0}^{\theta K} \cap \Gamma_{0}^{K}$ est alors non vide,
d'o{\`u} il vient $\theta K=K$. Par suite, $\varepsilon_K$ est r{\'e}elle.
\end{proof}

\section{Utilisation des transformations de Cayley}\label{part_cayley}

On utilise dans cette partie les transformations de
Cayley afin d'obtenir la proposition \ref{calcul_k}. On trouve davantage
de pr{\'e}cisions concernant ces transformations dans \cite{Knapp}. Si $\mathfrak{h}_0$ et $\mathfrak{h}^{'}_0$ sont deux sous-alg{\`e}bres de Cartan de $\mathfrak{g}_0$
stables par $\theta$, leurs complexifi{\'e}es $\mathfrak{h}$ et $\mathfrak{h}^{'}$ sont
conjugu{\'e}es dans $\mathfrak{g}$. Les transformations de Cayley permettent de
construire explicitement, et {\'e}tapes par
{\'e}tapes, un automorphisme de $\mathfrak{g}$ qui conjugue ces deux
sous-alg{\`e}bres. Il existe deux types de transformations de Cayley {\`a}
partir d'une sous-alg{\`e}bre de Cartan $\theta$-stable\,:
\begin{description} 
\item[(i)] avec une racine imaginaire non-compacte $\beta$, on
  construit une nouvelle sous-alg{\`e}bre de Cartan dont l'intersection
  avec $\mathfrak{p}_0$ augmente de $1$ en dimension. On notera
  ${\bf c}_{\beta}$ la transformation correspondante.   
\item[(ii)] avec une racine r{\'e}elle $\alpha$, on
  construit une nouvelle sous-alg{\`e}bre de Cartan dont l'intersection
  avec $\mathfrak{p}_0$ diminue de $1$ en dimension. On notera
  ${\bf d}_{\alpha}$ la transformation correspondante.
\end{description}

Seules les transformations de type ${\bf d}_{\alpha}$ vont intervenir 
dans la suite. On donne ici les r{\'e}sultats n{\'e}cessaires concernant ces
transformations. Soit $\mathfrak{h}_0$ une sous-alg{\`e}bre de Cartan de $\mathfrak{g}_0$
stable par $\theta$. Si $\alpha$ est une racine
r{\'e}elle de $\Delta$, on pose\,:
$${\bf d}_{\alpha}={\rm Ad}(\exp i \frac{\pi}{4}(\theta
X_{\alpha}-X_{\alpha})) .$$
Pour $\beta$ dans $\Delta$, on note ${\bf d}_{\alpha}(\beta)$ la forme
lin{\'e}aire de ${\bf d}_{\alpha}(\mathfrak{h})$ qui {\`a} $H$ dans ${\bf
  d}_{\alpha}(\mathfrak{h})$ associe $\beta({\bf
  d}_{\alpha}^{-1}(H))$; c'est un {\'e}l{\'e}ment du syst{\`e}me de racines associ{\'e} au couple
$(\mathfrak{g},{\bf d}_{\alpha}(\mathfrak{h}))$. On a
$$\mathfrak{g}_0 \cap  {\bf d}_{\alpha}(\mathfrak{h})=\ker
(\alpha_{|_{\mathfrak{h}_{0}}}) \oplus \mathbb{R}( X_{\alpha} +\theta
X_{\alpha}),$$
d'o{\`u} on tire la relation\,:
\begin{eqnarray}\label{rel_dim}
\dim ({\bf d}_{\alpha}(\mathfrak{h}) \cap \mathfrak{p})=\dim
(\mathfrak{h} \cap \mathfrak{p})-1 .
\end{eqnarray}
On dispose en outre de la relation suivante\,:
\begin{eqnarray}\label{1}
{\bf d}_{\alpha}(H_{\alpha}) = i \mu (X_{\alpha} +\theta
 X_{\alpha}),
\end{eqnarray}
o{\`u} $\mu$ est r{\'e}el non nul. Enfin, compte tenu de l'expression de ${\bf d}_{\alpha}$, il
est clair que si $\beta$ est une racine fortement orthogonale {\`a}
$\alpha$ et diff{\'e}rente de $\alpha$, alors ${\bf
  d}_{\alpha}(H_{\beta})=H_{\beta}$. 

\begin{lemme}\label{cayley} Soit $\mathfrak{h}_0$ une sous-alg{\`e}bre de Cartan de $\mathfrak{g}_0$
stable par $\theta$. Soit $K$ dans $\mathcal{K}_{{\rm
      r\acute{e}el}}(\Pi)$. Alors la sous-alg{\`e}bre ${\bf d}_{\varepsilon_K}(\mathfrak{h}_0)$ est une sous-alg{\`e}bre de Cartan de
$\mathfrak{g}_0$ qui est stable par $\theta$. De plus, on a\,:
\begin{eqnarray*}
\mathcal{K}^{''}({\bf d}_{\varepsilon_K}(\Pi)) & = & \{ \varepsilon_L \
| \ L \in \mathcal{K}^{''}(\Pi) \} \setminus \{\varepsilon_K\},\\
\mathcal{K}^{'}({\bf d}_{\varepsilon_K}(\Pi)) & = & \{ \varepsilon_L \
| \ L \in \mathcal{K}^{'}(\Pi) \} \cup \{{\bf
  d}_{\varepsilon_K}(\varepsilon_K)\}.
\end{eqnarray*}
De plus, si $(\mathfrak{h}_0,\Pi)$ a la propri{\'e}t{\'e} $(P)$, alors
la sous-alg{\`e}bre $({\bf d}_{\varepsilon_K}(\mathfrak{h}_0),{\bf
  d}_{\varepsilon_K}(\Pi))$ a la propri{\'e}t{\'e} $(P)$ et, si
$(\mathfrak{h}_0,\Pi)$ satisfait {\`a} la condition $({\tt *})$, alors
la sous-alg{\`e}bre $({\bf d}_{\varepsilon_K}(\mathfrak{h}_0),{\bf
  d}_{\varepsilon_K}(\Pi))$ satisfait {\`a} la condition $({\tt *})$.
\end{lemme}

\begin{proof} La sous-alg{\`e}bre ${\bf
  d}_{\varepsilon_K}(\mathfrak{h}_0)$ est une sous-alg{\`e}bre de Cartan de
$\mathfrak{g}_0$ et l'expression de ${\bf d}_{\varepsilon_K}$ montre
qu'elle est stable par $\theta$. Puisque les racines $\varepsilon_K$
et $\varepsilon_L$ sont fortement orthogonales, pour tout $L$ dans
$\mathcal{K}(\Pi)$ diff{\'e}rent de $K$, il est clair que
  ${\bf d}_{\varepsilon_K}(\varepsilon_L)=\varepsilon_L$, pour tout $L$ dans $\mathcal{K}(\Pi)$ diff{\'e}rent de $K$. Par
  ailleurs, d'apr{\`e}s la relation (\ref{1}), ${\bf d}_{\varepsilon_K}(H_{\varepsilon_K})$
appartient {\`a} $i \mathbb{R} (X_{\varepsilon_K} +\theta
 X_{\varepsilon_K})$, donc $H_{{\bf d}_{\varepsilon_K}(\varepsilon_K)}={\bf
 d}_{\varepsilon_K}(H_{\varepsilon_K})$ appartient {\`a} $\mathfrak{k}$. Par suite,
${\bf d}_{\varepsilon_K}(\varepsilon_K)$ est une racine imaginaire de
${\bf d}_{\varepsilon_K}(\Delta)$, d'apr{\`e}s le lemme \ref{struct1}
(iii). Les expressions de $\mathcal{K}^{''}({\bf d}_{\varepsilon_K}(\Pi))$ et
$\mathcal{K}^{'}({\bf d}_{\varepsilon_K}(\Pi))$ sont alors claires.\\
 
Compte tenu des relations pr{\'e}c{\'e}dentes, il est clair que si $(\mathfrak{h}_0,\Pi)$ a la propri{\'e}t{\'e} $(P)$, alors
la sous-alg{\`e}bre $({\bf d}_{\varepsilon_K}(\mathfrak{h}_0),{\bf
  d}_{\varepsilon_K}(\Pi))$ a la propri{\'e}t{\'e} $(P)$. Supposons que
$(\mathfrak{h}_0,\Pi)$ satisfait {\`a} la condition $({\tt *})$. D'apr{\`e}s
le lemme \ref{Kostant_Cartan2} (i), 
les racines $\varepsilon_K$ et $\alpha$ sont fortement
orthogonales, pour tout $\alpha$ dans $\Delta_{+}^{'}$. Par suite, on a\,:
$${\bf
  d}_{\varepsilon_K}(\Delta_{+}^{'})=\Delta_{+}^{'} \cup {\bf
  d}_{\varepsilon_K}(\varepsilon_K).$$

De plus, comme ${\bf
  d}_{\varepsilon_K}(\alpha)=\alpha$, pour tout $\alpha \in
\Delta_{+}^{'}$, il est clair que $({\bf d}_{\varepsilon_K}(\mathfrak{h}_0),{\bf
  d}_{\varepsilon_K}(\Pi))$ satisfait {\`a} la condition $({\tt
  *})$, puisque $(\mathfrak{h}_0,\Pi)$
satisfait {\`a} la condition $({\tt *})$.
\end{proof}

\begin{proposition}\label{calcul_k} On as\,:
$$\dim (\widehat{\mathfrak{h}} \cap \mathfrak{p}) - \# \mathcal{K}_{{\rm r\acute{e}el}}(\widehat{\Pi}) =\dim \widehat{\mathfrak{a}} - 
\# \mathcal{K}_{{\rm r\acute{e}el}}(\widehat{\Pi}) \geq {\rm rg} \; \mathfrak{g} - {\rm rg} \;
\mathfrak{k} .$$
Si de plus la condition $({\tt *})$ est
satisfaite,  alors on a les {\'e}galit{\'e}s suivantes\,:
$$\dim (\widehat{\mathfrak{h}} \cap \mathfrak{p}) - \# \mathcal{K}_{{\rm r\acute{e}el}}(\widehat{\Pi}) =\dim \widehat{\mathfrak{a}} - 
\# \mathcal{K}_{{\rm r\acute{e}el}}(\widehat{\Pi}) = {\rm rg} \; \mathfrak{g} - {\rm rg} \;
\mathfrak{k} .$$
\end{proposition}

\begin{proof} Notons tout d'abord que si $\mathfrak{h}_0$ est une
  sous-alg{\`e}bre de Cartan $\theta$-stable, on a\,:
$$\dim (\mathfrak{h} \cap \mathfrak{p}) \geq {\rm rg} \; \mathfrak{g} - {\rm rg} \;
\mathfrak{k},$$
car $\dim \mathfrak{t} \leq {\rm rg} \;
\mathfrak{k}$.\\
\\
\underline{$1^{\textrm{{\`e}re}}$ {\'e}tape}\,: si $\mathfrak{h}_0$ est une sous-alg{\`e}bre de Cartan de $\mathfrak{g}_0$
stable par $\theta$ et maximalement
compacte, alors $\dim (\mathfrak{h} \cap \mathfrak{p})={\rm rg} \; \mathfrak{g} - {\rm rg} \;
\mathfrak{k}$. D'apr{\`e}s le lemme \ref{max_compacte}, l'ensemble $ \mathcal{K}_{{\rm r\acute{e}el}}(\Pi)$ est 
vide et la relation,
$$\dim (\mathfrak{h} \cap \mathfrak{p}) - \# \mathcal{K}_{{\rm r\acute{e}el}}(\Pi) ={\rm rg} \; \mathfrak{g} - {\rm rg} \;
\mathfrak{k},$$
est satisfaite.\\
\\
\underline{$2^{\textrm{{\`e}me}}$ {\'e}tape}\,: soit $\mathfrak{h}_0$ une sous-alg{\`e}bre de Cartan de $\mathfrak{g}_0$
stable par $\theta$ telle que $(\mathfrak{h}_0,\Pi)$ a la
propri{\'e}t{\'e} $(P)$ et qui n'est pas maximalement
compacte. Si l'ensemble $\mathcal{K}_{{\rm
      r\acute{e}el}}(\Pi)$ est vide, alors on a l'in{\'e}galit{\'e}\,:
$$\dim (\mathfrak{h} \cap \mathfrak{p}) - \#
\mathcal{K}_{{\rm r\acute{e}el}}(\Pi) =\dim
(\mathfrak{h} \cap \mathfrak{p}) \geq {\rm rg} \; \mathfrak{g} - {\rm rg} \;
\mathfrak{k}.$$
Sinon, soit  $K$ dans $\mathcal{K}_{{\rm
      r\acute{e}el}}(\Pi)$. D'apr{\`e}s le lemme \ref{cayley}, la
  sous-alg{\`e}bre ${\bf
  d}_{\varepsilon_K}(\mathfrak{h}_0)$ est une sous-alg{\`e}bre de Cartan
stable par $\theta$ et on a la relation,
$$\# \mathcal{K}_{{\rm r\acute{e}el}}({\bf d}_{\varepsilon_K}(\Pi))=\#
\mathcal{K}_{{\rm r\acute{e}el}}(\Pi) - 1 .$$
Par ailleurs, d'apr{\`e}s la relation (\ref{rel_dim}), on a\,:
$$\dim ({\bf
  d}_{\varepsilon_K}(\mathfrak{h}) \cap \mathfrak{p})=\dim
  (\mathfrak{h} \cap \mathfrak{p})-1.$$
On en d{\'e}duit l'{\'e}galit{\'e},
$$\dim
  (\mathfrak{h} \cap \mathfrak{p})-\#
\mathcal{K}_{{\rm r\acute{e}el}}(\Pi)=\dim ({\bf
  d}_{\varepsilon_K}(\mathfrak{h}) \cap \mathfrak{p})-\#
  \mathcal{K}_{{\rm r\acute{e}el}}({\bf d}_{\varepsilon_K}(\Pi)).$$
Par ailleurs, $({\bf
  d}_{\varepsilon_K}(\mathfrak{h}_0),{\bf d}_{\varepsilon_K}(\Pi))$ poss{\`e}de la propri{\'e}t{\'e} $(P)$, d'apr{\`e}s le lemme
\ref{cayley}. De plus, toujours d'apr{\`e}s le lemme
\ref{cayley}, si $(\mathfrak{h}_0,\Pi)$ v{\'e}rifie la condition
$({\tt *})$, alors  $({\bf
  d}_{\varepsilon_K}(\mathfrak{h}_0),{\bf d}_{\varepsilon_K}(\Pi))$
la satisfait aussi.\\

Prouvons alors la proposition. Si $\widehat{\mathfrak{h}}_0$ est
maximalement compacte, on a l'{\'e}galit{\'e}\,:
$$\dim (\widehat{\mathfrak{h}} \cap \mathfrak{p}) - \#
\mathcal{K}_{{\rm r\acute{e}el}}(\widehat{\Pi}) ={\rm rg} \; \mathfrak{g} - {\rm rg} \;
\mathfrak{k},$$
d'apr{\`e}s la premi{\`e}re {\'e}tape. Sinon, on applique la deuxi{\`e}me {\'e}tape 
autant de fois n{\'e}cessaires pour obtenir une sous-alg{\`e}bre de Cartan
maximalement compacte. On applique la premi{\`e}re {\'e}tape {\`a} la sous-alg{\`e}bre de Cartan ainsi
obtenue. D'apr{\`e}s ces deux {\'e}tapes, on obtient l'in{\'e}galit{\'e} souhait{\'e}e.\\

Si de plus la condition $({\tt *})$ est
satisfaite, alors on obtient une {\'e}galit{\'e} car {\`a}
chaque {\'e}tape la nouvelle sous-alg{\`e}bre de Cartan $\theta$-stable
obtenue satisfiait la propri{\'e}t{\'e} $(P)$ et v{\'e}rifie la condition
$({\tt *})$. L'ensemble
$\mathcal{K}_{{\rm r\acute{e}el}}(\Pi)$ est alors non vide tant que
cette sous-alg{\`e}bre n'est pas maximalement compacte, d'apr{\`e}s le
lemme \ref{Kostant_Cartan2} (iii).
\end{proof}

\section{Formes lin{\'e}aires stables et indice de $\mathfrak{b}$}\label{part_stable}

\subsection{} Si $\mathfrak{q}$ est une alg{\`e}bre de Lie complexe et si $\varphi$ est une forme lin{\'e}aire sur $\mathfrak{q}$, on
d{\'e}signe par $\mathfrak{q}_{\varphi}$ l'ensemble des $s$ de $\mathfrak{q}$ tels que
$\varphi([\mathfrak{q},s])=0$. Autrement dit $\mathfrak{q}_{\varphi}=\{s \in \mathfrak{q} \; | \; ({\rm
  ad}^{*} s) \cdot \varphi=0 \}$, o{\`u} ${\rm ad}^{*} : \mathfrak{q} \rightarrow \mathfrak{gl}(\mathfrak{q}^{*})$ est la
repr{\'e}sentation coadjointe de $\mathfrak{q}$. On rappelle que l'{\it indice de
$\mathfrak{q}$}, not{\'e} ${\rm ind}~\mathfrak{q}$, est d{\'e}fini par\,:
$${\rm ind}~\mathfrak{q}= \min\limits_{\varphi \in \mathfrak{q}^{*}} \dim \mathfrak{q}_{\varphi} \;
.$$
On dit que l'{\'e}l{\'e}ment $\varphi$ de $\mathfrak{q}^*$ est {\it r{\'e}gulier} si $\dim
\mathfrak{q}_{\varphi}={\rm ind}~\mathfrak{q}$. L'ensemble des {\'e}l{\'e}ments
r{\'e}guliers de $\mathfrak{q}^{*}$ est un ouvert non vide de $\mathfrak{q}^*$. 

La notion de formes lin{\'e}aires stables est introduite dans \cite{Kosmann}. Rappelons qu'un {\'e}l{\'e}ment $\varphi$ de $\mathfrak{q}^*$
est dit {\it stable} s'il existe un voisinage $V$ de $\varphi$ dans $\mathfrak{q}^*$
tel que, pour tout $\psi$ de $V$, $\mathfrak{q}_{\varphi}$ et $\mathfrak{q}_{\psi}$ soient
conjugu{\'e}s par un {\'e}l{\'e}ment du groupe adjoint alg{\'e}brique de $\mathfrak{q}$. En particulier,
si $\varphi$ est une forme lin{\'e}aire stable, alors c'est un {\'e}l{\'e}ment
r{\'e}gulier de $\mathfrak{q}^*$. Lorsque $\mathfrak{q}$ poss{\`e}de une forme lin{\'e}aire stable,
l'indice de $\mathfrak{q}$ est donn{\'e} par la dimension du stabilisateur de cette
forme. En g{\'e}n{\'e}ral, $\mathfrak{q}$ ne poss{\`e}de pas de forme lin{\'e}aire stable; on
trouve des exemples d'alg{\`e}bres de Lie ne poss{\'e}dant pas de forme
lin{\'e}aire stable dans \cite{Kosmann} ou dans \cite{TauvelYu}. La
proposition qui suit est d{\'e}montr{\'e}e dans \cite{TauvelYu}, Th{\'e}or{\`e}me 1.7.\\

\begin{proposition}\label{stable} Soit $\mathfrak{q}$ une alg{\`e}bre de Lie, et
  soit $\varphi$
  un {\'e}l{\'e}ment de $\mathfrak{q}^*$. Les
  conditions suivantes sont {\'e}quivalentes\,: \begin{description}
\item[{\rm (i)}] On a la relation\,: $[\mathfrak{q},\mathfrak{q}_{\varphi}] \cap \mathfrak{q}_{\varphi}=\{0\}$,
\item[{\rm (ii)}] La forme lin{\'e}aire $\varphi$ est stable.\\

\end{description}
\end{proposition} 

\subsection{} On reprend les notations des parties pr{\'e}c{\'e}dentes. Rappelons que le
couple $(\widehat{\mathfrak{h}}_0,\widehat{\Pi})$ a la propri{\'e}t{\'e}
$(P)$, d'apr{\`e}s la proposition \ref{Kostant_Cartan1}. On va construire une forme lin{\'e}aire r{\'e}guli{\`e}re sur
$\mathfrak{b}$. Posons
\begin{eqnarray*}
u & = & \sum\limits_{K \in \mathcal{K}^{''}(\widehat{\Pi}) }
X_{-\varepsilon_{K}} = \displaystyle{\frac{1}{2}} \sum\limits_{K \in \mathcal{K}^{''}(\widehat{\Pi})} (X_{-\varepsilon_{K}}
+X_{-\varepsilon_{\theta K}}),\\
& = & \sum\limits_{K \in \mathcal{K}_{{\rm r\acute{e}el}}(\widehat{\Pi})} X_{-\varepsilon_{K}}
+\sum\limits_{K \in \mathcal{K}_{{\rm comp}}^+(\widehat{\Pi})} (X_{-\varepsilon_{K}}
+X_{-\varepsilon_{\theta K}}) .
\end{eqnarray*}

L'{\'e}l{\'e}ment $u$
  appartient {\`a} $\mathfrak{n}_-$. Pour $K$ dans $ \mathcal{K}^{''}(\widehat{\Pi})$, on pose\,:
$$\Gamma_{1}^{K}=\{\alpha \in \Gamma_{0}^{K} \ | \ \varepsilon_K-
\alpha \in \Delta_{+}^{'}\}.$$

\begin{lemme}\label{gamma1} \begin{itemize} \item[{\rm (i)}] On a\,:
    $\mathcal{K}(\widehat{\Pi}^{'})=\mathcal{K}^{'}(\widehat{\Pi})
    \iff \forall K \in \mathcal{K}^{''}(\widehat{\Pi}), \
  \Gamma_{1}^{K}=\emptyset.$

\item[{\rm (ii)}] Soit $K$ dans $\mathcal{K}^{''}(\widehat{\Pi})$ tel
  que $\Gamma_{1}^{K} \not=\emptyset$.  Alors $\theta K
  \subset_{\atop\hspace{-0.3cm}\not=} K$. En particulier, $K$ appartient {\`a} $\mathcal{K}_{{\rm comp}}(\widehat{\Pi})$.
\end{itemize}
\end{lemme}

\begin{proof} (i)  Supposons
  $\mathcal{K}(\widehat{\Pi}^{'})=\mathcal{K}^{'}(\widehat{\Pi})$. Alors il r{\'e}sulte du lemme
  \ref{Kostant_Cartan2} (i) que, pour tout $K$ dans
  $\mathcal{K}^{''}(\widehat{\Pi})$, l'ensemble $\Gamma_{1}^{K}$ est
  vide. R{\'e}ciproquement, si $\mathcal{K}(\widehat{\Pi}^{'}) \not=\mathcal{K}^{'}(\widehat{\Pi})$,
  alors il existe $L$ dans $\mathcal{K}(\widehat{\Pi}^{'})$ tel que
  $\varepsilon_L \not\in \mathcal{K}^{'}(\widehat{\Pi})$. Posons
  $\beta=\varepsilon_L$. Il existe alors $K$ dans $\mathcal{K}^{''}$,
  tel que $\beta \in  \Gamma_{0}^{K}$. Il est alors clair que $\alpha=\varepsilon_K-
\beta$ appartient $\Gamma_{1}^{K}$.\\
\\
(ii) Soit $K \in \mathcal{K}^{''}(\widehat{\Pi})$ tel que
$\Gamma_{1}^{K}$ est non vide et soit $\alpha \in
\Gamma_{1}^{K}$. Alors il existe $\beta \in \Gamma_{0}^{K}$ tel que
$$\alpha +\beta =\varepsilon_K.$$
La relation
$$- \theta \alpha=\varepsilon_{\theta K} + \beta$$
assure alors la relation $\theta K
  \subset K$, d'apr{\`e}s le lemme
  \ref{kostant2} (iv). L'{\'e}galit{\'e} $\theta K=K$ est impossible car
  $\varepsilon_{ K} + \beta$ n'est pas une racine, d'o{\`u} (ii).
\end{proof}

\begin{lemme}\label{equivalence} Soit $x$ un {\'e}l{\'e}ment de $\mathfrak{b}$. Les conditions suivantes sont {\'e}quivalentes\,:\begin{description}

\item[{\rm (i)}] L'{\'e}l{\'e}ment $x$ s'{\'e}crit sous la forme 
$$x=h+ \sum\limits_{K \in \mathcal{K}_{{\rm comp}}^+(\widehat{\Pi})}\left( a_K (X_{\varepsilon_{K}} -X_{\varepsilon_{\theta K}}) +
\sum\limits_{\alpha \in \Gamma_{1}^{K}} a_{\alpha} X_{\alpha} \right),$$
avec $h$ dans $\widehat{\mathfrak{a}}$ tel que $\varepsilon_{K}(h)=0$, pour tout $K
\in \mathcal{K}^{''}(\widehat{\Pi})$, et $a_K, a_{\alpha} \in \mathbb{C}$,  pour tout $K
\in  \mathcal{K}_{{\rm comp}}^+(\widehat{\Pi})$ et $\alpha \in
\Gamma_{1}^{K}$,

\item[{\rm (ii)}] Le crochet $[x,u]$ appartient au sous-espace $\mathfrak{n} \oplus \mathfrak{m}$.\\ 

\end{description}

\end{lemme}

\begin{remarque}\label{borel} Notons que si $\mathfrak{m}$ est nul,
  $\mathfrak{b}$ est une sous-alg{\`e}bre de Borel de $\mathfrak{g}$ et ce lemme n'est rien d'autre
que le lemme 2.5 de \cite{TauvelYu}. La d{\'e}monstration qui suit
reprend d'ailleurs pour une large part celle de \cite{TauvelYu}.
\end{remarque}

\begin{proof} ${\rm (i)} \Rightarrow {\rm (ii)}$\,: si ${\rm (i)}$ est v{\'e}rifi{\'e},
$x$ s'{\'e}crit 
$$x=h+ \sum\limits_{K \in \mathcal{K}_{{\rm comp}}^{+}(\widehat{\Pi})}
\left( a_K (X_{\varepsilon_{K}} -X_{\varepsilon_{\theta K}}) +
\sum\limits_{\alpha \in \Gamma_{1}^{K}} a_{\alpha} X_{\alpha} \right),$$
avec $h$ dans $\widehat{\mathfrak{a}}$ tel que $\varepsilon_{K}(h)=0$, pour tout $K
\in \mathcal{K}^{''}(\widehat{\Pi})$ et $a_K, a_{\alpha} \in \mathbb{C}$,  pour tout $K
\in  \mathcal{K}_{{\rm comp}}^+(\widehat{\Pi})$ et $\alpha \in
\Gamma_{1}^{K}$. Posons $x'=h+ \sum\limits_{K \in \mathcal{K}_{{\rm
      comp}}^{+}(\widehat{\Pi}) }
 a_K (X_{\varepsilon_{K}} -X_{\varepsilon_{\theta K}})$ et $x''=\sum\limits_{K \in \mathcal{K}_{{\rm comp}}^{+}(\widehat{\Pi})}
\sum\limits_{\alpha \in \Gamma_{1}^{K} } a_{\alpha} X_{\alpha}$ de
sorte que $x=x'+x''$. On a 
$$[x'',u]=\displaystyle{\frac{1}{2}} \sum\limits_{K \in
  \mathcal{K}^{''}(\widehat{\Pi})} \sum\limits_{L \in \mathcal{K}_{{\rm comp}}^{+}(\widehat{\Pi})}
\sum\limits_{\alpha \in \Gamma_{1}^{L}} a_{\alpha }
[X_{\alpha},X_{-\varepsilon_{K}}].$$ Soit $L,K$ dans $\mathcal{K}^{''}(\widehat{\Pi})$ et $\alpha \in \Gamma_{1}^{L}$. Si
$K \not= L$, alors $[X_{\alpha},X_{-\varepsilon_{K}}]$ appartient au
sous-espace $\sum\limits_{\beta \in \Delta_+} \mathfrak{g}^{\beta}$
qui est contenu dans $\mathfrak{n} \oplus \mathfrak{m}$. Si $K=L$,
alors l'{\'e}l{\'e}ment 
$[X_{\alpha},X_{-\varepsilon_{L}}]$ appartient {\`a} $\mathfrak{m}$, par
d{\'e}finition de l'ensemble $\Gamma_{1}^{L}$. Par suite\,: 
$[x'',u] \in \mathfrak{n} \oplus \mathfrak{m}$.\\

Par ailleurs, on a
$$\begin{array}{rccl}
[x',u] = & \displaystyle{\frac{1}{2}} \sum\limits_{K \in
  \mathcal{K}^{''}(\widehat{\Pi})} ( & -\varepsilon_{K}(h)X_{-\varepsilon_{K}}-\varepsilon_{\theta
  K}(h)X_{-\varepsilon_{\theta K}} & + a_K(H_{\varepsilon_{K}}-H_{\varepsilon_{\theta K}})),\\
= & 
\displaystyle{\frac{1}{2}} \sum\limits_{K \in \mathcal{K}^{''}(\widehat{\Pi})} ( & -\varepsilon_{K}(h)X_{-\varepsilon_{K}} +\theta \varepsilon_{K}(h)X_{-\varepsilon_{\theta K}} & + a_K(H_{\varepsilon_{K}}+\theta
  H_{\varepsilon_{K}})),\\
= & \displaystyle{\frac{1}{2}} \sum\limits_{K \in \mathcal{K}^{''}(\widehat{\Pi})} (
& -\varepsilon_{K}(h)X_{-\varepsilon_{K}} -\varepsilon_{K}(h)X_{-\varepsilon_{\theta K}}  & + a_K(H_{\varepsilon_{K}}+\theta
  H_{\varepsilon_{K}} )),\\
= &\displaystyle{\frac{1}{2}} \sum\limits_{K \in \mathcal{K}^{''}(\widehat{\Pi})} ( &
0  & +a_K(H_{\varepsilon_{K}}+\theta
  H_{\varepsilon_{K}})) .
\end{array}$$
On en d{\'e}duit que $[x',u]$ appartient {\`a} $\widehat{\mathfrak{t}} \subset
\mathfrak{m}$. Finalement, $[x,u]$ appartient {\`a} $\mathfrak{n} \oplus \mathfrak{m}$.\\
\\
${\rm (ii)} \Rightarrow {\rm (i)}$\,: l'{\'e}l{\'e}ment $x$ appartient {\`a}
$\mathfrak{b}$. Puisque $\mathfrak{b}=\widehat{\mathfrak{a}} \oplus
\sum\limits_{\alpha \in \widehat{\Delta}_{+}^{''}}
\mathfrak{g}^{\alpha}=\widehat{\mathfrak{a}} \oplus (\sum\limits_{K
  \in \mathcal{K}^{''}(\widehat{\Pi})} \sum\limits_{\alpha \in
  \Gamma^{K}} \mathfrak{g}^{\alpha})$, l'{\'e}l{\'e}ment $x$ s'{\'e}crit sous la forme\,: 
\begin{eqnarray*}
x   = h &+ &\sum\limits_{K \in \mathcal{K}_{{\rm r\acute{e}el}}(\widehat{\Pi})} a_{\varepsilon_K}
X_{\varepsilon_{K}} \\
&  + &\sum\limits_{K \in \mathcal{K}_{{\rm comp}}^{+}(\widehat{\Pi})}
(a_{\varepsilon_K} X_{\varepsilon_{K}}+ a_{\varepsilon_{\theta K}} X_{\varepsilon_{\theta K}}) 
+ \sum\limits_{K \in \mathcal{K}^{''}(\widehat{\Pi})} \sum\limits_{\alpha \in \Gamma_{0}^{K}} a_{\alpha}
X_{\alpha},
\end{eqnarray*}
avec $h$ dans $\widehat{\mathfrak{a}}$ et $a_{\alpha}$ dans $\mathbb{C}$, pour 
 tout $\alpha$ de $\Delta_{+}^{''}$. On en d{\'e}duit que le crochet
 $[x,u]$ s'{\'e}crit sous la forme\,:
$$[x,u] =  
\sum\limits_{K \in \mathcal{K}_{{\rm r\acute{e}el}}(\widehat{\Pi})} a_{\varepsilon_K}  H_{\varepsilon_K} 
+\sum\limits_{K \in \mathcal{K}_{{\rm comp}}^{+}(\widehat{\Pi})}(a_{\varepsilon_K} 
H_{\varepsilon_K} +  a_{\varepsilon_{\theta K}}  H_{\varepsilon_{\theta K}})
+ Y,$$
o{\`u} $Y$ est un {\'e}l{\'e}ment du sous-espace $\mathfrak{m} \oplus \mathfrak{n}
\oplus \mathfrak{n}_-$. Pour $K$ dans $\mathcal{K}_{{\rm
    r\acute{e}el}}(\widehat{\Pi})$, la racine $\varepsilon_K$ est r{\'e}elle
donc l'{\'e}l{\'e}ment $H_{\varepsilon_K}$
appartient {\`a} $\widehat{\mathfrak{a}}$. Pour $K$ dans $\mathcal{K}_{{\rm comp}}(\widehat{\Pi})$, la projection de $H_{\varepsilon_K}$ sur
$\widehat{\mathfrak{a}}$ selon
la d{\'e}composition $\widehat{\mathfrak{a}} \oplus \widehat{\mathfrak{t}}$ est
$\displaystyle{\frac{H_{\varepsilon_K}-\theta H_{\varepsilon_K}}{2}}$. Par
suite, de la
relation $[x,u] \in \mathfrak{m} \oplus \mathfrak{n}$, on tire la
relation\,:
$$ \sum\limits_{K \in \mathcal{K}_{{\rm r\acute{e}el}}(\widehat{\Pi})} a_{\varepsilon_K}  H_{\varepsilon_K} 
+\sum\limits_{K \in \mathcal{K}_{{\rm comp}}^{+}(\widehat{\Pi})}(a_{\varepsilon_K} 
\displaystyle{\frac{H_{\varepsilon_K}-\theta H_{\varepsilon_K}}{2}}  +  a_{\varepsilon_{\theta K}}
\displaystyle{\frac{H_{\varepsilon_{\theta K}}
-\theta H_{\varepsilon_{\theta K}}}{2}})=0 . $$ 
Puisque les racines $\varepsilon_K$ sont deux {\`a} deux fortement
orthogonales, on a $\varepsilon_K(H_{\varepsilon_L})=0$, pour $L\not=K$ et il vient
$a_{\varepsilon_K}=0$, pour tout $K$ de $ \mathcal{K}_{{\rm
    r\acute{e}el}}(\widehat{\Pi})$. En utilisant de plus les
relations $\theta H_{\varepsilon_K}=-H_{\varepsilon_{\theta K}}$ et 
$\theta H_{\varepsilon_{\theta K}}=-H_{\varepsilon_K}$, on obtient $a_{\varepsilon_K}+a_{\varepsilon_{\theta K}}=0$, pour tout $K$ de
$\mathcal{K}_{{\rm comp}}^{+}(\widehat{\Pi})$. On en d{\'e}duit que $x$ s'{\'e}crit\,:
$$x=h
+\sum\limits_{K \in \mathcal{K}_{{\rm comp}}^{+}(\widehat{\Pi})}
a_{\varepsilon_K} (X_{\varepsilon_{K}}-X_{\varepsilon_{\theta K}}) 
+ \sum\limits_{K \in \mathcal{K}^{''}(\widehat{\Pi})} \sum\limits_{\alpha \in \Gamma_{0}^{K}} a_{\alpha}
X_{\alpha} .$$

Soit $K$ un {\'e}l{\'e}ment de $\mathcal{K}^{''}(\widehat{\Pi})$ et soit
$\alpha$ dans $\Gamma_{0}^{K} \setminus \Gamma_{1}^{K}$. Alors la
racine $\beta=\varepsilon_{K}-\alpha$ de
$\Gamma_0^{K}$ appartient {\`a} $\widehat{\Delta}_{+}^{''}$. On a
\begin{eqnarray*}
[x,u] & & =\displaystyle{\frac{1}{2}} (\lambda a_{\alpha} X_{-\beta} 
 + a_{\alpha} 
\sum\limits_{K' \in \mathcal{K}^{''}(\widehat{\Pi})\setminus\{K\}}
  [X_{\alpha},X_{-\varepsilon_{K'}}]\\
& &+   \sum\limits_{K' \in \mathcal{K}^{''}(\widehat{\Pi}), \gamma \in \Delta_{+}^{''}
  \setminus\{ \alpha \}} a_{\gamma} [X_{\gamma},X_{-\varepsilon_{K'}}])
  -  \sum\limits_{K' \in  \mathcal{K}^{''}(\widehat{\Pi})} \varepsilon_{K'}(h)
X_{-\varepsilon_{K'}}\\
& &+ \sum\limits_{K \in \mathcal{K}_{{\rm comp}}^{+}(\widehat{\Pi})} a_{\varepsilon_K} 
(\displaystyle{\frac{H_{\varepsilon_K}+\theta H_{\varepsilon_K}}{2}}),
\end{eqnarray*} 
o{\`u} $\lambda$ est un scalaire non nul.\\

Supposons par
l'absurde $a_{\alpha} \not=0$. On a $\beta \not= \varepsilon_{K^{'}}$,
pour $K^{'}$ dans $\mathcal{K}^{''}(\widehat{\Pi})$, car $\langle
\beta, {\varepsilon_K}^{\vee} \rangle=1$ et $\langle \varepsilon_{K^{'}},
{\varepsilon_K}^{\vee} \rangle \in \{0,2\}$. Puisque la racine $\beta$ appartient {\`a} $\widehat{\Delta}_{+}^{''}$, l'{\'e}l{\'e}ment
 $\displaystyle{\frac{1}{2}} \lambda a_{\alpha} X_{-\beta}$ qui
intervient dans l'expression pr{\'e}c{\'e}dente de $[x,u]$ est un {\'e}l{\'e}ment non
nul de $\mathfrak{n}_{-}$. Comme $[x,u]$ appartient {\`a} $\mathfrak{n}
\oplus \mathfrak{m}$, il existe n{\'e}cessairement $K'$ dans
$\mathcal{K}^{''}(\widehat{\Pi})$ et $\gamma$ dans $\widehat{\Delta}_{+}^{''}
  \setminus\{ \alpha \}$ tels que
  $\beta=\varepsilon_{K'}-\gamma$. On a\,: $K\not= K'$, car $\gamma \not=
  \alpha$. Soit $K''$ dans $\mathcal{K}^{''}(\widehat{\Pi})$ tel que
  $\gamma$ appartient {\`a} $\Gamma^{K^{''}}$. Or d'apr{\`e}s le lemme
  \ref{kostant2}, (iii) et (iv), la relation
  $\beta+\gamma=\varepsilon_{K'}$ entra{\^\i}ne que $K=K^{''}$, puis que
  $K^{'}=K$, d'o{\`u} la contradiction.

On en d{\'e}duit que $a_{\alpha}$ est nul, pour tout  $K$ dans
$\mathcal{K}^{''}(\widehat{\Pi})$ et tout $\alpha$ dans
$\Gamma_{0}^{K} \setminus \Gamma_{1}^{K}$. Il r{\'e}sulte alors du lemme
\ref{gamma1} (ii) que
$x$ s'{\'e}crit sous la forme
$$x=h+ \sum\limits_{K \in \mathcal{K}_{{\rm comp}}^+(\widehat{\Pi})}\left( a_K (X_{\varepsilon_{K}} -X_{\varepsilon_{\theta K}}) +
\sum\limits_{\alpha \in \Gamma_{1}^{K}} a_{\alpha} X_{\alpha}
\right).$$
De l'{\'e}galit{\'e} 
$$[x,u]=- \sum\limits_{K' \in  \mathcal{K}^{''}(\widehat{\Pi})} \varepsilon_{K'}(h)
X_{-\varepsilon_{K'}}
+\sum\limits_{K \in \mathcal{K}_{{\rm comp}}^{+}(\widehat{\Pi})} a_K 
(\displaystyle{\frac{H_{\varepsilon_K}+\theta H_{\varepsilon_K}}{2}}),$$
on tire la relation $\varepsilon_{K'}(h)=0$, pour tout $K^{'}$ de
$\mathcal{K}^{''}(\widehat{\Pi})$, d'o{\`u} (i).
\end{proof}

On note $\kappa$ la forme de Killing de $\mathfrak{g}$. Pour $v$ un {\'e}l{\'e}ment de $\mathfrak{g}$, on note $\varphi_v$ la
forme lin{\'e}aire sur $\mathfrak{g}$ d{\'e}finie par\,:
$$\varphi_v(y)=\kappa(v,y),$$
pour $y$ dans $\mathfrak{g}$. Si $\mathfrak{r}$ est une sous-alg{\`e}bre
de $\mathfrak{g}$, $\mathfrak{r}_{\varphi_v}$ d{\'e}signe le stabilisateur dans
$\mathfrak{r}$ de la restriction de $\varphi_v$ {\`a} $\mathfrak{r}$. La
proposition suivante d{\'e}crit le stabilisateur dans $\mathfrak{b}$ de la
forme lin{\'e}aire $\varphi_u$.\\ 

\begin{proposition}\label{stab_u}
Le stabilisateur $\mathfrak{b}_{\varphi_u}$ de la restriction de $\varphi_u$
{\`a} $\mathfrak{b}$ est donn{\'e} par\,:
$$\mathfrak{b}_{\varphi_u}=( \bigcap\limits_{K \in
  \mathcal{K}(\widehat{\Pi})} \ker
  {\varepsilon_K}_{|_{\widehat{\mathfrak{a}}}} ) 
  \oplus \sum\limits_{\mathcal{K}_{{\rm
  comp}}^+(\widehat{\Pi})} \left( \mathbb{C}(X_{\varepsilon_{K}} - X_{\varepsilon_{\theta
  K}}) \oplus \sum\limits_{\alpha \in \Gamma_{1}^{K}} \mathbb{C}
X_{\alpha} \right)$$
et on a\,:
$$\dim \mathfrak{b}_{\varphi_u}=\dim \widehat{\mathfrak{a}}-\#
  \mathcal{K}_{{\rm r\acute{e}el}}(\widehat{\Pi}) + \sum\limits_{\mathcal{K}_{{\rm
  comp}}^+(\widehat{\Pi})} \#  \Gamma_{1}^{K}.$$
\end{proposition}

\begin{proof} On a 
$$\mathfrak{b}_{\varphi_u} =  \{x \in \mathfrak{b} \ | \ \kappa(u,[x,y])
= 0, \ \forall y \in \mathfrak{b} \}\\
=  \{x \in \mathfrak{b} \ | \ [x,u] \in {\mathfrak{b}}^{\perp} \},$$
o{\`u} ${\mathfrak{b}}^{\perp}$ d{\'e}signe l'orthogonal de $\mathfrak{b}$
dans $\mathfrak{g}$ pour la forme de Killing. La relation 
${\mathfrak{b}}^{\perp}=\mathfrak{n} \oplus \mathfrak{m}$ et le lemme
\ref{equivalence} entra{\^\i}nent que l'{\'e}l{\'e}ment $x$ appartient {\`a}
$\mathfrak{b}_{\varphi_u}$ si, et seulement si, il s'{\'e}crit sous la forme 
$$x=h+ \sum\limits_{K \in \mathcal{K}_{{\rm comp}}^{+}(\widehat{\Pi})}
\left( a_K (X_{\varepsilon_{K}} -X_{\varepsilon_{\theta K}}) +
\sum\limits_{\alpha \in \Gamma_{1}^{K}} a_{\alpha} X_{\alpha} \right),$$
avec $h$ dans $\widehat{\mathfrak{a}}$ tel que $\varepsilon_{K}(h)=0$, pour tout $K
\in \mathcal{K}^{''}(\widehat{\Pi})$ et $a_K, a_{\alpha} \in \mathbb{C}$,  pour tout $K
\in  \mathcal{K}_{{\rm comp}}^+(\widehat{\Pi})$ et $\alpha \in
\Gamma_{1}^{K}$

On dispose de l'{\'e}galit{\'e}\,: $\bigcap\limits_{K \in
  \mathcal{K}^{''}(\widehat{\Pi})} \ker {\varepsilon_K}_{|_{\widehat{\mathfrak{a}}}}=\bigcap\limits_{K \in
  \mathcal{K}(\widehat{\Pi})} \ker
  {\varepsilon_K}_{|_{\widehat{\mathfrak{a}}}}$. En effet, pour $K$ dans
$\mathcal{K}^{'}(\widehat{\Pi})$, la racine $\varepsilon_K$ est imaginaire donc s'annule sur
$\widehat{\mathfrak{a}}$. La premi{\`e}re assertion de la proposition est alors claire.\\  

Montrons\,: $\dim (\bigcap\limits_{K \in
  \mathcal{K}(\widehat{\Pi})} \ker
  {\varepsilon_K}_{|_{\widehat{\mathfrak{a}}}} )=\dim \widehat{\mathfrak{a}}-(\#
  \mathcal{K}_{{\rm r\acute{e}el}}(\widehat{\Pi})+\# \mathcal{K}_{{\rm
  comp}}^+(\widehat{\Pi}))$. Il suffit de montrer que la famille $\{
  {\varepsilon_K}_{|_{\widehat{\mathfrak{a}}}} \ | \ K \in  \mathcal{K}_{{\rm
  r\acute{e}el}}(\widehat{\Pi}) \cup \mathcal{K}_{{\rm
  comp}}^+(\widehat{\Pi}) \}$ forme une base du sous-espace de
  $\widehat{\mathfrak{a}}^{*}$ engendr{\'e} par les {\'e}l{\'e}ments
  ${\varepsilon_K}_{|_{\widehat{\mathfrak{a}}}}$, pour $K$ dans $
  \mathcal{K}(\widehat{\Pi})$. Si $K$ appartient {\`a}
  $\mathcal{K}^{'}(\widehat{\Pi})$, la racine $\varepsilon_K$ est
  imaginaire donc s'annule sur $\widehat{\mathfrak{a}}$. Si $K$ appartient {\`a}
  $\mathcal{K}_{{\rm comp}}^+(\widehat{\Pi})$, alors pour tout $H$
  dans 
  $\widehat{\mathfrak{a}}$, on a\,: $\varepsilon_{\theta K}(H)=-\theta \varepsilon_K
  (H)=  \varepsilon_K (-\theta H)= \varepsilon_K (H)$, d'o{\`u} ${\varepsilon_{\theta
  K}}_{|_{\widehat{\mathfrak{a}}}}={\varepsilon_K}_{|_{\widehat{\mathfrak{a}}}}$. On en d{\'e}duit
  que la famille pr{\'e}c{\'e}dente est g{\'e}n{\'e}ratrice. Montrons qu'elle est
  libre\,: soit 
$$\sum\limits_{ K \in \mathcal{K}_{{\rm r\acute{e}el}}(\widehat{\Pi})}
  a_K  {\varepsilon_K}_{|_{\widehat{\mathfrak{a}}}} + \sum\limits_{ K \in \mathcal{K}_{{\rm
  comp}}^+(\widehat{\Pi}) } a_K {\varepsilon_K}_{|_{\widehat{\mathfrak{a}}}}=0,$$
une combinaison lin{\'e}aire nulle dans
  $\widehat{\mathfrak{a}}^*$. Puisque la famille $\{ \varepsilon_K \ | \
  K \in  \mathcal{K}(\widehat{\Pi}) \}$ est un ensemble de racines
  deux {\`a} deux fortement orthogonales, l'{\'e}valuation du membre de gauche dans l'expression pr{\'e}c{\'e}dente en l'{\'e}l{\'e}ment
  $\displaystyle{\frac{H_{\varepsilon_K}-\theta
  H_{\varepsilon_K}}{2}}=\displaystyle{\frac{H_{\varepsilon_K}+H_{\varepsilon_{\theta
  K}}}{2}}$ de $\widehat{\mathfrak{a}}$, donne\,: $a_K=0$, pour
  tout $K$ de $\mathcal{K}_{{\rm r\acute{e}el}}(\widehat{\Pi})$ et
  tout $K$ de $\mathcal{K}_{{\rm
  comp}}^+(\widehat{\Pi})$. On a obtenu\,:
\begin{eqnarray*}
\dim \mathfrak{b}_{\varphi_u} & = & (\dim \widehat{\mathfrak{a}}-\# \mathcal{K}_{{\rm r\acute{e}el}}(\widehat{\Pi}) -
 \# \mathcal{K}_{{\rm comp}}^+(\widehat{\Pi}))+\# \mathcal{K}_{{\rm
 comp}}^+(\widehat{\Pi}) +\sum\limits_{\mathcal{K}_{{\rm
  comp}}^+(\widehat{\Pi})} \#  \Gamma_{1}^{K}\\
& = & \dim \widehat{\mathfrak{a}}-\# \mathcal{K}_{{\rm
 r\acute{e}el}}(\widehat{\Pi})+\sum\limits_{\mathcal{K}_{{\rm
  comp}}^+(\widehat{\Pi})} \#  \Gamma_{1}^{K},
\end{eqnarray*}
d'o{\`u} la proposition.
\end{proof}

\begin{remarque}\label{rq_stab} L'expression de
  $\mathfrak{b}_{\varphi_u}$ obtenue dans la proposition pr{\'e}c{\'e}dente permet d'obtenir une condition n{\'e}cessaire et
  suffisante pour que la restriction de $\varphi_u$ {\`a} $\mathfrak{b}$ soit stable. Pour $a$ et $b$ dans $\mathbb{C}$ et $K$ dans
  $\mathcal{K}_{{\rm comp}}^+(\widehat{\Pi})$, on a\,:
\begin{eqnarray}\label{crochet}
[H,a X_{\varepsilon_K} +b X_{\varepsilon_{\theta K}} ] 
=  \varepsilon_K(H) (a
X_{\varepsilon_K} +b X_{\varepsilon_{\theta K}}),
\end{eqnarray}
pour tout $H$ dans $\widehat{\mathfrak{a}}$. En particulier, les
{\'e}l{\'e}ments $X_{\varepsilon_{K}} - X_{\varepsilon_{\theta
  K}}$, pour $K$ dans $\mathcal{K}_{{\rm comp}}^+(\widehat{\Pi})$,
appartiennent {\`a} l'intersection $[\mathfrak{b},\mathfrak{b}_{\varphi_u}]
\cap \mathfrak{b}_{\varphi_u}$. De l'expression de $\mathfrak{b}_{\varphi_u}$
obtenue dans la proposition \ref{stab_u} et de la relation
$[\mathfrak{b},\mathfrak{b}]=\mathfrak{n}$, on tire alors l'{\'e}galit{\'e}\,:
$$[\mathfrak{b},\mathfrak{b}_{\varphi_u}] \cap \mathfrak{b}_{\varphi_u}= \sum\limits_{\mathcal{K}_{{\rm
  comp}}^+(\widehat{\Pi})} \mathbb{C}(X_{\varepsilon_{K}} - X_{\varepsilon_{\theta
  K}}) .$$
Il r{\'e}sulte de la proposition \ref{stable} que la restriction de
  $\varphi_u$ {\`a} $\mathfrak{b}$ n'est pas stable si l'ensemble $\mathcal{K}_{{\rm
  comp}}(\widehat{\Pi})$ n'est pas vide. R{\'e}ciproquement, il est clair
que si  l'ensemble $\mathcal{K}_{{\rm
  comp}}(\widehat{\Pi})$ est vide, alors $\varphi_u$ {\`a}
$\mathfrak{b}$-stable, d'apr{\`e}s la proposition \ref{stab_u}.
\end{remarque}

Posons
$$\mathfrak{r}=\sum\limits_{\alpha \in \widehat{\Delta}_{+}^{'}}
\mathfrak{g}^{\alpha}
 \textrm{ \  \  et \ \ } 
\mathfrak{r}_{-}=\sum\limits_{\alpha \in \widehat{\Delta}_{+}^{'}}
\mathfrak{g}^{-\alpha} .$$
Ainsi on a\,: $\mathfrak{m}=\mathfrak{r}_{-} \oplus \widehat{\mathfrak{t}} \oplus
\mathfrak{r}$. Posons aussi\,:
$$\widetilde{\mathfrak{b}}=\mathfrak{b} \oplus \widehat{\mathfrak{t}} \oplus
\mathfrak{r},$$
de sorte que $\widetilde{\mathfrak{b}}$ est une sous-alg{\`e}bre de Borel
de $\mathfrak{g}$. Posons enfin
$$\widetilde{u}=u + \sum\limits_{K \in \mathcal{K}^{'}(\widehat{\Pi})} X_{-\varepsilon_K} .$$
D'apr{\`e}s la partie 2 de \cite{TauvelYu} ou d'apr{\`e}s la remarque
  \ref{rq_stab} appliqu{\'e}e au cas $\mathfrak{m}=0$ (remarque
  \ref{borel}), la restriction de la forme $\varphi_{\widetilde{u}}$ {\`a}
  $\widetilde{\mathfrak{b}}$ est stable pour
  $\widetilde{\mathfrak{b}}$. Notons que $\varphi_u$ n'est rien d'autre que la restriction {\`a}
$\mathfrak{b}$ de la forme lin{\'e}aire $\varphi_{\widetilde{u}}$. On est
d{\'e}sormais en mesure de d{\'e}montrer la relation (\ref{relation2}) annonc{\'e}e en introduction\,:\\

\begin{theoreme}\label{formule_indice} On a\,:
$${\rm ind} \; \mathfrak{b} \geq {\rm rg} \; \mathfrak{g} - {\rm rg} \;
\mathfrak{k}.$$
De plus, l'{\'e}galit{\'e} 
$${\rm ind} \; \mathfrak{b}={\rm rg} \; \mathfrak{g} - {\rm rg} \;
\mathfrak{k}$$
a lieu si, et seulement si, la condition
$({\tt *})$ est satisfaite.
\end{theoreme}

\begin{proof} D'apr{\`e}s la proposition \ref{stab_u} et la proposition
  \ref{calcul_k}, on a $\dim \mathfrak{b}_{\varphi_u} \geq {\rm rg} \; \mathfrak{g} - {\rm rg} \;
\mathfrak{k}$. De plus, il r{\'e}sulte du lemme \ref{gamma1} (i) et de la
proposition \ref{calcul_k} que la relation $\dim \mathfrak{b}_{\varphi_u}={\rm rg} \; \mathfrak{g} - {\rm rg} \;
\mathfrak{k}$ a lieu si, et seulement si, la condition
$\mathcal{K}(\widehat{\Pi}^{'})=\mathcal{K}^{'}(\widehat{\Pi})$ est
remplie. Il suffit donc de prouver que la restriction de
$\varphi_u$ {\`a} $\mathfrak{b}$ est $\mathfrak{b}$-r{\'e}guli{\`e}re.\\

Soit $\widetilde{B}$ le groupe adjoint alg{\'e}brique
  de $\widetilde{\mathfrak{b}}$. Puisque $\mathfrak{n}$ est un id{\'e}al
  de $\widetilde{\mathfrak{b}}$ contenu dans le radical nilpotent
  $\mathfrak{n} \oplus \mathfrak{r}$ de $\widetilde{\mathfrak{b}}$, il r{\'e}sulte de la proposition
  40.6.3 de \cite{Tauvel} que l'orbite de la forme lin{\'e}aire 
  ${\varphi_{\widetilde{u}}}_{|_{\mathfrak{n}}}={\varphi_u}_{|_{\mathfrak{n}}}$
  de $\mathfrak{n}^*$ sous l'action de $\widetilde{B}$ est ouverte dans
  $\mathfrak{n}^*$. Le dual de $\mathfrak{n}$ s'identifie via
  la forme de Killing de $\mathfrak{g}$ au sous-espace
  $\mathfrak{n}_{-}$. On en d{\'e}duit que l'ensemble
$$\widetilde{V}=\{v \in \widehat{\mathfrak{a}} \oplus \mathfrak{n}_{-}
\ | \ (\varphi_{{\rm pr}_{ \mathfrak{n}_-}(v)})_{|_{\mathfrak{n}}}
\in  \widetilde{B} \cdot {\varphi_u}_{|_{\mathfrak{n}}} \}$$
est un ouvert non vide de $\widehat{\mathfrak{a}} \oplus
\mathfrak{n}_{-}$, o{\`u} ${\rm pr}_{ \mathfrak{n}_-}$ est la
  projection de $\widehat{\mathfrak{a}} \oplus  \mathfrak{n}_{-}$ sur
  $\mathfrak{n}_{-}$ parall{\`e}lement {\`a} $\widehat{\mathfrak{a}}$. Par
  ailleurs, le dual de $\mathfrak{b}$ s'identifie {\`a} $\widehat{\mathfrak{a}} \oplus \mathfrak{n}_{-}$ via la
forme de Killing de $\mathfrak{g}$. On en d{\'e}duit que l'ensemble
$$V=\{v \in \widehat{\mathfrak{a}} \oplus \mathfrak{n}_{-} \ | \
  {\varphi_v}_{|_{\mathfrak{b}}}  \textrm{ est }  \mathfrak{b}-\textrm{r{\'e}guli{\`e}re} \}$$
est un ouvert non vide de $\widehat{\mathfrak{a}} \oplus
  \mathfrak{n}_{-}$. L'intersection $\widetilde{V} \cap V$ est alors
  non vide. Soit $v$ dans cette intersection. Puisque $v$ appartient {\`a}
  $\widetilde{V}$, il existe un {\'e}l{\'e}ment $\rho$ dans $\widetilde{B}$
  tel que,
$$(\varphi_{{\rm pr}_{
    \mathfrak{n}_-}(v)})_{|_{\mathfrak{n}}}=\rho({\varphi_u}_{|_{\mathfrak{n}}}).$$
En particulier, pour tout $x$ dans $\mathfrak{n}$, on a\,:
$$\langle {\rm pr}_{
    \mathfrak{n}_-}(v), \rho(x) \rangle =\langle u,x \rangle .$$
Comme $v$ appartient {\`a} $V$, on dipose par ailleurs des relations suivantes,\:
$$ \dim \mathfrak{b}_{\varphi_v}= {\rm ind} \; \mathfrak{b} \leq \dim \mathfrak{b}_{\varphi_u}.$$
Il reste donc {\`a} prouver la relation, $\dim \mathfrak{b}_{\varphi_u}
\leq \dim \mathfrak{b}_{\varphi_v}$. La sous-alg{\`e}bre $ \mathfrak{b}$ est un id{\'e}al de
$\widetilde{\mathfrak{b}}$. Par suite, la sous-alg{\`e}bre $\rho(
\mathfrak{b}_{\varphi_u} )$ est contenue dans $\mathfrak{b}$. Prouvons
alors l'inclusion,
$$\rho( \mathfrak{b}_{\varphi_u} ) \subset
\mathfrak{b}_{\varphi_v}.$$
On en d{\'e}duira le r{\'e}sultat, car $\dim \rho( \mathfrak{b}_{\varphi_u} )
= \dim \mathfrak{b}_{\varphi_u}$. Cela revient {\`a} prouver la relation\,:
$$\langle [v,\rho( \mathfrak{b}_{\varphi_u} )] ,\mathfrak{b} \rangle =\{0\} .$$
Puisque $[\widehat{\mathfrak{a}},\mathfrak{b}] \subset \mathfrak{n}
\subset {\mathfrak{b}}^{\perp}$, cela revient {\`a} prouver la
relation\,:
$$\langle [{\rm pr}_{
    \mathfrak{n}_-}(v) ,\rho( \mathfrak{b}_{\varphi_u} )]
    ,\mathfrak{b} \rangle =\{0\} .$$
On a\,:
\begin{eqnarray*} 
\langle [{\rm pr}_{
    \mathfrak{n}_-}(v) ,\rho( \mathfrak{b}_{\varphi_u} )]
    ,\mathfrak{b} \rangle 
 & = & \langle {\rm pr}_{
    \mathfrak{n}_-}(v) , [\rho( \mathfrak{b}_{\varphi_u} )
    ,\mathfrak{b}] \rangle \\
& =&  \langle {\rm pr}_{
    \mathfrak{n}_-}(v) , \rho ([\mathfrak{b}_{\varphi_u} 
    ,\mathfrak{b}] )\rangle \\
& =&  \langle u , [\mathfrak{b}_{\varphi_u} 
    ,\mathfrak{b}] \rangle,  \textrm{ car }  [\mathfrak{b}_{\varphi_u} 
    ,\mathfrak{b}] \subset \mathfrak{n}\\
& =& \{0\}.
\end{eqnarray*}
Ceci termine la d{\'e}monstration, d'apr{\`e}s ce qui pr{\'e}c{\`e}de. 
\end{proof}

\subsection{}  La proposition suivante pr{\'e}cise la remarque \ref{rq_stab}\,:\\

\begin{proposition}\label{caract} La sous-alg{\`e}bre $\mathfrak{b}$ de
  $\mathfrak{g}$ poss{\`e}de une forme lin{\'e}aire stable si, et seulement si,
  l'ensemble $\mathcal{K}_{{\rm comp}}(\widehat{\Pi})$ est vide.
\end{proposition}

\begin{proof} Si l'ensemble $\mathcal{K}_{{\rm comp}}(\widehat{\Pi})$ est vide, on a d{\'e}j{\`a}
not{\'e} (remarque \ref{rq_stab}) que la restriction de $\varphi_u$ {\`a}
$\mathfrak{b}$ est
stable. R{\'e}ciproquement, supposons que
$\mathfrak{b}$ poss{\`e}de une forme lin{\'e}aire stable et montrons que l'ensemble
  $\mathcal{K}_{{\rm comp}}(\widehat{\Pi})$ est vide. On reprend les
  notations de la d{\'e}monstration
  du th{\'e}or{\`e}me \ref{formule_indice} et on pose\,:
$$V'=\{v \in \widehat{\mathfrak{a}} \oplus \mathfrak{n}_{-} \ | \ {\varphi_v}_{|_{\mathfrak{b}}} \textrm{ est }
  \mathfrak{b}-\textrm{stable} \} .$$
D'apr{\`e}s l'hypoth{\`e}se, l'ensemble $V'$ est un ouvert non vide de
$\widehat{\mathfrak{a}} \oplus \mathfrak{n}_{-}$ et l'intersection $V' \cap
  \widetilde{V}$ est alors non vide. Soit $v$ dans cette
  intersection. Alors il existe un {\'e}l{\'e}ment $\rho$ dans $\widetilde{B}$
  tel que,
$$(\varphi_{{\rm pr}_{
    \mathfrak{n}_-}(v)})_{|_{\mathfrak{n}}}
=\rho({\varphi_u}_{|_{\mathfrak{n}}}),$$
et il r{\'e}sulte de la d{\'e}monstration  du th{\'e}or{\`e}me \ref{formule_indice} la
    relation,
$$\mathfrak{b}_{\varphi_v}=\rho( \mathfrak{b}_{\varphi_u} ).$$
Par suite, $ {\varphi_v}_{|_{\mathfrak{b}}}$ est $\mathfrak{b}$-stable
si, et seulement si, $ {\varphi_u}_{|_{\mathfrak{b}}}$ l'est, d'apr{\`e}s
la proposition \ref{stable}. On
d{\'e}duit alors de la remarque \ref{rq_stab}, que l'ensemble
$\mathcal{K}_{{\rm comp}}(\widehat{\Pi})$ est n{\'e}cessairement vide.
\end{proof}

\section{Alg{\`e}bres de Lie quasi-r{\'e}ductives}\label{part_quasi}

\subsection{} Soit $\mathfrak{q}$ une
alg{\`e}bre de Lie alg{\'e}brique de centre $\mathfrak{z}$. La notion
d'alg{\`e}bre de Lie quasi-r{\'e}ductive a {\'e}t{\'e} introduite par
M. Duflo pour son importance dans l'analyse sur les groupes de Lie. On
en rappelle ici la d{\'e}finition.

\begin{definition}\label{red} On dit qu'une forme lin{\'e}aire $f$ de
  $\mathfrak{q}^*$ est {\rm r{\'e}ductive} si l'image de la sous-alg{\`e}bre
  $\mathfrak{q}_f / \mathfrak{z}$ dans $\mathfrak{gl}(\mathfrak{q})$
  par la repr{\'e}sentation adjointe de $\mathfrak{q}$ est une
  sous-alg{\`e}bre de Lie r{\'e}ductive.

On dit que l'alg{\`e}bre de Lie $\mathfrak{q}$ est {\rm quasi-r{\'e}ductive} si elle
poss{\`e}de une forme lin{\'e}aire r{\'e}ductive.
\end{definition}

Il est clair que si $\mathfrak{q}$ est r{\'e}ductive, alors $\mathfrak{q}$
est quasi-r{\'e}ductive; en effet la forme lin{\'e}aire nulle est r{\'e}ductive pour $\mathfrak{q}$. 

\begin{lemme}\label{tore} Une forme lin{\'e}aire est r{\'e}guli{\`e}re et
  r{\'e}ductive pour $\mathfrak{q}$ si, et seulement si, la sous-alg{\`e}bre
  $\mathfrak{q}_f / \mathfrak{z}$ est un tore de $\mathfrak{q}$.
\end{lemme}

\begin{proof} Si $\mathfrak{q}_f / \mathfrak{z}$ est un tore de
  $\mathfrak{q}$, alors $f$ est $\mathfrak{q}$-stable donc $f$ est
  $\mathfrak{q}$-r{\'e}guli{\`e}re. De plus il est clair que $f$ est r{\'e}ductive
  pour $\mathfrak{q}$.

R{\'e}ciproquement, si $f$ est $\mathfrak{q}$-r{\'e}guli{\`e}re et $\mathfrak{q}$-r{\'e}ductive, alors $\mathfrak{q}_f$ est une sous-alg{\`e}bre
de Lie commutative, d'apr{\`e}s  \cite{Tauvel2}. Il est alors bien connu (voir par
exemple \cite{Tauvel}, Th{\'e}or{\`e}me 20.5.10) que l'image de $\mathfrak{q}_f /
\mathfrak{z}$ dans $\mathfrak{gl}(\mathfrak{q})$
  par la repr{\'e}sentation adjointe de $\mathfrak{q}$  est une
  sous-alg{\`e}bre de Lie r{\'e}ductive si, et seulement si, elle est form{\'e}e
  d'{\'e}l{\'e}ments semi-simples.   
\end{proof}

Le th{\'e}or{\`e}me suivant est connu. On l'{\'e}nonce sans d{\'e}monstration\,:

\begin{theoreme}\label{quasi} On suppose que $\mathfrak{q}$
est quasi-r{\'e}ductive. Alors l'ensemble des formes lin{\'e}aires r{\'e}guli{\`e}res et
  r{\'e}ductives pour $\mathfrak{q}$ est un ouvert dense de
  $\mathfrak{q}^*$ form{\'e} de formes lin{\'e}aires $\mathfrak{q}$-stables. 
\end{theoreme}

En particulier, une alg{\`e}bre de Lie
quasi-r{\'e}ductive poss{\`e}de une forme lin{\'e}aire stable. 

\subsection{} On reprend les notations des parties pr{\'e}c{\'e}dentes. La
proposition suivante pr{\'e}cise la proposition \ref{caract}.

\begin{proposition}\label{caract2} La sous-alg{\`e}bre $\mathfrak{b}$ de
  $\mathfrak{g}$ est quasi-r{\'e}ductive si, et seulement si,
  l'ensemble $\mathcal{K}_{{\rm comp}}(\widehat{\Pi})$ est vide.\\
\end{proposition}

\begin{proof} Si $\mathcal{K}_{{\rm comp}}(\widehat{\Pi})$ n'est pas
  vide, $\mathfrak{b}$ ne poss{\`e}de pas de forme stable d'apr{\`e}s la
  proposition \ref{caract}. Il r{\'e}sulte donc du th{\'e}or{\`e}me \ref{quasi}
  que $\mathfrak{b}$ n'est pas quasi-r{\'e}ductive. R{\'e}ciproquement, si
  $\mathcal{K}_{{\rm comp}}(\widehat{\Pi})$ est vide, la restriction {\`a}
  $\mathfrak{b}$ de la forme
  $\varphi_u$ est $\mathfrak{b}$-stable et son stabilisateur est une
  sous-alg{\`e}bre commutative form{\'e}e d'{\'e}l{\'e}ments semi-simples, d'apr{\`e}s
  la proposition \ref{stab_u}. La forme lin{\'e}aire $\varphi_u$ est donc
  r{\'e}ductive pour $\mathfrak{b}$ d'apr{\`e}s le lemme \ref{tore}, d'o{\`u} la proposition. 
\end{proof}

La sous-alg{\`e}bre $\mathfrak{m} \oplus \widehat{\mathfrak{a}} \oplus
\mathfrak{n}$ est une sous-alg{\`e}bre parabolique de
  $\mathfrak{g}$. Les sous-alg{\`e}bres paraboliques de
  $\mathfrak{g}$ obtenues ainsi {\`a} partir d'une d{\'e}composition d'Iwasawa sont
  dites {\it minimales}. La proposition suivante est connue. On en
  redonne ici une d{\'e}monstration.

\begin{proposition}\label{m_quasi} La sous-alg{\`e}bre parabolique
  $\mathfrak{m} \oplus \widehat{\mathfrak{a}} \oplus \mathfrak{n}$ de
  $\mathfrak{g}$ est quasi-r{\'e}ductive. 
\end{proposition}

\begin{proof} Montrons que la sous-alg{\`e}bre de Lie $\mathfrak{q}_0=\mathfrak{m}_0 \oplus \widehat{\mathfrak{a}}_0 \oplus
  \mathfrak{n}_0$ de
  $\mathfrak{g}_0$ 
   est quasi-r{\'e}ductive. Soit $\varphi$ une forme lin{\'e}aire
   $\mathfrak{q}_0$-r{\'e}guli{\`e}re. Il suffit de prouver que son
   stabilisateur n'est form{\'e} que d'{\'e}l{\'e}ments semi-simples de $\mathfrak{g}_0$. Supposons par
  l'absurde le contraire. Comme le stabilisateur
  de $\varphi$ dans $\mathfrak{q}_0$  est une alg{\`e}bre de
  Lie alg{\'e}brique, il contient les composantes nilpotentes et
  semi-simples de ses {\'e}l{\'e}ments. Par suite, il contient un {\'e}l{\'e}ment
  nilpotent non nul $x$. La composante de
  $x$ selon $\widehat{\mathfrak{a}}_0\oplus\mathfrak{m}_0$ dans la d{\'e}composition $\mathfrak{q}_0=\mathfrak{m}_0 \oplus \widehat{\mathfrak{a}}_0 \oplus
  \mathfrak{n}_0$ est donc un {\'e}l{\'e}ment 
  nilpotent de $\mathfrak{g}_0$. Comme $[\mathfrak{m}_0,\mathfrak{m}_0]$ est une
sous-alg{\`e}bre de Lie compacte de $\mathfrak{g}_0$ (\cite{Knapp},
Chapitre VI), la sous-alg{\`e}bre $[\mathfrak{m}_0,\mathfrak{m}_0]$ est form{\'e}e d'{\'e}l{\'e}ments semi-simples de
$\mathfrak{g}_0$. Il en r{\'e}sulte que $x$ est contenu dans
$\mathfrak{n}_0$. 

Notons encore $\varphi$ la forme lin{\'e}aire sur $\mathfrak{q}$ obtenue {\`a}
partir de $\varphi$ par $\mathbb{C}$-lin{\'e}arit{\'e}. Soit alors $v$
l'{\'e}l{\'e}ment de $\mathfrak{m} \oplus \widehat{\mathfrak{a}} \oplus
  \mathfrak{n}_-$ qui s'identifie {\`a} $\varphi$ par la forme de
  Killing de $\mathfrak{g}$. Il r{\'e}sulte de \cite{Tauvel}, Proposition
  40.6.3, qu'il y a dans le dual de
$\mathfrak{n}$ une $\mathfrak{q}$-orbite ouverte. Par suite, on
peut supposer que $v$ s'{\'e}crit sous la forme $v=u+v'$, avec $v'$ dans
$\mathfrak{m} \oplus \widehat{\mathfrak{a}}$. Rappelons que $u= \sum\limits_{K \in
    \mathcal{K}^{''}(\widehat{\Pi})} X_{-\varepsilon_K}$ est introduit
  dans la partie \ref{part_stable}. Puisque $[x,v']$ appartient {\`a}
  $\mathfrak{n}$, on a\,: $[x,u] \in \mathfrak{n}$. Il r{\'e}sulte alors de la d{\'e}monstration de la proposition
\ref{stab_u} que $x$ est nul. On a ainsi obtenu une contradiction.
\end{proof}

\subsection{} D'apr{\`e}s la proposition \ref{m_quasi}, les sous-alg{\`e}bres paraboliques
minimales d'une alg{\`e}bre de simple complexe sont quasi-r{\'e}ductives. Il n'est en g{\'e}n{\'e}ral pas facile de
savoir si une sous-alg{\`e}bre parabolique est quasi-r{\'e}ductive. Les
r{\'e}sultats obtenus, entre autres, dans \cite{Panyushev3} et
\cite{Dvorsky} apportent de nombreuses r{\'e}ponses pour les cas
classiques. Mais le probl{\`e}me demeure pour une large part dans les cas
exceptionnels \cite{Moreau5}. On donne
dans ce paragraphe une caract{\'e}risation pour certaines sous-alg{\'e}bres
paraboliques.

On suppose que $\mathfrak{g}$ est une alg{\`e}bre de Lie simple complexe. On fixe une sous-alg{\`e}bre de Cartan $\mathfrak{h}$ de
$\mathfrak{g}$. Soit alors $\Delta$ le syst{\`e}me de racines associ{\'e} {\`a}
$(\mathfrak{g},\mathfrak{h})$, $\Delta_{+}$ un syst{\`e}me de racines
positives de $\Delta$, et $\Pi$ la base de $\Delta_{+}$. Soit $\mathfrak{b}$ la
sous-alg{\`e}bre de Borel standard (relativement {\`a} $\Pi$) de
$\mathfrak{g}$. Autrement dit, avec les notations des parties
pr{\'e}c{\'e}dentes, on se place dans le cas o{\`u} $\mathfrak{g}_0$ est une forme
r{\'e}elle split de
$\mathfrak{g}$. On reprend alors les notations pr{\'e}c{\'e}dentes. 

Notons  $E_{\Pi}$ le sous-espace de $\mathfrak{h}^*$ engendr{\'e} par les
{\'e}l{\'e}ments $\varepsilon_{K}$,  pour $K$  dans
$\mathcal{K}(\Pi)$. Introduisons trois sous-ensembles dans $\Delta_+$\,: 
\begin{eqnarray*}
\Delta_{+}^1 & = & \Delta_{+} \setminus  E_{\Pi}\\
\Delta_{+}^2 & = & \Delta_{+}  \cap \{ \ \varepsilon_{K} \ ,  \ K \in \mathcal{K}(\Pi) \ \}\\
\Delta_{+}^3 & = & \Delta_{+}  \cap  \left( E_{\Pi} \setminus\{ \
  \varepsilon_{K} \ ,  \ K \in \mathcal{K}(\Pi) \ \} \right)
\end{eqnarray*}
On a clairement\,:
$$\Delta_{+} = \Delta_{+}^1 \cup\Delta_{+}^2\cup
\Delta_{+}^3 .$$

Pour $\alpha \in \Delta_+$, on d{\'e}signe par $K_{\alpha}$ l'unique
{\'e}l{\'e}ment $K$ de $\mathcal{K}(\Pi)$ d{\'e}fini par le lemme \ref{kostant2}
(ii) tel que  $\alpha \in
\Gamma^{K_{\alpha}}$, et par $\mathcal{K}'(\alpha)$ l'ensemble des
{\'e}l{\'e}ments  $M \in \mathcal{K}(\Pi)$ tels que $\varepsilon_M+\alpha$ est
une racine.  D'apr{\`e}s le lemme \ref{kostant2},
(iv), on a\,: $\forall M \in \mathcal{K}'(\alpha), \ M
\subset_{\atop\hspace{-0.3cm}\not=} K_{\alpha}$. Notons $\Delta_{+}^{3'}$ l'ensemble des {\'e}l{\'e}ments  $\alpha$ dans $\Delta_{+}^{3}$
de la forme 
$$\alpha =  \displaystyle{\frac{1}{2}}(\varepsilon_{K_{\alpha}}-\varepsilon_{K_{\alpha}^{'}}),$$
o{\`u} $K_{\alpha}^{'}$ appartient {\`a} $\mathcal{K}'(\alpha)$. 

\begin{lemme}\label{cond_root} Soit $\alpha$ dans $\Pi$, $L \in
  \mathcal{K}(\Pi)$, $M \in
  \mathcal{K}'(\alpha)$ et $\gamma \in \Delta_+$, $\gamma \not=\alpha$.
 
\begin{itemize} \item[{\rm (i)}] $\varepsilon_L-\alpha$ est une racine
  si, et seulement si, $L=K_{\alpha}$,

\item[{\rm (ii)}] si $\alpha$ appartient {\`a} $\Delta_{+}^{3'}$, alors
  $\mathcal{K}'(\alpha)$ est r{\'e}duit {\`a} $\{K_{\alpha}^{'}\}$,

\item[{\rm (iii)}] si $\alpha$ appartient {\`a} $\Delta_{+}^{3} \setminus
\Delta_{+}^{3'}$, alors, pour tout $N \in
  \mathcal{K}(\Pi)$, $\varepsilon_M+\alpha \not= \varepsilon_N+\gamma$.
\end{itemize}
\end{lemme}

\begin{proof} (i) Supposons $L \not= K_{\alpha}$. Si
  $\varepsilon_L-\alpha$ est une racine, il r{\'e}sulte de la d{\'e}finition
  de $K_{\alpha}$ que $\alpha-\varepsilon_L$ est une racine
  positive, ce qui est impossible puisque $\alpha$ appartient {\`a} la base $\Pi$.\\
\\
(ii) Les tables \ref{classic} et
\ref{except} permettent de v{\'e}rifier que l'ensemble $\Pi \cap
\Delta_{+}^{3'}$ est non vide si, et seulement si, $\mathfrak{g}$ est
de type $F_4$, $C_l$ ou $B_l$, avec $l$ impair. On
v{\'e}rifie alors pour chacun de ces cas l'assertion (ii). On
reprend les notations des tables \ref{classic} et
\ref{except}.\\
\\ 
$F_4$\,: $\Pi \cap \Delta_{+}^{3'}=\{\beta_3=\displaystyle{\frac{1}{2}} (\varepsilon_{K_3}
-\varepsilon_{K_4})\}$, avec $K_{3}=\{\beta_2,\beta_3\}$ et
$K_4=\{\beta_2\}$. L'assertion (ii) est alors claire.\\
\\
$C_l$\,: $\Pi \cap \Delta_{+}^{3'}=\{\beta_i=\displaystyle{\frac{1}{2}} (\varepsilon_{K_i} -
\varepsilon_{K_{i+1}}), \ 1 \leq i \leq l-1\}$, 
o{\`u} $K_{j}=\{\beta_j, \ldots,\beta_l\}$, pour $j=1,\ldots,l$. Comme
$\varepsilon_{K_j} + \beta_{i}$ n'est pas une racine pour $j >i+1$,
on a (ii) pour $C_l$.\\
\\
$B_l$, $l$ impair\,: $\Pi \cap
\Delta_{+}^{3'}=\{\beta_{l}=\displaystyle{\frac{1}{2}}(\varepsilon_{K_{l}}-\varepsilon_{K_{l-1}})\}$,
avec $K_{l-1}=\{\beta_{l-1},\beta_{l}\}$ et $K_{l}=\{\beta_{l}\}$. L'assrtion (ii) est alors claire.\\
\\  
(iii) Supposons que la situation 
\begin{eqnarray}\label{situation}
\varepsilon_M+\alpha =
\varepsilon_N+\gamma.
\end{eqnarray} 
ait lieu, pour un certain $N$ dans $\mathcal{K}(\Pi)$. D'apr{\`e}s le lemme \ref{kostant2} (iv), les racines
$\varepsilon_M+\alpha$ et $\varepsilon_N+\gamma$ appartiennent {\`a}
$\Gamma_{0}^{K_{\alpha}}$. Comme $M \not=N$, il r{\'e}sulte alors de (i) que $\langle \alpha,
\varepsilon_{N}^{\vee} \rangle$ est n{\'e}gatif ou nul, d'o{\`u} $\langle \gamma,\varepsilon_{N}^{\vee}\rangle  \in\{-2,-3\}$. Si
$\langle \gamma,\varepsilon_{N}^{\vee}\rangle  =-3$, alors
$\mathfrak{g}$ est n{\'e}cssairement de type $G_2$. Or la situation
(\ref{situation}) n'a
pas lieu dans $G_2$. Si $\langle \gamma,\varepsilon_{N}^{\vee} \rangle = -2$, alors 
$\mathfrak{g}$ est n{\'e}cessairement de type $B_l$, $C_l$ ou $F_4$. On
v{\'e}rifie alors, gr{\^a}ce aux tables \ref{classic} et \ref{except}, que la
situation (\ref{situation}) n'a pas lieu dans chacun de ces types. 
\end{proof} 

Le lemme suivant, pr{\'e}sent{\'e} dans \cite{TauvelYu2} (Lemme 4.5) et
d{\'e}montr{\'e} dans \cite{Dixmier} (Lemme 1.12.2), intervient dans la
d{\'e}monstration du th{\'e}or{\`e}me \ref{car_sea}.

\begin{lemme}\label{hyperplan} Soit $V$ un espace vectoriel de
  dimension finie, $V'$ un hyperplan de $V$, $\Phi$ une forme bilin{\'e}aire
  sur $V$, et $\Phi'$ sa restriction {\`a} $V'$. On note $N$ et $N'$ les
  noyaux de $\Phi$ et $\Phi'$. \begin{itemize} \item[{\rm (i)}] Si $N
  \subset N'$, alors $N$ est un hyperplan de $N'$,

\item[{\rm (ii)}] Si $N \not\subset N'$, on a $N'=N \cap V'$, et $N'$
  est un hyperplan de $N$.\\

\end{itemize} 

\end{lemme}

Pour $S \subset \Pi$, on d{\'e}signe par $\mathfrak{p}_S$ la sous-alg{\`e}bre
parabolique {\it standard} (relativement {\`a} $\Pi$) $
\mathfrak{p}_S=\mathfrak{l}_S \oplus \mathfrak{u}_S$, o{\`u}
$\mathfrak{l}_S=\bigoplus\limits_{\alpha \in \Delta^S}
\mathfrak{g}^{\alpha}$ est la partie Levi de $\mathfrak{p}_S$, et o{\`u}
$\mathfrak{u}_S=\bigoplus\limits_{\alpha \in \Delta_+ \setminus \Delta^S}
\mathfrak{g}^{\alpha}$ est son radical nilpotent. Comme les
sous-alg{\`e}bres paraboliques de $\mathfrak{g}$ sont conjugu{\'e}es aux
sous-alg{\`e}bres paraboliques standards, il suffit de consid{\'e}rer les
sous-alg{\`e}bres paraboliques standards. Dans le th{\'e}or{\`e}me suivant une
donne une description des sous-alg{\`e}bres paraboliques quasi-r{\'e}ductives
standards pour $S=\{\alpha\}$, avec $\alpha \in \Pi$.

\begin{theoreme}\label{car_sea} Soit $\alpha$ une racine simple. Alors
  la sous-alg{\`e}bre parabolique  $\mathfrak{p}_{\{\alpha\}}$ est
  quasi-r{\'e}ductive si, et seulement si, $\alpha$ appartient {\`a}  $\Delta_{+}^{1} \cup \Delta_{+}^{2}
 \cup \Delta_{+}^{3'}$.
\end{theoreme}

\begin{proof}  Pour simplifer les notations, notons $\mathfrak{q}$ la
  sous-alg{\`e}bre $\mathfrak{p}_{\{\alpha\}}$,
  $\mathfrak{l}_{\mathfrak{q}}$ sa partie Levi, et 
  $\mathfrak{u}_{\mathfrak{q}}={\mathfrak{q}}^{\perp}$ son radical
  nilpotent.\\
 
Si $\alpha$ appartient $\Delta_{+}^{1} \cup
  \Delta_{+}^{2}$, $\mathfrak{p}_{\alpha}$ est  quasi-r{\'e}ductive
  d'apr{\`e}s la d{\'e}monstration de la proposition 4.4 de \cite{Tauvel}. Supposons que  $\alpha$
  appartient {\`a} $\Delta_{+}^{3'}$. La sous-alg{\`e}bre de Borel  $\mathfrak{b}=\mathfrak{h} \oplus \sum\limits_{\beta \in \Delta_+}
\mathfrak{g}^{\beta}$ est un  hyperplan de  $\mathfrak{q}$. Posons
$u=\sum\limits_{K \in \mathcal{K}(\Pi)}
  X_{-\varepsilon_K}$. 
Le stabilisateur de 
$\varphi_{u}$ dans $\mathfrak{b}$ est donn{\'e} par la relation suivante\,:
$$\mathfrak{b}_{\varphi_{u}}=\bigcap\limits_{K \in \mathcal{K}(\Pi)} \ker
\varepsilon_K,$$
d'apr{\`e}s \cite{TauvelYu} ou d'apr{\`e}s la proposition \ref{stab_u}. Par cons{\'e}quent, $\mathfrak{b}_{\varphi_u}$
est contenu dans le stabilisateur  $\mathfrak{q}_{\varphi_{u}}$ de
$\varphi_{u}$, car 
$[u,\mathfrak{b}_{\varphi_{u}}]=\{0\}$. Il r{\'e}sulte alors du lemme \ref{hyperplan} que
$\mathfrak{b}_{\varphi_{u}}$ est un  hyperplan de
$\mathfrak{q}_{\varphi_{u}}$. Montrons que   $X_{-\alpha}+b
X_{\alpha}$ appartient {\`a} $\mathfrak{q}_{\varphi_{u}}$, o{\`u}
$b$ est {\`a} d{\'e}finir. Puisque $\alpha$ appartient {\`a} $\Delta_{+}^{3'}$, on
a $\alpha=\displaystyle{\frac{1}{2}}(
\varepsilon_{K_{\alpha}}-\varepsilon_{K_{\alpha}^{'}})$, pour 
$K_{\alpha}^{'}$ dans $\mathcal{K}^{'}(\alpha)$. De plus,
d'apr{\`e}s le lemme \ref{cond_root} (i),
$\varepsilon_M - \alpha$ est une racine si, et seulement si, 
$M=K_{\alpha}$. Par suite, on a\,:
\begin{eqnarray*}
[u,X_{-\alpha}+b
X_{\alpha}] & = & \sum\limits_{M \in \mathcal{K}(\Pi)}
[X_{-\varepsilon_M},X_{-\alpha}] + b \sum\limits_{M \in \mathcal{K}(\Pi)}
[X_{-\varepsilon_M},X_{\alpha}] \\
&= & \lambda X_{-(\varepsilon_L+\alpha)}+ b \mu  X_{-(\varepsilon_K-\alpha)},
\end{eqnarray*}
o{\`u} $\lambda$ et $\mu$ sont des complexes non nuls. Puisque
$\varepsilon_L+\alpha=\varepsilon_K-\alpha$, on peut choisir $b$ de
sorte que $[u,X_{-\alpha}+b
X_{\alpha}]=0$. Par cons{\'e}quent la relation  
$$ \mathfrak{q}_{\varphi_{u}}=\mathbb{C} (X_{-\alpha}+b
X_{\alpha}) \oplus \mathfrak{b}_{\varphi_{u}},$$
a lieu, pour $\mathfrak{b}_{\varphi_{u}}$ un hyperplan
$\mathfrak{q}_{\varphi_{u}}$. L'alg{\`e}bre de Lie  
$\mathfrak{q}_{\varphi_{u}}$ est  ab{\'e}lienne,
car $\alpha$ appartient {\`a} $E_{\Pi}$. De plus,
$\mathfrak{q}_{\varphi_{u}}$ est form{\'e} d'{\'e}l{\'e}ments semi-simples. Ainsi
$\varphi_{u}$ est 
$\mathfrak{q}$-r{\'e}ductive, d'apr{\`e}s le lemme \ref{tore} et $\mathfrak{q}$ est quasi-r{\'e}ductive.\\

R{\'e}ciproquement, supposons que  $\alpha$ appartienne {\`a}  $\Delta_{+}^{3} \setminus
\Delta_{+}^{3'} $. Supposons par l'absurde que $\mathfrak{q}$ est 
quasi-r{\'e}ductive. On cherche {\`a} aboutir {\`a} une contradiction. Le radical nilpotent $\mathfrak{u}_{\mathfrak{q}}$
de $\mathfrak{q}$ est un id{\'e}al de $\mathfrak{b}$ contenu dans le
radical nilpotent $\mathfrak{u}$ de
$\mathfrak{b}$. Soit $B$ le groupe adjoint alg{\'e}brique de
$\mathfrak{b}$. Il r{\'e}sulte de la proposition
40.6.3 de \cite{Tauvel} que la $B$-orbite de la forme lin{\'e}aire  ${\varphi_{u}}_{|_{\mathfrak{u}_{\mathfrak{q}}}}$
  est un ouvert dense de $\mathfrak{u}_{\mathfrak{q}}^*$. Le dual de
  $\mathfrak{u}_{\mathfrak{q}}$ s'identifie via la forme  de Killing
  de 
  $\mathfrak{g}$ au sous-espace  $\mathfrak{u}_{\mathfrak{q}_-}=\sum\limits_{\beta \in \Delta_{+}
    \setminus \{ \alpha \}} \mathfrak{g}^{-\beta}$. Par suite, l'ensemble
$$\widetilde{W}=\{w \in 
\mathfrak{b}_{-} \oplus \mathfrak{g}^{\alpha}
\ | \ (\varphi_{ {\rm pr}_{\mathfrak{u}_{\mathfrak{q}_-} } (w) })_{|_{
    \mathfrak{u}_{ \mathfrak{q} } } } \in  B \cdot {\varphi_{ u }
}_{|_{\mathfrak{u}_{\mathfrak{q} } }} \}$$
est un ouvert dense de $\mathfrak{b}_- \oplus
\mathfrak{g}^{\alpha}$, o{\`u} ${\rm
  pr}_{\mathfrak{u}_{\mathfrak{q}_-}}$ est la  projection de 
   $\mathfrak{q}_-$ sur
  $\mathfrak{u}_{\mathfrak{q}_-}$ selon la  d{\'e}composition 
    $\mathfrak{q}_-=\mathfrak{l}_{\mathfrak{q}}\oplus \mathfrak{u}_{\mathfrak{q}_-}$.\\

Par ailleurs, puisque le dual de $\mathfrak{q}$ s'identifie via la
forme de Killing de $\mathfrak{g}$ {\`a} 
$\mathfrak{q}_{-}=\mathfrak{b}_{-} \oplus \mathfrak{g}_{\alpha}$, l'ensemble  
$$W=\{w \in \mathfrak{b}_{-} \oplus \mathfrak{g}_{\alpha}
\ | \ \varphi_w \textrm{ est } \mathfrak{q}-\textrm{r{\'e}guli{\`e}re et }
\mathfrak{q}-\textrm{r{\'e}ductive}\}$$
est un ouvert dense de $\mathfrak{b}_{-} \oplus
\mathfrak{g}_{\alpha}$, d'apr{\`e}s l'hypoth{\`e}se. L'intersection $W
\cap \widetilde{W}$ est alors non vide. Soit $w$ dans cette
intersection. On montre alors sans difficult{\'e} que $w$ peut
s'{\'e}crire sous la forme\,:
$$w=u+a X_{-\alpha} + h + b X_{+\alpha},$$
avec $a,b$ dans $\mathbb{C}$ et $h$ dans $\mathfrak{h}$. Comme 
$\alpha$ appartient {\`a}  $E_{\Pi}$, on a 
$[\mathfrak{b}_{\varphi_{u}},w]=\{0\}$, d'o{\`u} l'inclusion
$\mathfrak{b}_{\varphi_{u}} \subset
\mathfrak{q}_{\varphi_w}$. Il r{\'e}sulte alors du lemme 
\ref{hyperplan} que 
$\mathfrak{b}_{\varphi_{u}}$ est un  hyperplan
de $\mathfrak{q}_{\varphi_w}$. Soit alors $x$ dans $\mathfrak{q}$ tel
que  l'on ait la d{\'e}composition\,:
$$\mathfrak{q}_{\varphi_w} =\mathbb{C} x \oplus (\bigcap\limits_{K \in \mathcal{K}(\Pi)} \ker
\varepsilon_K ).$$
Puisque $w$ appartient {\`a} $W$, $\mathfrak{q}_{\varphi_w}$ est une
alg{\`e}bre de Lie ab{\'e}lienne form{\'e}e d'{\'e}l{\'e}ments  semi-simples. En
particulier, la composante $x_{\mathfrak{l}_{\mathfrak{q}}}$ de $x$ sur
 $\mathfrak{l}_{\mathfrak{q}}$ dans la d{\'e}composition
$\mathfrak{q}=\mathfrak{l}_{\mathfrak{q}}\oplus
\mathfrak{u}_{\mathfrak{q}}$ est semi-simple. L'{\'e}l{\'e}ment  $x$ s'{\'e}crit,
$$x=\lambda X_{-\alpha} +H + \mu X_{\alpha}+x_{\mathfrak{u}_{\mathfrak{q}}},$$
avec $\lambda, \mu$ dans $\mathbb{C}$, $H$ dans $\mathfrak{h}$ et
$x_{\mathfrak{u}_{\mathfrak{q}}}$ dans $ \mathfrak{u}_{\mathfrak{q}}$. Si
$\lambda$ {\'e}tait nul, alors  $\mu$ serait nul aussi, car la composante
de $x$ selon $\mathfrak{l}_{\mathfrak{q}}$ est semi-simple. Par suite, la relation $[x,w] \in
{\mathfrak{q}}^{\perp}$ impliquerait que 
$x_{\mathfrak{u}_{\mathfrak{q}}}$ appartient {\`a}
$\mathfrak{q}_{\varphi_w}$. Il en r{\'e}sulterait que 
$x_{\mathfrak{u}_{\mathfrak{q}}}$ est un {\'e}l{\'e}ment non nul de $\mathfrak{u}_{\mathfrak{q}}$, car $\bigcap\limits_{K \in \mathcal{K}(\Pi)} \ker
\varepsilon_K$ et $\mathbb{C} x$ sont en somme directe. Ceci serait
alors en contradication avec le fait que  $\mathfrak{q}_{\varphi_w}$
n'est form{\'e} que d'{\'e}l{\'e}ments semi-simples. Par suite $\lambda=0$ et  on
obtient de la m{\^e}me mani{\`e}re que $\mu$ n'est pas nul.  D'apr{\`e}s le lemme \ref{cond_root}
(i), on a alors
\begin{eqnarray*}
[x,w] & =&-\lambda \sum\limits_{L \in \mathcal{K}'(\alpha)} [X_{-\varepsilon_L},X_{-\alpha}]
+\sum\limits_{L \in \mathcal{K}(\Pi)} \varepsilon_K(H) X_{-\varepsilon_K}
-\mu [X_{-\varepsilon_{K_{\alpha}}},X_{\alpha}] \\
& & + (a \alpha(H) + \lambda \alpha(h)) X_{-\alpha} +( a \mu -b
\lambda) H_{\alpha} + (-b \alpha(H) - \mu \alpha(h)) X_{\alpha}\\
&& + [x_{\mathfrak{u}_{\mathfrak{q}}},w].
\end{eqnarray*}

Puisque  $[x,w]$ appartient {\`a} $\mathfrak{u}_{\mathfrak{q}}$, l'{\'e}l{\'e}ment
 $[X_{-\varepsilon_{K_{\alpha}}},X_{\alpha}]$ de $\mathfrak{u}_{\mathfrak{q}_-}$ doit
 {\^e}tre compens{\'e}. Supposons qu'il se compense gr{\^a}ce au terme 
$[x_{\mathfrak{u}_{\mathfrak{q}}},w]$. Alors, il existe $M \in \mathcal{K}(\Pi)$ et $\beta \in
\Delta_+ \setminus \{\alpha\}$ tels que 
$$\varepsilon_{K_{\alpha}}-\alpha=\varepsilon_M - \beta.$$
Comme $M \not=K_{\alpha}$, on a $\langle \beta,
\varepsilon_{K_{\alpha}}^{\vee}\rangle=-1$. Par cons{\'e}quent, $M$ 
appartient {\`a} $\mathcal{K}'(\alpha)$, car 
$\varepsilon_{K_{\alpha}}+\beta=\varepsilon_M+\alpha$ est une
racine. Ceci est en contradiction avec le lemme \ref{cond_root} (ii). Ce qui
pr{\'e}c{\`e}de prouve alors qu'il existe un {\'e}l{\'e}ment  $L$ dans
$\mathcal{K}^{'}(\alpha)$ tel que 
$\varepsilon_{K_{\alpha}}-\alpha=\varepsilon_L+\alpha$. Par suite,
$\alpha=\displaystyle{\frac{1}{2}}(\varepsilon_{K_{\alpha}}-\varepsilon_L)$,
ce qui est absurde car $\alpha$ n'appartient pas {\`a} $\Delta_+^{3'}$. On a
ainsi obtenu une contradiction.
\end{proof} 

Dans la table
\ref{colored_root}, on a repr{\'e}sent{\'e} en noire, pour chaque type
d'alg{\`e}bre de Lie simple complexe $\mathfrak{g}$, les racines simples
$\alpha$ telles que $\mathfrak{p}_{\{\alpha\}}$ est une  sous-alg{\`e}bre
quasi-r{\'e}ductive. On ommet ici les
calculs qui ont {\'e}t{\'e} effectu{\'e}s pour {\'e}tablir cette table. Ces
calculs sont d{\'e}taill{\'e}s dans \cite{Moreau5}.

Par ailleurs, le th{\'e}or{\`e}me \ref{car_sea} peut permettre de d{\'e}cider dans
certains cas si la condition $({\tt *})$ est satifaite (voir le cas de $EVI$).

\hspace{-4cm}
{\footnotesize
\begin{table}[h]
\begin{tabular}{|c|c|}
\hline
$A_l$, $l \geq 1$ &
\begin{pspicture}(-2.5,-0)(3.5,1)
\pscircle[fillstyle=solid, fillcolor=black](-2,0){1mm}
\pscircle[fillstyle=solid, fillcolor=black](-1,0){1mm}
\pscircle[fillstyle=solid, fillcolor=black](0,0){1mm}
\pscircle[fillstyle=solid, fillcolor=black](1,0){1mm}
\pscircle[fillstyle=solid, fillcolor=black](2,0){1mm}
\pscircle[fillstyle=solid, fillcolor=black](3,0){1mm}
\psline(-1.9,0)(-1.1,0)
\psline(-0.1,0)(-0.9,0)
\psline[linestyle=dotted](0.1,0)(0.9,0)
\psline(1.1,0)(1.9,0)
\psline(2.1,0)(2.9,0)

\rput[b](-2,0.2){$\beta_1$}
\rput[b](-1,0.2){$\beta_2$}
\rput[b](2,0.2){$\beta_{l-1}$}
\rput[b](3,0.2){$\beta_l$}
\end{pspicture} \\
&\\
\hline
$B_l$, $l \geq 2$, $l$ impair &
\begin{pspicture}(-2.5,-0)(3.5,1)
\pscircle[fillstyle=solid, fillcolor=black](-2,0){1mm}
\pscircle(-1,0){1mm}
\pscircle[fillstyle=solid, fillcolor=black](0,0){1mm}
\pscircle[fillstyle=solid, fillcolor=black](1,0){1mm}
\pscircle(2,0){1mm}
\pscircle[fillstyle=solid, fillcolor=black](3,0){1mm}
\psline(-1.9,0)(-1.1,0)
\psline(-0.1,0)(-0.9,0)
\psline[linestyle=dotted](0.1,0)(0.9,0)
\psline(1.1,0)(1.9,0)
\psline(2.1,0.05)(2.9,0.05)
\psline(2.1,-0.05)(2.9,-0.05)
\rput[b](-2,0.2){$\beta_1$}
\rput[b](-1,0.2){$\beta_2$}
\rput[b](2,0.2){$\beta_{l-1}$}
\rput[b](3,0.2){$\beta_l$}
\rput(2.5,0){$>$}
\end{pspicture}\\
&\\
\hline
$B_l$, $l \geq 2$, $l$ pair &
\begin{pspicture}(-2.5,-0)(3.5,1)
\pscircle[fillstyle=solid, fillcolor=black](-2,0){1mm}
\pscircle(-1,0){1mm}
\pscircle[fillstyle=solid, fillcolor=black](0,0){1mm}
\pscircle(1,0){1mm}
\pscircle[fillstyle=solid, fillcolor=black](2,0){1mm}
\pscircle[fillstyle=solid, fillcolor=black](3,0){1mm}
\psline(-1.9,0)(-1.1,0)
\psline(-0.1,0)(-0.9,0)
\psline[linestyle=dotted](0.1,0)(0.9,0)
\psline(1.1,0)(1.9,0)
\psline(2.1,0.05)(2.9,0.05)
\psline(2.1,-0.05)(2.9,-0.05)
\rput[b](-2,0.2){$\beta_1$}
\rput[b](-1,0.2){$\beta_2$}
\rput[b](2,0.2){$\beta_{l-1}$}
\rput[b](3,0.2){$\beta_l$}
\rput(2.5,0){$>$}
\end{pspicture}\\
&\\
\hline
$C_l$, $l \geq 3$&
\begin{pspicture}(-2.5,0)(3.5,1)
\pscircle[fillstyle=solid, fillcolor=black](-2,0){1mm}
\pscircle[fillstyle=solid, fillcolor=black](-1,0){1mm}
\pscircle[fillstyle=solid, fillcolor=black](0,0){1mm}
\pscircle[fillstyle=solid, fillcolor=black](1,0){1mm}
\pscircle[fillstyle=solid, fillcolor=black](2,0){1mm}
\pscircle[fillstyle=solid, fillcolor=black](3,0){1mm}
\psline(-1.9,0)(-1.1,0)
\psline(-0.1,0)(-0.9,0)
\psline[linestyle=dotted](0.1,0)(0.9,0)
\psline(1.1,0)(1.9,0)
\psline(2.1,0.05)(2.9,0.05)
\psline(2.1,-0.05)(2.9,-0.05)
\rput[b](-2,0.2){$\beta_1$}
\rput[b](-1,0.2){$\beta_2$}
\rput[b](2,0.2){$\beta_{l-1}$}
\rput[b](3,0.2){$\beta_l$}
\rput(2.5,0){$<$}
\end{pspicture}\\
&\\
\hline
$D_l$, $l \geq 5$, $l$ impair &
\begin{pspicture}(-2.5,-0.6)(3.5,1.1)
\pscircle[fillstyle=solid, fillcolor=black](-2,0){1mm}
\pscircle(-1,0){1mm}
\pscircle[fillstyle=solid, fillcolor=black](0,0){1mm}
\pscircle(1,0){1mm}
\pscircle[fillstyle=solid, fillcolor=black](2,0){1mm}
\pscircle[fillstyle=solid, fillcolor=black](3,0.7){1mm}
\pscircle[fillstyle=solid, fillcolor=black](3,-0.7){1mm}
\psline(-1.9,0)(-1.1,0)
\psline(-0.1,-0)(-0.9,0)
\psline[linestyle=dotted](0.1,0)(0.9,0)
\psline(1.1,0)(1.9,0)
\psline(2.07,0.05)(2.92,0.68)
\psline(2.07,-0.05)(2.92,-0.68)
\rput[b](-2,0.2){$\beta_1$}
\rput[b](-1,0.2){$\beta_2$}
\rput[b](2,0.2){$\beta_{l-2}$}
\rput[l](3.1,0.8){$\beta_l$}
\rput[l](3.1,-0.8){$\beta_{l-1}$}
\end{pspicture}\\
&\\
\hline
$D_l$, $l \geq 4$, $l$ pair &
\begin{pspicture}(-2.5,-0.6)(3.5,1.1)
\pscircle[fillstyle=solid, fillcolor=black](-2,0){1mm}
\pscircle(-1,0){1mm}
\pscircle[fillstyle=solid, fillcolor=black](0,0){1mm}
\pscircle[fillstyle=solid, fillcolor=black](1,0){1mm}
\pscircle(2,0){1mm}
\pscircle[fillstyle=solid, fillcolor=black](3,0.7){1mm}
\pscircle[fillstyle=solid, fillcolor=black](3,-0.7){1mm}
\psline(-1.9,0)(-1.1,0)
\psline(-0.1,-0)(-0.9,0)
\psline[linestyle=dotted](0.1,0)(0.9,0)
\psline(1.1,0)(1.9,0)
\psline(2.07,0.05)(2.92,0.68)
\psline(2.07,-0.05)(2.92,-0.68)
\rput[b](-2,0.2){$\beta_1$}
\rput[b](-1,0.2){$\beta_2$}
\rput[b](2,0.2){$\beta_{l-2}$}
\rput[l](3.1,0.8){$\beta_l$}
\rput[l](3.1,-0.8){$\beta_{l-1}$}
\end{pspicture}\\
&\\
\hline
$E_6$ & \begin{pspicture}(-2.5,-0.8)(2.5,1)
\pscircle[fillstyle=solid, fillcolor=black](-2,0){1mm}
\pscircle[fillstyle=solid, fillcolor=black](-1,0){1mm}
\pscircle[fillstyle=solid, fillcolor=black](0,0){1mm}
\pscircle[fillstyle=solid, fillcolor=black](1,0){1mm}
\pscircle[fillstyle=solid, fillcolor=black](2,0){1mm}
\pscircle(0,-1){1mm}

\psline(-1.9,0)(-1.1,0)
\psline(-0.1,0)(-0.9,0)
\psline(0.1,0)(0.9,0)
\psline(1.1,0)(1.9,0)
\psline(0,-0.1)(0,-0.9)

\rput[b](-2,0.2){$\beta_1$}
\rput[b](-1,0.2){$\beta_3$}
\rput[b](0,0.2){$\beta_4$}
\rput[b](1,0.2){$\beta_5$}
\rput[b](2,0.2){$\beta_6$}
\rput[r](-0.2,-1){$\beta_2$}
\end{pspicture}\\
&\\
\hline
$E_7$ &
\begin{pspicture}(-2.5,-0.8)(3.5,1)
\pscircle(-2,0){1mm}
\pscircle[fillstyle=solid, fillcolor=black](-1,0){1mm}
\pscircle(0,0){1mm}
\pscircle[fillstyle=solid, fillcolor=black](1,0){1mm}
\pscircle(2,0){1mm}
\pscircle[fillstyle=solid, fillcolor=black](3,0){1mm}

\pscircle[fillstyle=solid, fillcolor=black](0,-1){1mm}

\psline(-1.9,0)(-1.1,0)
\psline(-0.1,0)(-0.9,0)
\psline(0.1,0)(0.9,0)
\psline(1.1,0)(1.9,0)
\psline(2.1,0)(2.9,0)
\psline(0,-0.1)(0,-0.9)

\rput[b](-2,0.2){$\beta_1$}
\rput[b](-1,0.2){$\beta_3$}
\rput[b](0,0.2){$\beta_4$}
\rput[b](1,0.2){$\beta_5$}
\rput[b](2,0.2){$\beta_6$}
\rput[b](3,0.2){$\beta_7$}
\rput[r](-0.2,-1){$\beta_2$}
\end{pspicture}\\
&\\
\hline
$E_8$ & 
\begin{pspicture}(-2.5,-0.8)(4.5,1)
\pscircle(-2,0){1mm}
\pscircle[fillstyle=solid, fillcolor=black](-1,0){1mm}
\pscircle(0,0){1mm}
\pscircle[fillstyle=solid, fillcolor=black](1,0){1mm}
\pscircle(2,0){1mm}
\pscircle[fillstyle=solid, fillcolor=black](3,0){1mm}
\pscircle(4,0){1mm}
\pscircle[fillstyle=solid, fillcolor=black](0,-1){1mm}

\psline(-1.9,0)(-1.1,0)
\psline(-0.1,0)(-0.9,0)
\psline(0.1,0)(0.9,0)
\psline(1.1,0)(1.9,0)
\psline(2.1,0)(2.9,0)
\psline(3.1,0)(3.9,0)
\psline(0,-0.1)(0,-0.9)

\rput[b](-2,0.2){$\beta_1$}
\rput[b](-1,0.2){$\beta_3$}
\rput[b](0,0.2){$\beta_4$}
\rput[b](1,0.2){$\beta_5$}
\rput[b](2,0.2){$\beta_6$}
\rput[b](3,0.2){$\beta_7$}
\rput[b](4,0.2){$\beta_8$}
\rput[r](-0.2,-1){$\beta_2$}
\end{pspicture}\\
&\\
\hline
$F_4$ &
\begin{pspicture}(-2.5,0)(2.5,1)

\pscircle(-2,0){1mm}
\pscircle[fillstyle=solid, fillcolor=black](-1,0){1mm}
\pscircle[fillstyle=solid, fillcolor=black](0,0){1mm}
\pscircle[fillstyle=solid, fillcolor=black](1,0){1mm}

\psline(-1.9,0)(-1.1,0)
\psline(-0.09,0.05)(-0.91,0.05)
\psline(-0.09,-0.05)(-0.91,-0.05)
\psline(0.1,0)(0.9,0)

\rput[b](-2,0.2){$\beta_{1}$}
\rput[b](-1,0.2){$\beta_{2}$}
\rput[b](0,0.2){$\beta_{3}$}
\rput[b](1,0.2){$\beta_{4}$}

\rput(-0.5,0){$>$}
\end{pspicture}\\
&\\
\hline
$G_2$ &
\begin{pspicture}(-2.5,0)(1.5,1)
\pscircle(-1,0){1mm}
\pscircle[fillstyle=solid, fillcolor=black](0,0){1mm}

\psline(-0.09,0.06)(-0.91,0.06)
\psline(-0.09,0)(-0.91,0)
\psline(-0.09,-0.06)(-0.91,-0.06)

\rput[b](-1,0.2){$\beta_{1}$}
\rput[b](0,0.2){$\beta_{2}$}

\rput(-0.5,0){$>$}
\end{pspicture}\\
&\\
\hline
\end{tabular}
~\\
\\
\caption{\label{colored_root} Sous-alg{\`e}bres paraboliques
  quasi-r{\'e}ductives associ{\'e}es {\`a} une racine simple (Th{\'e}or{\`e}me \ref{car_sea})}
\end{table}}

\section{Calculs explicites dans les alg{\`e}bres de Lie simples r{\'e}elles}\label{part_calcul}
On suppose dans cette partie  que $\mathfrak{g}_0$ est une alg{\`e}bre de Lie
simple r{\'e}elle. 

\subsection{} Consid{\'e}rons les propri{\'e}t{\'e}s suivantes\,:\\

\begin{description} \item[(A)] L'indice de $\mathfrak{b}$ est donn{\'e}
  par la relation\,: ${\rm ind} \; \mathfrak{b}={\rm rg} \; \mathfrak{g} - {\rm rg} \;
\mathfrak{k}$. Autrement dit l'indice est additif dans la
d{\'e}composition $\mathfrak{g}=\mathfrak{k} \oplus \mathfrak{b}$.\\

\item[(B)] La sous-alg{\`e}bre $\mathfrak{b}$ poss{\`e}de une forme stable.\\

\item[(B)'] La sous-alg{\`e}bre $\mathfrak{b}$ est quasi-r{\'e}ductive.\\

\item[(C)] Il y a dans le dual $\mathfrak{b}^{*}$ de  $\mathfrak{b}$
  une $\mathfrak{b}$-orbite ouverte. Autrement dit l'indice de $\mathfrak{b}$ est nul.\\

\end{description}

D'apr{\`e}s les propositions \ref{caract} et \ref{caract2}, les conditions
(b) et (b)' sont {\'e}quivalentes. En outre, d'apr{\`e}s les parties pr{\'e}c{\'e}dentes, on dispose de crit{\`e}res simples pour chacune
de ces conditions. Pr{\'e}cis{\'e}ment, avec les notations des parties pr{\'e}c{\'e}dentes, on a\,:

$\mathfrak{b}$ v{\'e}rifie {\bf (A)} si, et seulement, si $\mathcal{K}(\widehat{\Pi}^{'})=\mathcal{K}^{'}(\widehat{\Pi})$ (Th{\'e}or{\`e}me \ref{formule_indice}),

$\mathfrak{b}$ v{\'e}rifie {\bf (B)} ou {\bf (B)'} si, et seulement, si $\mathcal{K}_{{\rm
    comp}}(\widehat{\Pi})$ est vide (Propositions \ref{caract} et \ref{caract2}) et, 

$\mathfrak{b}$ v{\'e}rifie {\bf (C)} si, et seulement, si 
($\mathcal{K}(\widehat{\Pi}^{'})=\mathcal{K}^{'}(\widehat{\Pi})$ et ${\rm rg} \; \mathfrak{g} = {\rm rg} \;
\mathfrak{k}$) (Th{\'e}or{\`e}me \ref{formule_indice}).\\

Il est clair que la condition {\bf (C)} implique la condition {\bf (B)}. En effet, si
l'indice de $\mathfrak{b}$ est nul, alors il existe un {\'e}l{\'e}ment de $\mathfrak{b}^*$ dont le stabilisateur est nul. Si {\bf (B)} est v{\'e}rifi{\'e}e, alors l'ensemble $\mathcal{K}_{{\rm
    comp}}(\widehat{\Pi})$ est vide, d'apr{\`e}s ce qui pr{\'e}c{\`e}de. Il
r{\'e}sulte alors du lemme \ref{gamma1} (ii) que la condition
$\mathcal{K}(\widehat{\Pi}^{'})=\mathcal{K}^{'}(\widehat{\Pi})$ est
automatiquement remplie. La condition {\bf (A)} est alors satisfaite, d'apr{\`e}s
ce qui pr{\'e}c{\`e}de. Par suite, {\bf (B)} implique {\bf (A)}.

\subsection{} On cherche maintenant pour quels types d'alg{\`e}bres de Lie simples 
r{\'e}elles la sous-alg{\`e}bre $\mathfrak{b}$ v{\'e}rifie les conditions {\bf (A)},
{\bf (B)}, {\bf (B)'} ou {\bf (C)}.  D'apr{\`e}s la classification des alg{\`e}bres de Lie simples 
r{\'e}elles obtenue par exemple dans \cite{Knapp}, Th{\'e}or{\`e}me
6.105, l'alg{\`e}bre de Lie $\mathfrak{g}_0$ est isomorphe {\`a} l'une des alg{\`e}bres de
Lie simples r{\'e}elles de
la liste suivante\,:\\

\begin{description}
\item[(a)] L'alg{\`e}bre de Lie $\mathfrak{s}^{\mathbb{R}}$, o{\`u} $\mathfrak{s}$ est
  simple complexe de type $A_n$, pour $n \geq 1$, $B_n$, pour $n \geq
  2$, $C_n$, pour $n \geq 3$, $D_n$, pour $n \geq 4$, $E_6$, $E_7$,
  $E_8$, $F_4$ ou $G_2$,

\item[(b)] La forme r{\'e}elle compacte d'une alg{\`e}bre de Lie
  $\mathfrak{s}$ comme en {\bf (a)},

\item[(c)] Les alg{\`e}bres de matrices classiques\,: \begin{itemize}
\item $\mathfrak{sl}(n,\mathbb{R})$, avec $n \geq 2$,
\item $\mathfrak{sl}(n,\mathbb{H})$, avec $n \geq 2$,
\item $\mathfrak{su}(p,q)$, avec $p \geq q > 0$, $p+q \geq 2$,
\item $\mathfrak{so}(p,q)$, avec $p > q > 0$, $p+q$ impair, $p+q \geq
  5$, ou $p > q > 0$, $p+q$ pair, $p+q \geq
  8$,
\item $\mathfrak{sp}(p,q)$, avec $p \geq q > 0$, $p+q \geq 3$,
\item $\mathfrak{sp}(n,\mathbb{R})$, avec $n \geq 3$,
\item $\mathfrak{so}^*(2n)$, avec $n \geq 4$,
\end{itemize}

\item[(d)] Les $12$ alg{\`e}bres de Lie simples exceptionnelles non
  complexes, non compactes $EI$, $EII$, $EIII$, $EIV$, $EV$, $EVI$,
  $EVII$, $EVIII$, $EIX$, $FI$, $FII$ et $G$.\\
 
\end{description}

Si la condition 
  $({\tt *})$ est
  satisfaite, alors le cardinal de $ \mathcal{K}_{{\rm
      comp}}(\widehat{\Pi})$ est donn{\'e} par la formule suivante\,:
\begin{eqnarray}\label{Kcomp}
\# \mathcal{K}_{{\rm comp}}(\widehat{\Pi})=k_{\mathfrak{g}} -
k_{\mathfrak{m}} + {\rm rg} \; \mathfrak{g} - {\rm rg} \; \mathfrak{k} - \dim
\widehat{\mathfrak{a}},
\end{eqnarray}
d'apr{\`e}s la proposition \ref{calcul_k}, car
  $k_{\mathfrak{m}}=\# \mathcal{K}(\widehat{\Pi}^{'})$. La relation
  (\ref{Kcomp}) permettra, lorsque la condition {\bf (A)} est remplie, de
  voir si $\mathfrak{b}$ satisfait de plus aux conditions {\bf (B)} et {\bf (B)'}.

\subsection*{ (a) Cas des formes r{\'e}elles compactes d'alg{\`e}bres de Lie simples
  complexes} 

Si $\mathfrak{g}_0$ est la forme r{\'e}elle compacte d'une alg{\`e}bre de Lie
  simple complexe, alors $\mathfrak{k}_0=\mathfrak{g}_0$ et la sous-alg{\`e}bre $\mathfrak{b}$ est
  nulle.\\ 

\subsection*{ (b) Cas des alg{\`e}bres de Lie r{\'e}elles sous-jacentes {\`a} une alg{\`e}bre
de Lie simple complexe} 

Si $\mathfrak{g}_0$ est l'alg{\`e}bre de Lie r{\'e}elle sous-jacente {\`a} une alg{\`e}bre
de Lie simple complexe, $\mathfrak{g}_0$ est de la forme
 $\mathfrak{s}^{\mathbb{R}}$ avec $\mathfrak{s}$ simple complexe. Soit $\mathfrak{u}_0$
 une forme r{\'e}elle compacte de $\mathfrak{s}$. La d{\'e}composition de Cartan de $\mathfrak{g}_{0}$
s'{\'e}crit $\mathfrak{g}_{0}=\mathfrak{u}_0 \oplus i \mathfrak{u}_0$ et $\theta$ est la conjugaison complexe relative {\`a} $\mathfrak{u}_0$. 
Ici, l'alg{\`e}bre $\mathfrak{k}_0$ est la forme r{\'e}elle compacte $\mathfrak{u}_0$, d'o{\`u}
${\rm rg} \ \mathfrak{k}_0={\rm rg} \ \mathfrak{u}_0={\rm rg} \
\mathfrak{s}$.  On a de plus, ${\rm rg} \ \mathfrak{g}_0={\rm rg \ } \mathfrak{g}=2 \ {\rm rg}
\ \mathfrak{s}$. Soit $\mathfrak{c}_0$ une sous-alg{\`e}bre de Cartan de
$\mathfrak{u}_0$. Alors on a\,: $\widehat{\mathfrak{a}}_0=i \mathfrak{c}_0$ et
$\mathfrak{m}_0=\mathfrak{c}_0$. Notons que la sous-alg{\`e}bre de Cartan
$\widehat{\mathfrak{h}}_0=\widehat{\mathfrak{a}}_0 \oplus i \widehat{\mathfrak{a}}_0$ est {\`a} la fois maximalement compact et maximalement
non-compact. Puisque $\mathfrak{m}_0$ est
une sous-alg{\`e}bre ab{\'e}lienne, la condition
$({\tt *})$ est
remplie. La relation (\ref{Kcomp}) donne alors\,:
$\# \mathcal{K}_{{\rm comp}}(\widehat{\Pi})= 2 k_{\mathfrak{s}} - 0 +
2 {\rm rg} \
\mathfrak{s} - {\rm rg} \
\mathfrak{s} - {\rm rg} \
\mathfrak{s}= 2 k_{\mathfrak{s}} \not=0$. D'apr{\`e}s la proposition
\ref{caract} et le th{\'e}or{\`e}me \ref{formule_indice}, la sous-alg{\`e}bre $\mathfrak{b}$ ne
poss{\`e}de pas de forme stable et son indice est {\'e}gal {\`a} ${\rm rg} \
\mathfrak{s}$. Ainsi $\mathfrak{g}_0$ satisfait {\`a} la condition {\bf (A)},
mais ne satisfait pas aux conditions {\bf (B)}, {\bf (B)'} et {\bf (C)}.\\

\subsection{(c) et (d) Cas des alg{\`e}bres de Lie simples r{\'e}elles, non complexes, non
  compactes} 

On suppose que l'alg{\`e}bre de Lie r{\'e}elle $\mathfrak{g}_0$ est l'une
  des alg{\`e}bres de la liste {\bf (c)} ou {\bf (d)}. 

On utilise ici les donn{\'e}es de \cite{Knapp}, Appendice
C. Soit $K$ le sous-groupe connexe de $G$ d'alg{\`e}bre de Lie
$\mathfrak{k}$. Lorsque le quotient $G/K$ est hermitien sym{\'e}trique, il
r{\'e}sulte par exemple de \cite{Rossi} que le
dual de $\mathfrak{b}$ poss{\`e}de une $\mathfrak{b}$-orbite ouverte. La
condition {\bf (C)} est alors satisfaite. On utilise en outre les notations et les donn{\'e}es des tables \ref{classic} et \ref{except} pour
{\'e}tudier la condition {\bf (A)}. 

On regroupe
  dans la table \ref{tableau} les donn{\'e}es n{\'e}cessaires concernant
  $\mathfrak{g}_0$ qui permettent d'{\'e}tudier les conditions {\bf (A)}, {\bf (B)}, {\bf (B)'} et
  {\bf (C)}. En particulier, ces donn{\'e}es permettent de calculer le cardinal de $\mathcal{K}_{{\rm
  comp}}(\widehat{\Pi})$ selon la formule (\ref{Kcomp}),
lorsque la condition
$({\tt *})$ est
remplie.\\

\begin{itemize} 
\item Pour $\mathfrak{sl}(n,\mathbb{R})$, avec $n \geq 2$, $\mathfrak{su}(p,p)$, avec $p \geq 1$, $\mathfrak{so}(p,p)$, $p \geq
  1$, $\mathfrak{so}(2p,2p+1)$, avec $p \geq 1 $, $\mathfrak{so}(2p,2p-1)$,  $\mathfrak{so}(2p-1,2p+1)$, $p \geq
  1$, $\mathfrak{sp}(n,\mathbb{R})$, avec $n \geq 1$, $EI$, $EII$, $EV$,
$EVIII$, $FI$ et $G$, la sous-alg{\`e}bre $\mathfrak{m}$ est
ab{\'e}lienne donc la condition $({\tt *})$ est remplie. On {\'e}tudie alors les conditions {\bf (B)} et {\bf (B)'} gr{\^a}ce
{\`a} la formule (\ref{Kcomp}). La condition {\bf (C)} est en outre
satisfaite d{\`e}s que ${\rm rg \ } \mathfrak{g}={\rm rg \ } \mathfrak{k}$.\\

\item $\mathfrak{sl}(n,\mathbb{H})$, avec $n \geq 2$\,: ici
  $\mathfrak{m} \simeq \mathfrak{su}(2)^n$. L'ensemble $\Pi^{'}$ a
  $n \geq 2$ composantes connexes. Il r{\'e}sulte alors de la table
  \ref{classic} que la condition $({\tt *})$ n'est pas 
remplie.\\

\item $\mathfrak{su}(p,q)$, avec $1\leq p \leq q$\,: puisque $G/K$ est
  Hermitien sym{\'e}trique, la condition {\bf (C)} est satisfaite. \\

\item $\mathfrak{so}(2p,2q+1)$, avec $1 \leq p \leq q$\,: $\mathfrak{m}
  \simeq \mathfrak{so}(2q-2p+1)$. Si $q > p$, on a 
  $\Pi^{'}=\{\beta_{2p+1},\ldots,\beta_{p+q}\}$. D'apr{\`e}s la table
  \ref{classic}, la condition $({\tt *})$ est remplie. La condition {\bf
    (C)} est de plus satisfaite car ${\rm rk \ } \mathfrak{g}={\rm rk
    \ } \mathfrak{k}$.\\
\\
Si $p=q$, alors $\mathfrak{m}$ est
  ab{\'e}lienne; la condition $({\tt *})$ est remplie et {\bf (C)} est
  alors satisfaite.\\

\item $\mathfrak{so}(2p,2q+1)$, avec $p > q \geq 0$\,: $\mathfrak{m}
  \simeq \mathfrak{so}(2p-2q-1)$. Si $p-q-1 \geq 2$, alors on a 
$\Pi^{'}=\{\beta_{2p+2},\ldots,\beta_{p+q}\}$, donc la condition $({\tt *})$ n'est pas 
remplie.\\ 
\\
Si $p=q+1$, alors $\mathfrak{m}$ est
  ab{\'e}lienne; la condition $({\tt *})$ est remplie et {\bf (C)} est
  alors satisfaite.\\
\\
Si $p=q+2$, alors  $\mathfrak{m}
  \simeq \mathfrak{so}(3)$ et $\Pi'=\{\beta_l\}$. La condition $({\tt *})$ n'est pas 
remplie d'apr{\`e}s la table \ref{classic}, car $l=2q+2$ est pair.\\

\item $\mathfrak{sp}(p,q)$, avec $1 \leq p \leq q$\,: $\mathfrak{m}
  \simeq \mathfrak{su}(2)^p \oplus \mathfrak{sp}(q-p)$. D{\`e}s que $p>1$ ou $q-p \geq 1$, alors $\Pi'$ n'est pas connexe. Il r{\'e}sulte alors de
  la table \ref{classic} que la condition $({\tt *})$ n'est pas
  remplie.\\
\\
Si $p=q=1$, alors $\mathfrak{g}$ est de type $B_2$ et $\Pi'=\{\beta_2\}$; la condition $({\tt *})$ n'est donc pas
  remplie.\\

\item $\mathfrak{so}(2p+1,2q+1)$, avec $0 \leq p \leq  q$, sauf
  $\mathfrak{so}(1,1)$ et $\mathfrak{so}(1,3)$\,: $\mathfrak{m}
  \simeq \mathfrak{so}(2q-2p)$. Si $q-p \geq 4$, alors $\Pi^{'}=\{\beta_{2p},\ldots,\beta_{p+q+1}\}$, donc la condition $({\tt *})$ n'est pas
remplie.\\
\\
Si $p=q$, alors $\mathfrak{m}$ est ab{\'e}lienne; la condition
$({\tt *})$ est donc remplie et on a ${\rm rk \ } \mathfrak{g}- {\rm rk \ }
\mathfrak{k}=1 >0$, donc {\bf (C)} n'a pas lieu. {\'E}tudions la
condition {\bf (B)}\,: on a, d'apr{\`e}s la relation (\ref{Kcomp})\,:
\begin{eqnarray*}
 \# \mathcal{K}_{{\rm comp}}(\widehat{\Pi}) & =& k_{\mathfrak{g}} -
k_{\mathfrak{m}} + {\rm rg} \; \mathfrak{g} - {\rm rg} \; \mathfrak{k} - \dim
\widehat{\mathfrak{a}},\\
& =& 2p - 0 + 2p+1 - 2p - (2p+1)= 0, 
\end{eqnarray*}
donc {\bf (B)} a lieu.\\
\\
Si $q=p+1$, alors $\mathfrak{m}$ est ab{\'e}lienne; la condition
$({\tt *})$ est donc remplie et on a ${\rm rk \ } \mathfrak{g}- {\rm rk \ }
\mathfrak{k}=1 >0$, donc {\bf (C)} n'a pas lieu. {\'E}tudions la
condition {\bf (B)}\,: on a, d'apr{\`e}s la relation (\ref{Kcomp})\,:
\begin{eqnarray*}
 \# \mathcal{K}_{{\rm comp}}(\widehat{\Pi}) & =& k_{\mathfrak{g}} -
k_{\mathfrak{m}} + {\rm rg} \; \mathfrak{g} - {\rm rg} \; \mathfrak{k} - \dim
\widehat{\mathfrak{a}},\\
& =& 2p+2 - 0 + 2p+2 - (2p+1) - (2p+1)=2 \not= 0, 
\end{eqnarray*}
donc {\bf (B)} n'a pas lieu non plus.\\
\\
Si $q=p+2$, alors $\mathfrak{m} \simeq \mathfrak{so}(4)$ est de type
$A_1 \times A_1$. On a $\Pi'=\{ \beta_{l-1},\beta_l\}$, donc  la condition $({\tt *})$ n'est pas
  remplie, d'apr{\`e}s la table \ref{classic}, car $l=2p+3$ est impair.\\
\\
Si $q=p+3$, alors $\mathfrak{m} \simeq \mathfrak{so}(6)$ est de type
$A_3$. On a $\Pi'=\{\beta_{l-2}, \beta_{l-1},\beta_l\}$, donc  la condition $({\tt *})$ n'est pas
  remplie, d'apr{\`e}s la table \ref{classic}, car $l=2p+4$ est pair.\\

\item $\mathfrak{so}(2p,2q)$, avec $1 \leq p \leq q$, sauf $\mathfrak{so}(2,2)$\,:  $\mathfrak{m}
  \simeq \mathfrak{so}(2q-2p)$.  Si $q-p \geq 4$, alors $\Pi^{'}=\{\beta_{2p+1},\ldots,\beta_{p+q}\}$, donc la condition $({\tt *})$ est 
remplie.\\
\\
Si $q=p$ ou si $q=p+1$, alors $\mathfrak{m}$ est ab{\'e}lienne et la condition $({\tt *})$ est 
remplie.\\
\\
Si $q=p+2$, $\mathfrak{m}$ est de type $A_1 \times A_1$. On a
$\Pi'=\{\beta_{l-1},\beta_l\}$. Alors $({\tt *})$ est 
remplie d'apr{\`e}s la table \ref{classic}, car $l=2p+2$ est pair.\\
\\
Si $q=p+3$, $\mathfrak{m}$ est de type $A_3$. On a
$\Pi'=\{\beta_{l-2},\beta_{l-1},\beta_l\}$. Alors $({\tt *})$ est 
remplie d'apr{\`e}s la table \ref{classic}, car $l=2p+3$ est impair.\\

Dans tous les cas la condition $({\tt *})$ est 
remplie et la condition {\bf (C)} est de plus satisfaite car ${\rm rk \ } \mathfrak{g}={\rm rk
    \ } \mathfrak{k}$.\\

\item $\mathfrak{so}^*(2n)$\,: puisque $G/K$ est
  Hermitien sym{\'e}trique, la condition {\bf (C)} est satisfaite. \\

\item $EIII$\,: puisque $G/K$ est
  Hermitien sym{\'e}trique, la condition {\bf (C)} est satisfaite. \\

\item $EIV$\,:  $\mathfrak{m} \simeq \mathfrak{so}(8)$. On a
   $\Pi^{'}=\{\beta_{2},\beta_{3},\beta_4,\beta_5\}$
  donc la condition $({\tt *})$ n'est pas remplie, d'apr{\`e}s la table
  \ref{except}.\\

\item $EVI$\,: $\Pi^{'}$ est form{\'e} de
  3 racines deux {\`a} deux fortement orthogonales. D'apr{\`e}s la
  proposition \ref{m_quasi}, la sous-alg{\`e}bre $\mathfrak{m} \oplus
  \widehat{\mathfrak{a}} \oplus \mathfrak{n}$ est quasi-r{\'e}ductive. On
  utilise alors la table
  \ref{colored_root} pour en d{\'e}duire que $\Pi^{'} \subset
  \{\beta_{2},\beta_{3},\beta_5,\beta_7\}$. La
  condition $({\tt *})$ est alors remplie,  d'apr{\`e}s la table \ref{except}.\\

\item $EVII$, $EIX$\,:  $\mathfrak{m} \simeq \mathfrak{so}(8)$. On a
   $\Pi^{'}=\{\beta_{2},\beta_{3},\beta_4,\beta_5\}$. Il r{\'e}sulte alors
   de la table \ref{except} que la condition $({\tt *})$ est 
remplie.\\

\item $FII$\,: $\mathfrak{m} \simeq \mathfrak{so}(7)$.  On a
   $\Pi^{'}=\{\beta_{1},\beta_{2},\beta_3\}$,
  donc la condition $({\tt *})$ n'est pas remplie d'apr{\`e}s la table \ref{except}.\\

\end{itemize}

Dans la table \ref{tableau}, on donne, pour chaque alg{\`e}bre des listes {\bf (c)} et {\bf (d)}, le type du complexifi{\'e} $\mathfrak{g}$ et
  son rang, la sous-alg{\`e}bre $\mathfrak{k}_0$ et
  son rang, la dimension de $\widehat{\mathfrak{a}}$, $k_\mathfrak{g}$, la sous-alg{\`e}bre
  $\mathfrak{m}_0$ et $k_{\mathfrak{m}}$. Dans la derni{\`e}re colonne, on
indique la propri{\'e}t{\'e} la plus forte satisfaite par $\mathfrak{b}$. La
mention \guillemotleft rien\guillemotright \  signifie que l'indice n'est pas additif dans la d{\'e}composition
$\mathfrak{g}=\mathfrak{k} \oplus \mathfrak{b}$.\\ 

Les cas de $EI$ et
  $\mathfrak{sl}(n,\mathbb{R})$ montrent que la condition {\bf (B)} n'implique
  pas {\bf (C)}. Le cas de $\mathfrak{so}(2p-1,2p+1)$,
  $p \geq 1$ et des alg{\`e}bres de Lie simples complexes montrent que
  la condition {\bf (A)} n'implique pas {\bf (B)}.

\footnotesize
{\tiny
\begin{sidewaystable}[ht]
\begin{center}
\hspace{-5cm}
\begin{tabular}{|l|l|l|l|l|l|l|l|l|}
\hline
\textbf{$\mathfrak{g}_0$\,:} & \textbf{${\rm rg} \; \mathfrak{g}$\,:} &
\textbf{$\mathfrak{k}_0$\,:} &
\textbf{${\rm rg} \; \mathfrak{k}$\,:} &  
\textbf{$\dim \widehat{\mathfrak{a}}$\,:} & 
\textbf{$k_{\mathfrak{g}}$\,:} &
\textbf{$\mathfrak{m}_0$\,:} &
\textbf{$k_{\mathfrak{m}}$\,:} & 
\textbf{$\mathfrak{b}$}\\
\hline
$\mathfrak{sl}(n,\mathbb{R}), n \geq 2$ & $n-1$ &
$\mathfrak{so}(n)$ & $\left[ \frac{n}{2} \right]$ & 
$n-1$ & $\left[ \frac{n}{2} \right]$ &
 $0$ & $0$ &  {\bf (B)} si $n >2$,\\
 &  &
 &  & 
 &  &
  &  &  {\bf (C)} si $n=2$\\
&&&&&&&&\\
$\mathfrak{sl}(n,\mathbb{H}), n \geq 2$ & $2n-1$ & $\mathfrak{sp}(n)$ &
$n$ & $n-1$ & $n$ & $\mathfrak{su}(2)^{n}$ & $n$ & rien\\ 
$\mathfrak{su}(p,q), 1 \leq p < q$ 
& $p+q-1$ &
$\mathfrak{s}(\mathfrak{u}(p) \oplus \mathfrak{u}(q))$ &
$p+q-1$  & $p$ & $\left[ \frac{p+q}{2} \right]$ & $\mathbb{R} \oplus
\mathfrak{su}(q-p)$ & $\left[
\frac{q-p}{2} \right]$  & {\bf (C)}
\\
$\mathfrak{so}(2p,2q+1),$
& $p+q$ &
$\mathfrak{so}(2p) \oplus \mathfrak{so}(2q+1)$ &
$p+q$ & $2p$ & $p+q$ & $\mathfrak{so}(2q-2p+1)$ & $q-p$ & {\bf (C)}
\\
$1 \leq p \leq q$ &&&&&&&&\\
$\mathfrak{so}(2p,2q+1)$,
& $p+q$ &
$\mathfrak{so}(2p) \oplus \mathfrak{so}(2q+1)$ &
$p+q$ & $2q+1$ & $p+q$ & $\mathfrak{so}(2p-2q-1)$ & $p-q-1$ & rien si
$p \not= q+1$, \\
$0 \leq q < p$, &&&&&&&& {\bf (C)} si $p=q+1$.\\
&&&&&&&&\\
$\mathfrak{sp}(p,q), 1 \leq p \leq q$ 
&  $p+q$ &
$\mathfrak{sp}(p) \oplus \mathfrak{sp}(q)$ &
$p+q$ & $p$ & $p+q$ & $\mathfrak{su}(2)^{p} \oplus \mathfrak{sp}(q-p)$
& $q$ & rien \\
$\mathfrak{sp}(n,\mathbb{R}), n \geq 1$ & $n$ &
$\mathfrak{u}(n)$ &
$n$ & $n$ & $n$ & $0$ & $0$ & {\bf (C)}\\
&&&&&&&&\\
$\mathfrak{so}(2p+1,2q+1)$,
&  $p+q+1$
\ &
$\mathfrak{so}(2p+1) \oplus \mathfrak{so}(2q+1)$, \ &
$p+q$ & $2p+1$ & $2 \left[ \frac{p+q+1}{2} \right]$ &
$\mathfrak{so}(2q-2p)$ & $2 \left[ \frac{q-p}{2} \right]$ & rien si 
$q \geq p+2$,\\
$1 \leq p < q$, sauf  &&&&&&&& {\bf (A)} si $q=p+1$,\\
$\mathfrak{so}(1,1)$ ou $\mathfrak{so}(1,3)$ &&&&&&&& {\bf (B)} si $q=p$ . \\
&&&&&&&&\\
$\mathfrak{so}(2p,2q)$,
&  $p+q$ &
$\mathfrak{so}(2p) \oplus \mathfrak{so}(2q)$ &
$p+q$ & $2p$ & $2 \left[ \frac{p+q}{2} \right]$ &
$\mathfrak{so}(2q-2p)$ & $2 \left[ \frac{q-p}{2} \right]$ & {\bf (C)}
\\
$1 \leq p \leq  q$, sauf $\mathfrak{so}(2,2)$ &&&&&&&& \\
$\mathfrak{so}^*(2n), n \geq 3$   
&  $n$ &
$\mathfrak{u}(n)$ &
$n$ & $\left[\frac{n}{2}\right]$ & $2 \left[ \frac{n}{2} \right] $ &
$\mathfrak{su}(2)^{\left[\frac{n}{2}\right]}$  si $n$ paire et, &
$\left[\frac{n}{2}\right]$ & {\bf (C)}
\\
&   & & &  &  & $\mathfrak{su}(2)^{\left[\frac{n}{2}\right]} \oplus \mathbb{R}$ sinon & & \\
\hline
$EI$ &  $6$ & $\mathfrak{sp}(4)$ & $4$ & $6$ & $4$ & $0$ & $0$
&{\bf (B)}
\\
$EII$ &  6 & $\mathfrak{su}(6) \oplus \mathfrak{su}(2)$ & 6 & 4 &
$4$ & $\mathbb{R}^2$ & $0$ & {\bf (C)}
\\
$EIII$ &  $6$ & $\mathfrak{so}(10) \oplus \mathbb{R}$ & 6 & 2  &
$4$  & $\mathfrak{su}(4) \oplus \mathbb{R}$ &  $2$ & {\bf (C)}
\\
$EIV$ &  $6$ & $\mathfrak{f}_4$ & $4$ & $2$ & $4$ &
$\mathfrak{so}(8)$ & $4$ & rien
\\
$EV$ & $7$ & $\mathfrak{su}(8)$ & $7$ & $7$ & $7$ & $0$ & $0$ & {\bf (C)}
\\
$EVI$ & $7$ & $\mathfrak{so}(12) \oplus \mathfrak{su}(2)$ &
$7$ & $4$ & $7$  & $\mathfrak{su}(2)^3$ & $3$ & {\bf (C)}
\\
$EVII$ &  $7$ & $\mathfrak{e}_6 \oplus \mathbb{R}$ & $7$ & $3$ & $7$ &
$\mathfrak{so}(8)$ & $4$ & {\bf (C)}
\\
$EVIII$ &  $8$ & $\mathfrak{so}(16)$ & $8$ & $8$ & $8$ & $0$ &
$0$ & {\bf (C)}
\\
$EIX$ &  $8$ & $\mathfrak{e}_7 \oplus \mathfrak{su}(2)$ & $8$ &
$4$ & $8$ & 
$\mathfrak{so}(8)$ & $4$ & {\bf (C)}
\\
$FI$ &  4 & $\mathfrak{sp}(3)  \oplus \mathfrak{su}(2)$ & 4 & 4 &
$4$ & $0$ & $0$ & {\bf (C)}
\\
$FII$ &  4 & $\mathfrak{so}(9)$ & 4 & 1 &
$4$ & $\mathfrak{so}(7)$ & $3$ & rien
\\
$G$ &  2 & $\mathfrak{su}(2)  \oplus \mathfrak{su}(2)$ & 2 & 2 &
$2$ & $0$ &  $0$ & {\bf (C)}\\
\hline
\end{tabular}
\end{center}
~\\
\caption{\label{tableau} Donn{\'e}es concernant les alg{\`e}bres de
  Lie simples r{\'e}elles, non complexes, non compactes.}
\end{sidewaystable}
}

\nocite{*}
\bibliographystyle{plain}
\bibliography{biblio_rec}

\end{document}